\documentclass[smallextended]{svjour3}       
\smartqed  
\usepackage{graphicx}
\usepackage{amssymb}
\usepackage{latexsym}

\usepackage{amsmath,amstext,amsthm,amssymb,bbm,amsbsy}
\usepackage{mathtools}
\usepackage{tabularx}
\usepackage{subcaption}
\usepackage{dsfont}
\usepackage{xcolor}

\usepackage[ruled,linesnumbered]{algorithm2e}
\usepackage{adjustbox}

\usepackage[backref=page]{hyperref}
\renewcommand*{\backref}[1]{}
\renewcommand*{\backrefalt}[4]{%
    \ifcase #1 (Not cited.)%
    \or        (Cited on page~#2.)%
    \else      (Cited on pages~#2.)%
    \fi}

\DeclareMathOperator*{\argmin}{arg\,min}

\def\bx {\boldsymbol{x}}
\def\by {\boldsymbol{y}}
\def\bz {\boldsymbol{z}}
\def\bv {\boldsymbol{v}}
\def\bu {\boldsymbol{u}}
\def\bv {\boldsymbol{v}}
\def\bd {\boldsymbol{d}}
\def\bp {\boldsymbol{p}}

\def\ba {\boldsymbol{a}}
\def\bb {\boldsymbol{b}}

\def\balp {\boldsymbol{\alpha}}
\def\bbet {\boldsymbol{\beta}}
\def\R {\mathbb{R}}
\def\N {\mathbb{N}}
\def\dom {\mathrm{dom~}}

\def\co {\mathrm{co~}}

\def\unitsim {\Delta}
\def\Rn {\R^n}

\def\ind {I}

\def\opttraj {\gamma}
\def\opttrajhd {\boldsymbol\opttraj}
\def\initpos {u}
\def\initposhd {\bu}

\newcommand{\initcond}{\Phi}
\newcommand{\potentfn}{K}
\newcommand{\valuefn}{V}
\newcommand{\matP}{P}

\def\proxpty{y} 

\def\vz{\bv_0}

\def\opttrajgeneral {\opttrajhd}

\def\admmb{w}
\def\admmbb{\boldsymbol{\admmb}}
\def\Sanaly{\valuefn}
\def\gmanaly{\opttrajhd}
\def\Snum{\hat{\valuefn}}
\def\gmnum{\hat{\opttrajhd}}

\begin{document}

\title{Lax-Oleinik-type formulas and efficient algorithms for certain high-dimensional optimal control problems\thanks{Research supported by 
 DOE-MMICS SEA-CROGS DE-SC0023191 and AFOSR MURI FA9550-20-1-0358. P.C. is supported by the SMART Scholarship, which is funded by USD/R\&E (The Under Secretary of Defense-Research and Engineering), National Defense Education Program (NDEP) / BA-1, Basic Research.
Authors' names are given in last/family name alphabetical order.}
}

\titlerunning{}        

\author{Paula Chen \and J\'er\^ome Darbon$^*$  
        \and Tingwei Meng
}

\institute{Paula Chen \at
        Division of Applied Mathematics, Brown University, Providence, RI, USA  \\
              \email{paula\_chen@brown.edu}           
              \and
        J\'er\^ome Darbon$^*$ (Corresponding author) \at
              Division of Applied Mathematics, Brown University, Providence, RI, USA \\
              \email{jerome\_darbon@brown.edu}           
           \and
        Tingwei Meng
        \at
      Department of Mathematics, UCLA, Los Angeles, CA, USA \\
      \email{tingwei@math.ucla.edu}
}

\date{Received: date / Accepted: date}

\maketitle

\begin{abstract}
Two of the main challenges in optimal control are solving problems with state-dependent running costs and developing efficient numerical solvers that are computationally tractable in high dimension. In this paper, we provide analytical solutions to certain optimal control problems whose running cost depends on the state variable and with constraints on the control. We also provide Lax-Oleinik-type representation formulas for the corresponding Hamilton-Jacobi partial differential equations with state-dependent Hamiltonians. 
Additionally, we present an efficient, grid-free numerical solver based on our representation formulas, which is shown to scale linearly with the state dimension, and thus, to overcome the curse of dimensionality. 
Using existing optimization methods and the min-plus technique, we extend our numerical solvers to address more general classes of convex and nonconvex initial costs.
We demonstrate the capabilities of our numerical solvers using implementations on a central processing unit (CPU) and a field-programmable gate array (FPGA). In several cases, our FPGA implementation obtains over a 10 times speedup compared to the CPU, which demonstrates the promising performance boosts FPGAs can achieve. Our numerical results show that our solvers have the potential to serve as a building block for solving broader classes of high-dimensional optimal control problems in real-time.
\end{abstract}

\section{Introduction}
Optimal control problems are an important class of optimization problems with many applications, including robot manipulator control~\cite{lewis2004robot,Jin2018Robot,Kim2000intelligent,Lin1998optimal,Chen2017Reachability}, humanoid robot control~\cite{Khoury2013Optimal,Feng2014Optimization,kuindersma2016optimization,Fujiwara2007optimal,fallon2015architecture,denk2001synthesis}, and trajectory planning~\cite{Coupechoux2019Optimal,Rucco2018Optimal,Hofer2016Application,Delahaye2014Mathematical,Parzani2017HJB,Lee2021Hopf}.
A general continuous optimal control problem with a finite time horizon $t\in(0,+\infty)$ reads as follows:
\begin{equation}\label{eqt:intro_optctrl}
    \valuefn(\bx,t) = \min\left\{\int_0^t
    \ell(\bx(s), s, \balp(s)) ds
    + \initcond(\bx(0)) \right\},
\end{equation}
subject to the constraint that a trajectory $\bx(\cdot)\colon [0,t]\to \Rn$ satisfies the following backward ordinary differential equation (ODE):
\begin{equation*}
    \begin{dcases}
    \dot{\bx}(s) = f(\bx(s),s,\balp(s)) & s\in (0,t),\\
    \bx(t) = \bx.
    \end{dcases}
\end{equation*}
In optimal control, we refer to $\ell\colon \Rn\times [0,t]\times A\to\R$ as the running cost (where the control space $A$ is a subset of a Euclidean space), $\initcond\colon \Rn\to\R$ as the initial cost, and the objective function of the minimization problem in~\eqref{eqt:intro_optctrl} as the cost of a control $\balp\colon [0,t]\to A$ and the corresponding trajectory $\bx(\cdot)$.

Under some assumptions, the value function $\valuefn$, as defined in~\eqref{eqt:intro_optctrl}, solves the following Hamilton-Jacobi partial differential equation (HJ PDE):
\begin{equation} \label{eqt: HJ}
\begin{dcases} 
\frac{\partial \valuefn}{\partial t}(\bx,t)+H(\bx,t, \nabla_{\bx}\valuefn(\bx,t)) = 0 & \bx\in\mathbb{R}^{n}, t\in(0,+\infty),\\
\valuefn(\bx,0)=\initcond(\bx) & \bx\in\mathbb{R}^{n},
\end{dcases}
\end{equation}
where the Hamiltonian $H\colon \Rn\times [0,T]\times \Rn\to\R$ is defined using the functions $f$ and $\ell$ in the optimal control problem~\eqref{eqt:intro_optctrl} and the initial condition is given by the initial cost $\initcond$.
Moreover, the optimal control in the optimal control problem can be recovered from the spatial gradient $\nabla_{\bx}\valuefn(\bx,t)$ of the viscosity solution $\valuefn$ to the HJ PDE.
This relationship between optimal control problems and HJ PDEs is well-known (see~\cite{Bardi1997Optimal}, for instance) and demonstrates that solving the optimal control problem~\eqref{eqt:intro_optctrl} and solving the HJ PDE~\eqref{eqt: HJ} go hand-in-hand.

One of the main challenges in optimal control and the study of HJ PDEs is handling problems involving running costs and Hamiltonians that depend on the state variable $\bx$. Typically, Hopf and Lax-Oleinik formulas are used to represent the solution of HJ PDEs. While Hopf and Lax-Oleinik representation formulas are computationally tractable for solving high-dimensional optimal control problems (see~\cite{darbon2015convex,darbon2019decomposition,Darbon2016Algorithms,yegorov2017perspectives}), they only apply to state-independent Hamiltonians.
In general, when the running cost and corresponding Hamiltonian depend on the state variable $\bx$ (which occurs in many practical applications, including~\cite{Kim2000intelligent,Chen2017Reachability,fallon2015architecture,Hofer2016Application,Lin1998optimal,Coupechoux2019Optimal}), there are no known representation formulas that are computable in high dimensions.
Instead, approximation or discretization algorithms are often used to solve high-dimensional optimal control problems with state-dependent running costs. 

A popular building block used in such algorithms is the linear-quadratic regulator (LQR), see, for instance,~\cite{Li2004iterative,Sideris2005efficient,McEneaney2006maxplus,Coupechoux2019Optimal}.
LQR solves a class of optimal control problems~\eqref{eqt:intro_optctrl} where the function $f$ is linear in $(\bx,\balp)$ with coefficients depending on $t$, the running cost $\ell$ is quadratic in $(\bx,\balp)$ with coefficients depending on $t$ and satisfying some positive definiteness assumptions, and the initial cost $\initcond$ is a certain second-order polynomial. This class of optimal control problems has analytical solutions that are easy to solve numerically and computationally tractable in high dimensions, both of which contribute to LQR's popularity. 
Although there is no constraint on the control in the LQR problem, several algorithms have been proposed that use LQR as a building block to solve optimal control problems with constraints on the control, including~\cite{Chen2019Autonomous,Chen2017Constrained,Ma2020alternating,Burachik2014Duality,Jaddu2002Spectral,Park2008LQ,Cannon2006Efficient,Aipanov2014Analytical}.
However, these algorithms also employ other approximation and discretization methods, such as time discretizations, splitting methods, barrier methods, and spectral methods. Hence, these algorithms do not directly provide the solution to these optimal control problems with constraints on the control. Recently, a generalized Lax formula was proposed in~\cite{Lee2021Computationally} for general optimal control problems with state-dependent running costs and state constraints. However, the numerical algorithm proposed in~\cite{Lee2021Computationally} still requires discretization in time and thus, also does not directly solve these state-constrained optimal control problems.

In this paper, we provide the analytical solution to a class of optimal control problems with running cost quadratic in the state variable and certain constraints on the control. We present a representation formula which solves these problems exactly, without discretizations or approximations of the optimal control problem. Note that this differs from numerical algorithms in the existing literature, which only approximate the solution. For quadratic initial costs, we show that we can compute our representation formula exactly and efficiently. For more general initial costs, we demonstrate how our solver for the quadratic case can be used as a building block for proximal point-based methods, such as the Alternating Direction Method of Multipliers (ADMM), to numerically compute the representation formulas in these more general cases. Therefore, we provide both theoretical guarantees and efficient numerical methods for this class of problems. As a result, our numerical methods have the potential to complement LQR as a building block in numerical algorithms for solving certain control-constrained optimal control problems.

Another major challenge in optimal control and the study of HJ PDEs is handling high dimensions. Many practical engineering applications involve high dimensions. For example, robot manipulator control problems involve multiple joints and end effectors, each of which yields several degrees of freedom, including velocities, angles, and positions. In turn, each of these degrees of freedom results in a state variable.
As a result, robot manipulator control problems generally have high-dimensional state spaces (with dimension usually greater than five).
However, when the dimension is high (say, greater than five), standard grid-based numerical algorithms such as ENO \cite{Osher1991High}, WENO \cite{Jiang2000Weighted}, and DG \cite{Hu1999Discontinuous} are no longer feasible to apply.
This infeasibility is due to the ``curse of dimensionality"~\cite{bellman1961adaptive}, i.e., 
the complexity of such grid-based methods scales exponentially with dimension.
Hence, efficiently solving optimal control problems and HJ PDEs in high dimensions remains an important but challenging problem. 
Previously, several methods have been proposed to overcome the curse of dimensionality when solving high-dimensional optimal control problems and their associated HJ PDEs, which include, but are not limited to, 
optimization methods \cite{darbon2015convex,darbon2019decomposition,Darbon2016Algorithms,yegorov2017perspectives,darbon2021hamilton,Lee2021Computationally},
max-plus methods \cite{akian2006max,akian2008max,dower2015maxconference,Fleming2000Max,gaubert2011curse,McEneaney2006maxplus,McEneaney2007COD,mceneaney2008curse,mceneaney2009convergence},
tensor decomposition techniques \cite{dolgov2019tensor,horowitz2014linear,todorov2009efficient},
model order reduction \cite{alla2017error,kunisch2004hjb},
polynomial approximation \cite{kalise2019robust,kalise2018polynomial},
sparse grids \cite{bokanowski2013adaptive,garcke2017suboptimal,kang2017mitigating},
dynamic programming and reinforcement learning \cite{alla2019efficient,bertsekas2019reinforcement,zhou2021actor}, and neural networks \cite{bachouch2018deep,bansal2020deepreach,Djeridane2006Neural,jiang2016using,Han2018Solving,hure2018deep,hure2019some,lambrianides2019new,Niarchos2006Neural,reisinger2019rectified,royo2016recursive,Sirignano2018DGM,Li2020generating,darbon2020overcoming,Darbon2021Neural,nakamurazimmerer2021adaptive,NakamuraZimmerer2021QRnet,jin2020learning,JIN2020Sympnets,darbon2021neuralcontrol,onken2021neural}. 

In this paper, instead of considering general optimal control problems, we focus on a class of optimal control problems with particular control constraints and a running cost that is quadratic in the state variable. 
We derive analytical representation formulas and then use these representation formulas to design efficient numerical solvers for these optimal control problems and the corresponding HJ PDEs in high dimensions (e.g., in spatial dimension 16).
We demonstrate the efficiency of our numerical solvers using both a central processing unit (CPU) implementation and a field-programmable gate array (FPGA) implementation. FPGAs are arrays of programmable logic blocks and memory elements connected by reconfigurable interconnects (we refer the reader to \cite{KastnerFPGA} for a brief overview of FPGAs). Although CPU implementations are standard in scientific computing for fast computations, real-time optimal control applications often require strict constraints on power/energy consumption, computational resource availability, and/or computational speed that cannot be met by a CPU. 
FPGAs offer more flexibility than CPUs in designing implementations that meet these constraints. Thus, although our CPU results already illustrate the efficiency of our numerical solvers, our FPGA results demonstrate the potential performance boosts that FPGAs are able to achieve over CPUs, while also tracking the specific amounts of logic and memory resources consumed. As such, our numerical solvers show promise in their ability to solve certain high-dimensional optimal control applications involving state-dependent running costs and control constraints in real time.

In this paper, we present analytical solutions to certain classes of optimal control problems and their corresponding HJ PDEs, where the associated running cost and Hamiltonian are state-dependent. We also provide efficient numerical solvers, which have the potential to solve these problems in high dimension and in real time. The organization of this paper is as follows. In Section~\ref{sec: HJ_optctrl}, we present the class of optimal control problems and HJ PDEs considered in this paper as well as the analytical solutions of these problems. 
In Section~\ref{subsec: optctrlHJ_1d}, the one-dimensional problems are analyzed. 
In Section~\ref{subsec: optctrlHJ_hd}, we consider a class of separable high-dimensional problems. In Section~\ref{subsec:optctrlHJ_hd_general}, we provide a Lax-Oleinik-type representation formula for our most general high-dimensional case.
In Section~\ref{sec: ADMM}, we propose efficient numerical solvers for these problems and present some high-dimensional numerical results. More specifically, in Section~\ref{sec:quad_convex}, we present an efficient exact solver for quadratic initial costs, which will be used later as a building block for more general initial costs. In Section~\ref{sec:numerical_convex}, we present an ADMM algorithm using the building block from the quadratic case to solve certain optimal control problems with more general convex initial costs. In Section~\ref{sec: ADMM_nonconvex}, we generalize our proposed methods to certain nonconvex initial costs using min-plus techniques. In each of these subsections, we present corresponding high-dimensional numerical results using both a CPU and an FPGA implementation. In Section~\ref{sec:conclusion}, we make some concluding remarks and list some possible future directions.
Finally, some technical lemmas and computations for the proofs and the numerical methods are provided in the Appendix.

\section{Analytical solutions}\label{sec: HJ_optctrl}
In this section, we provide the analytical solutions to certain optimal control problems, where the running cost depends on the state variable. We also provide a Lax-Oleinik-type representation formula for the corresponding HJ PDEs, where the Hamiltonian depends on the state variable. 

Specifically, we consider the following problem.
Let $\vz$ be a vector in $\Rn$, $\{a_i\}_{i=1}^n$ and $\{b_i\}_{i=1}^n$ be positive scalars, 
$\matP$ be an invertible matrix with $n$ rows and $n$ columns,
and $M:= \matP\matP^T$ be a symmetric positive definite matrix 
which defines a norm $\|\bx\|_M:=\sqrt{\langle\bx,M\bx\rangle}$. We use $\langle \cdot, \cdot\rangle$ to denote the Euclidean inner product in $\Rn$ whose associated $\ell^2$-norm is denoted by $\|\cdot\|$. 
In this paper, we use bold characters to denote high-dimensional vectors in $\Rn$, and we use $x_i$ to denote the $i$-th component of a high-dimensional vector $\bx\in\Rn$.
Our goal is to solve the following optimal control problem:
\begin{equation}\label{eqt: result_optctrl2_hd_general}
\begin{split}
    \valuefn(\bx,t) = \min \Bigg\{\int_0^t \frac{1}{2}\|\bx(s) -\vz\|_M^2  ds + \initcond(\bx(0)) \colon \bx(t) = \bx,\quad \\ \matP^T\dot{\bx}(s)\in\prod_{i=1}^n [-b_i,a_i] \,\,\forall s\in(0,t)\Bigg\},
\end{split}
\end{equation}
where $\vz$, $\{a_i\}_{i=1}^n$, $\{b_i\}_{i=1}^n$, and $\matP$ satisfy the above assumptions, the time horizon is $t>0$, the terminal position is $\bx\in\Rn$, and the initial cost is a lower semi-continuous function $\initcond\colon\Rn\to\R$. In the optimal control problem, the function $\bx(\cdot)\colon [0,t]\to\Rn$ is assumed to be a Lipschitz function, and $\dot{\bx}(s)$ denotes its derivative at time $s$, which exists at $s\in(0,t)$ almost everywhere. Any trajectory $s\mapsto \bx(s)$ is called feasible if it satisfies the constraints in the problem~\eqref{eqt: result_optctrl2_hd_general}.
To avoid the ambiguity of a trajectory and a vector, we use $\bx(\cdot)$ to denote the trajectory, which is a function of the time variable, and we use $\bx$ to denote the vector.

In the literature, it is well-known that optimal control problems are highly related to HJ PDEs. Specifically, the optimal values in the optimal control problems are equal to the viscosity solutions to the corresponding HJ PDEs, and the optimal controls are related to the spatial gradient of the viscosity solutions to the HJ PDEs (see \cite{Bardi1997Optimal}).
In our case, the optimal control problem~\eqref{eqt: result_optctrl2_hd_general} corresponds to the following HJ PDE:
\begin{equation} \label{eqt: HJhd_2_general}
\begin{dcases} 
\frac{\partial \valuefn}{\partial t}(\bx,t) + \potentfn(\nabla_{\bx}\valuefn(\bx,t)) - \frac{1}{2}\|\bx-\vz\|_M^2 = 0 & \bx\in\mathbb{R}^n, t\in(0,+\infty),\\
\valuefn(\bx,0)=\initcond(\bx) & \bx\in\mathbb{R}^n,
\end{dcases}
\end{equation}
where $\potentfn\colon\Rn\to[0,+\infty)$ is a piecewise affine 1-homogeneous convex function related to the parameters $\{a_i\}_{i=1}^n$, $\{b_i\}_{i=1}^n$, and $\matP$ in the optimal control problem~\eqref{eqt: result_optctrl2_hd_general}.

In Section~\ref{subsec: optctrlHJ_1d}, we provide the 
analytical solution
for the one-dimensional case. In Section~\ref{subsec: optctrlHJ_hd}, we solve a simple separable high-dimensional case where $\vz$ is the zero vector and $\matP$ is the identity matrix. In Section~\ref{subsec:optctrlHJ_hd_general}, we present the analytical solution to the high-dimensional optimal control problem~\eqref{eqt: result_optctrl2_hd_general} and the HJ PDE~\eqref{eqt: HJhd_2_general}.

\subsection{One-dimensional case}\label{subsec: optctrlHJ_1d}

In this section, we consider the one-dimensional optimal control problem, which reads:
\begin{equation}\label{eqt: result_optctrl2_1d}
    \valuefn(x,t) = \min \left\{\int_0^t \frac{x(s)^2}{2} ds + \initcond(x(0)) \colon \dot{x}(s) \in[-b,a]\,\, \forall s\in(0,t), \,\, x(t) = x\right\},
\end{equation}
where $x\in\R$ and $t>0$ denote the terminal position and time horizon, respectively, and $a,b>0$ are positive scalars which give the restrictions on the velocity $\dot{x}(s)$.
The corresponding HJ PDE reads:
\begin{equation} \label{eqt: result_HJ2_1d}
\begin{dcases} 
\frac{\partial \valuefn}{\partial t}(x,t)+ \potentfn(\nabla_{x}\valuefn(x,t)) - \frac{x^2}{2} = 0 & x\in\mathbb{R}, t\in(0,+\infty),\\
\valuefn(x,0)=\initcond(x) & x\in\mathbb{R},
\end{dcases}
\end{equation}
where $\potentfn\colon \R\to[0,+\infty)$ is the 1-homogeneous convex function defined by:
\begin{equation}\label{eqt: result_defV}
    \potentfn(x) = \begin{dcases}
    ax & x\geq 0,\\
    -bx & x< 0,
    \end{dcases}
\end{equation}
where $a$ and $b$ are the positive parameters in~\eqref{eqt: result_optctrl2_1d}.

In what follows, we will present the analytical solutions to the optimal control problem~\eqref{eqt: result_optctrl2_1d} and the HJ PDE~\eqref{eqt: result_HJ2_1d}.
First, we start with the case when $\initcond(x) = \ind_{\{\initpos\}}(x)$ for some $\initpos\in\R$, where $\ind_{\{\initpos\}}$ denotes the indicator function of the set $\{\initpos\}$. Recall that the indicator function $\ind_C\colon\Rn\to\R\cup\{+\infty\}$ of a set $C\subseteq \Rn$ is defined by:
\begin{equation*}
    \ind_C(x) := \begin{dcases}
    0 & x\in C,\\
    +\infty & x\not\in C.
    \end{dcases}
\end{equation*}
Under this initial cost, the initial position in the optimal control problem~\eqref{eqt: result_optctrl2_1d} is fixed to be the point $\initpos$. Thus, the optimal control problem~\eqref{eqt: result_optctrl2_1d} becomes
\begin{equation}\label{eqt: optctrl_1d_indicatorcost}
    \min \left\{\int_0^t \frac{x(s)^2}{2} ds  \colon \dot{x}(s) \in [-b,a]\,\, \forall s\in(0,t), \,\, x(0)=\initpos, \,x(t) = x\right\}.
\end{equation}
We denote by $\valuefn(x,t; \initpos,a,b)$ the optimal value in the minimization problem~\eqref{eqt: optctrl_1d_indicatorcost}, which is a function of $(x,t)\in\R\times [0,+\infty)$ with parameters $\initpos\in\R$, $a>0$, and $b>0$. The variables $x$ and $t$ are the terminal position and time horizon, respectively, the parameter $\initpos$ denotes the parameter in the initial cost $\ind_{\{\initpos\}}$, and the parameters $a$ and $b$ are the positive parameters in~\eqref{eqt: optctrl_1d_indicatorcost}. 
The optimal trajectory is denoted by $s\mapsto \opttraj(s;x,t,\initpos,a,b)$. In the notation of the optimal trajectory, there are five parameters: $x$, $t$, $\initpos$, $a$, and $b$, which have the same meaning as the corresponding variables or parameters in the notation of the optimal value $\valuefn(x,t; \initpos,a,b)$. 

\begin{figure}[htbp]
    \centering
    \includegraphics[width=0.7\textwidth]{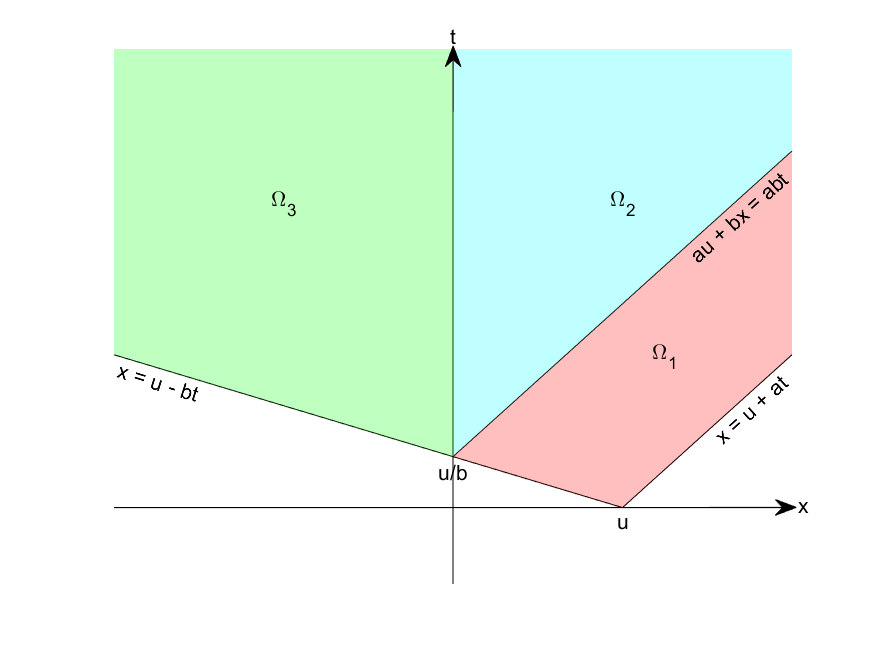}
    \caption{An illustration of a two-dimensional slice of the three sets $\Omega_1,\Omega_2,\Omega_3$ on the $xt$-plane.}
    \label{fig: plot_of_domain}
\end{figure}

Now, we define the functions $\valuefn$ and $\opttraj$. In Lemma~\ref{lem: optctrl2_traj_x0}, we will prove that the following definitions are indeed the optimal value and optimal trajectory in the corresponding optimal control problem~\eqref{eqt: optctrl_1d_indicatorcost}.
When $\initpos \geq 0$, we define the function $\R\times [0,+\infty)\ni(x,t)\mapsto \valuefn(x,t; \initpos,a,b)\in \R\cup\{+\infty\}$ by
\begin{equation}\label{eqt: result_S2_1d}
    \valuefn(x,t; \initpos,a,b) :=\begin{dcases}
    \frac{\initpos^3}{6b} + \frac{x^3}{6a} - \left(\frac{1}{6a} + \frac{1}{6b}\right)\left(\frac{a\initpos + bx - abt}{a+b}\right)^3 & \text{if } (x,t,\initpos)\in \Omega_1,\\
    \frac{\initpos^3}{6b} + \frac{x^3}{6a} & \text{if } (x,t,\initpos)\in \Omega_2,\\ 
    \frac{\initpos^3}{6b} - \frac{x^3}{6b} & \text{if } (x,t,\initpos)\in \Omega_3,\\
    +\infty & \text{otherwise},
    \end{dcases}
\end{equation}
where the three sets $\Omega_i \subset \R\times [0,+\infty)\times [0,+\infty)$ with $i=1,2,3$ are defined by:
\begin{equation}\label{eqt:def_dom_omegai}
\begin{split}
    \Omega_1 &:= \left\{(x,t,\initpos)\colon 0\leq t<\frac{\initpos}{b},\,\, \initpos-bt\leq x\leq \initpos+at\right\} \\
    &\quad\quad \quad \quad \bigcup \left\{(x,t,\initpos)\colon t\geq \frac{\initpos}{b}, \,\, at-\frac{a\initpos}{b}\leq x\leq \initpos + at\right\},\\
    \Omega_2 &:= \left\{(x,t,\initpos)\in\R\times [0,+\infty)\times [0,+\infty)\colon t\geq \frac{\initpos}{b}, \,\, 0\leq x< at-\frac{a\initpos}{b}\right\},\\
    \Omega_3 &:= \left\{(x,t,\initpos)\in\R\times [0,+\infty)\times [0,+\infty)\colon t\geq \frac{\initpos}{b}, \,\, \initpos - bt\leq x< 0\right\}.
\end{split}
\end{equation}
An illustration of a two-dimensional slice of these three sets for a fixed $\initpos > 0$ is shown in Figure~\ref{fig: plot_of_domain}. After some computation, we conclude that the domain of the function $(x,t)\mapsto \valuefn(x,t;\initpos,a,b)$ is
\begin{equation}\label{eqt:domS_1d_upos}
    \dom \left((x,t)\mapsto \valuefn(x,t;\initpos,a,b)\right) = \left\{(x,t)\in \Rn\times [0,+\infty)\colon \initpos-bt\leq x\leq \initpos+at\right\}.
\end{equation}

To define the trajectory $s\mapsto \opttraj(s;x,t,\initpos,a,b)$, we consider three cases, which correspond to the first three lines in~\eqref{eqt: result_S2_1d}, respectively. We define the trajectory
$[0,t]\ni s\mapsto \opttraj(s;x,t,\initpos,a,b)\in\R$ in these three cases as follows:
\begin{enumerate}
    \item When $(x,t,\initpos)\in\Omega_1$ holds, which corresponds to the first line in~\eqref{eqt: result_S2_1d}, we define $\opttraj(s;x,t,\initpos,a,b)$ by:
    \begin{equation}\label{eqt: optctrl2_defx_1}
        \opttraj(s;x,t,\initpos,a,b) := \begin{dcases}
        \initpos - bs & s\in\left[0, \frac{-x+\initpos + at}{a+b}\right),\\
        a(s-t) +x & s\in\left[\frac{-x+\initpos + at}{a+b}, t\right].
        \end{dcases}
    \end{equation}
    The trajectory $\opttraj(s;x,t,\initpos,a,b)$ is non-negative for all $s\in[0,t]$. The velocity $\frac{d}{ds}\opttraj(s;x,t,\initpos,a,b)$ is $-b$ in the first line and $a$ in the second line. In other words, the controlled object goes left in the first time period, and it goes right in the second time period.
    
    \item 
    When $(x,t,\initpos)\in\Omega_2$ holds, which corresponds to the second line in~\eqref{eqt: result_S2_1d}, we define $\opttraj(s;x,t,\initpos,a,b)$ by: 
    \begin{equation}\label{eqt: optctrl2_defx_2}
        \opttraj(s;x,t,\initpos,a,b) := \begin{dcases}
        \initpos - bs & s\in\left[0, \frac{\initpos}{b}\right),\\
        0 & s\in \left[\frac{\initpos}{b}, t - \frac{x}{a}\right),\\
        a\left(s-t\right) + x & s\in\left[t - \frac{x}{a}, t\right].
        \end{dcases}
    \end{equation}
    The trajectory $\opttraj(s;x,t,\initpos,a,b)$ is  non-negative for all $s\in[0,t]$. The velocity $\frac{d}{ds}\opttraj(s;x,t,\initpos,a,b)$ is $-b$ in the first line, $0$ in the second line, and $a$ in the third line.
    In other words, the controlled object first goes to the left, then stays at zero, and then goes to the right. 
    
    \item When $(x,t,\initpos)\in\Omega_3$ holds, which corresponds to the third line in~\eqref{eqt: result_S2_1d}, we define $\opttraj(s;x,t,\initpos,a,b)$ by:
    \begin{equation}\label{eqt: optctrl2_defx_3}
        \opttraj(s;x,t,\initpos,a,b) := \begin{dcases}
        \initpos - bs & s\in\left[0, \frac{\initpos}{b}\right),\\
        0 & s\in \left[\frac{\initpos}{b}, t - \frac{|x|}{b}\right),\\
        -b\left(s-t\right) +x & s\in\left[t - \frac{|x|}{b}, t\right].
        \end{dcases}
    \end{equation}
    In the first line, the trajectory is positive, and the velocity is $-b$. In the second line, the trajectory and the velocity are both zero. In the third line, the trajectory is negative, and the velocity is $-b$. 
    In other words, the controlled object first goes to the left, then stays at zero, and then goes to the left again. 
\end{enumerate}

So far, we have presented the representation formulas for the optimal values and the optimal trajectories of the one-dimensional optimal control problem~\eqref{eqt: optctrl_1d_indicatorcost}, where the initial position $\initpos$ is non-negative. Note that the domain of the value function $\valuefn$ is given in~\eqref{eqt:domS_1d_upos}, and the optimal trajectory $\opttraj$ is well-defined if and only if $(x,t)$ is in the domain of $\valuefn$.

Now, we consider the case when $\initpos$ is negative.
In this case, we define the optimal values and the optimal trajectories by symmetry.
To be specific, for any $\initpos<0$, the function
$\R\times [0,+\infty)\ni(x,t)\mapsto \valuefn(x,t; \initpos,a,b)\in \R\cup\{+\infty\}$ is defined by:
\begin{equation} \label{eqt: result_S2_1d_negative}
    \valuefn(x,t; \initpos, a,b) := \valuefn(-x,t; -\initpos, b,a) \quad\forall x\in\R,\, t\geq 0,
\end{equation}
where the right-hand side is defined in~\eqref{eqt: result_S2_1d}. After some computation, we conclude that the domain of the function $(x,t)\mapsto \valuefn(x,t;\initpos,a,b)$ also satisfies~\eqref{eqt:domS_1d_upos} in this case.
Similarly, for any $x\in\R$, $t\geq 0$, and $\initpos<0$, the optimal trajectory $[0,t]\ni s\mapsto \opttraj(s; x,t,\initpos,a,b)\in\R$ is defined by:
\begin{equation}\label{eqt: optctrl2_defx_neg}
    \opttraj(s; x,t,\initpos,a,b) := -\opttraj(s; -x,t,-\initpos,b,a) \quad\forall \,s\in[0,t],
\end{equation}
where the right-hand side is defined in~\eqref{eqt: optctrl2_defx_1},~\eqref{eqt: optctrl2_defx_2}, and~\eqref{eqt: optctrl2_defx_3} for different cases. Note that $\opttraj(s; -x,t,-\initpos,b,a)$ on the right-hand side of~\eqref{eqt: optctrl2_defx_neg} is well-defined if and only if the corresponding optimal value $\valuefn(-x,t; -\initpos, b,a)$ is finite. Hence, $\opttraj(s; x,t,\initpos,a,b)$ on the left-hand side of~\eqref{eqt: optctrl2_defx_neg} is well-defined if and only if $(x,t)$ is in the domain of the function $(x,t)\mapsto \valuefn(x,t;\initpos,a,b)$, which equals the set in~\eqref{eqt:domS_1d_upos}.

In the following lemma, we prove that the function $\valuefn$ and the trajectory $\opttraj$, whenever they are well-defined, are indeed the unique value function and optimal trajectory for the optimal control problem~\eqref{eqt: optctrl_1d_indicatorcost}.

\begin{lemma}\label{lem: optctrl2_traj_x0}
Let $a,b,t$ be positive scalars and $x,\initpos$ be real numbers satisfying $\initpos-bt\leq x\leq \initpos+at$. Define the function $[0,t]\ni s\mapsto \opttraj(s; x,t,\initpos,a,b)\in\R$ by~\eqref{eqt: optctrl2_defx_1},~\eqref{eqt: optctrl2_defx_2},~\eqref{eqt: optctrl2_defx_3}, and~\eqref{eqt: optctrl2_defx_neg} for different cases.
Then, $s\mapsto \opttraj(s; x,t,\initpos,a,b)$
is the unique optimal trajectory of the optimal control problem~\eqref{eqt: optctrl_1d_indicatorcost}. 
Moreover, the optimal value equals $\valuefn(x,t;\initpos,a,b)$, as defined in~\eqref{eqt: result_S2_1d} and~\eqref{eqt: result_S2_1d_negative}.
\end{lemma}
\begin{proof}
In this proof, whenever there is no ambiguity, we write $\opttraj(s)$ instead of $\opttraj(s;x,t,\initpos,a,b)$.
We assume $\initpos\geq 0$. For the case where $\initpos$ is negative, the proof is similar, so we omit it here. Let $s\mapsto x(s)$ be an arbitrary feasible trajectory. In other words, $s\mapsto x(s)$ is an absolutely continuous function satisfying $x(0)=\initpos$, $x(t)=x$, and $\dot{x}(s)\in[-b,a]$ for all $s\in(0,t)$. Our goal is to prove $\int_0^t \frac{\opttraj(s)^2}{2}ds \leq \int_0^t \frac{x(s)^2}{2}ds$. For this, it suffices to prove that $|\opttraj(s)|\leq |x(s)|$ holds for all $s\in[0,t]$.

First, assume $(x,t,\initpos)\in \Omega_1$. If $0\leq s< \frac{-x+\initpos + at}{a+b}$ holds, then we have that
\begin{equation}\label{eqt: lem31_1}
    x(s) = x(0) + \int_0^s \dot{x}(\tau)d\tau \geq \initpos -bs =\opttraj(s)\geq 0,
\end{equation}
where the first inequality holds since we have $\dot{x}(\tau)\geq -b$ for all $\tau\in[0,t]$ and $x(0)=\initpos$.
If $\frac{-x+\initpos + at}{a+b}\leq s\leq t$ holds, we obtain that
\begin{equation}\label{eqt: lem31_2}
    x(s) = x - \int_s^t \dot{x}(\tau)d\tau \geq x - a(t-s) = \opttraj(s)\geq 0,
\end{equation}
where the first inequality holds since we have $\dot{x}(\tau)\leq a$ for all $\tau\in[0,t]$.
As a result, we have shown that $|\opttraj(s)|\leq |x(s)|$ for all $s\in[0,t]$.

Next, we consider the case where $(x,t,\initpos)\in \Omega_2$. 
Note that~\eqref{eqt: lem31_1} and~\eqref{eqt: lem31_2} still hold for $s\in[0,\frac{\initpos}{b})$ and $s\in[t-\frac{x}{a},t]$, respectively. For $s\in[\frac{\initpos}{b}, t-\frac{x}{a})$, we have that
\begin{equation}\label{eqt: lem31_3}
|\opttraj(s)| = 0\leq |x(s)|.
\end{equation}
As a result, $|\opttraj(s)|\leq |x(s)|$ holds for all $s\in[0,t]$.

Now, we assume $(x,t,\initpos)\in \Omega_3$.
Note that~\eqref{eqt: lem31_1} and~\eqref{eqt: lem31_3} still hold for $s\in[0,\frac{\initpos}{b})$ and $s\in[\frac{\initpos}{b}, t-\frac{|x|}{b})$, respectively. For $s\in[t-\frac{|x|}{b},t]$, we have that
\begin{equation*}
    x(s) = x - \int_s^t \dot{x}(\tau)d\tau\leq x + b(t-s) = \opttraj(s) \leq 0,
\end{equation*}
where the first inequality holds since we have $-\dot{x}(\tau)\leq b$ for all $\tau\in[0,t]$.
As a result, we conclude that $|\opttraj(s)|\leq |x(s)|$ holds for all $s\in[0,t]$.

Therefore, $\opttraj$ is an optimal trajectory. Moreover, it is the unique optimal trajectory since the optimal control problem~\eqref{eqt: optctrl_1d_indicatorcost} is a convex optimization problem with a strictly convex objective function.
Finally, by Lemma~\ref{lem: appendix_optval2_equalS} (see Appendix \ref{sec: appendix_analyticsolns}), the cost of the trajectory $\opttraj$ equals $\valuefn(x,t;\initpos,a,b)$, and hence, $\valuefn(x,t;\initpos,a,b)$ is the optimal value in the problem~\eqref{eqt: optctrl_1d_indicatorcost}.
\end{proof}

\bigbreak
Now, we consider a general lower semi-continuous initial cost $\initcond\colon\R\to\R$.
We define the function $\R\times [0,+\infty)\ni (x,t)\mapsto \valuefn(x,t)\in \R$ as follows:
\begin{equation}\label{eqt: result_Lax2_1d}
    \valuefn(x,t) := \inf_{\initpos\in[x-at, x+bt]} \{\valuefn(x,t; \initpos,a,b) + \initcond(\initpos)\}\quad\forall x\in\R, t\geq 0,
\end{equation}
where the term $\valuefn(x,t;\initpos,a,b)$ is defined in~\eqref{eqt: result_S2_1d} and~\eqref{eqt: result_S2_1d_negative}. This is a Lax-Oleinik-type representation formula. In Propositions~\ref{prop: optctrl_1d} and~\ref{prop: HJ_1d}, we show that $\valuefn(x,t)$ is the optimal value in the optimal control problem~\eqref{eqt: result_optctrl2_1d} and that the function $\valuefn$ is the viscosity solution to the HJ PDE~\eqref{eqt: result_HJ2_1d}.

Let $u^*$ be a minimizer of the minimization problem in~\eqref{eqt: result_Lax2_1d}. Note that the minimizer $u^*$ exists since the objective function in the minimization problem in~\eqref{eqt: result_Lax2_1d} is a lower semi-continuous function with compact domain $[x-at,x+bt]$ (see~\cite[Theorem~1.9]{Rockafellar1998Variational}). However, the minimizer may be not unique. When there are multiple minimizers, let $u^*$ be one such minimizer.
We define the function $\opttraj(s; x,t)$ using $u^*$ as follows:
\begin{equation}\label{eqt: optctrl21d_traj_J}
    \opttraj(s; x,t) := \opttraj(s; x,t,u^*, a,b)\quad\forall s\in[0,t],
\end{equation}
where the term $\opttraj(s; x,t,u^*, a,b)$ on the right-hand side is defined in~\eqref{eqt: optctrl2_defx_1},~\eqref{eqt: optctrl2_defx_2}, \eqref{eqt: optctrl2_defx_3}, and~\eqref{eqt: optctrl2_defx_neg} for different cases.
Note that the minimizer $u^*$ satisfies $u^*\in [x-at, x+bt]$. In other words, $(x,t)$ is in the domain in~\eqref{eqt:domS_1d_upos} with $\initpos=u^*$, and hence the right-hand side in~\eqref{eqt: optctrl21d_traj_J} is well-defined.
We prove in Proposition~\ref{prop: optctrl_1d} that this trajectory $s\mapsto\opttraj(s; x,t)$ is indeed an optimal trajectory of the problem~\eqref{eqt: result_optctrl2_1d}.
When there are multiple minimizers in~\eqref{eqt: result_Lax2_1d}, by Proposition~\ref{prop: optctrl_1d}, for any minimizer $u^*$, the corresponding function $s\mapsto\opttraj(s; x,t)$ defined in~\eqref{eqt: optctrl21d_traj_J} with $u^*$ is one optimal trajectory.

Now, we present the main results for the one-dimensional problems~\eqref{eqt: result_optctrl2_1d} and~\eqref{eqt: result_HJ2_1d}.

\begin{proposition}\label{prop: optctrl_1d}
Let $\initcond\colon\R\to\R$ be a lower semi-continuous function and $a,b$ be some positive scalars. Then, for all $x\in\R$ and $t>0$, the function $[0,t]\ni s\mapsto \opttraj(s; x,t)\in\R$ defined in~\eqref{eqt: optctrl21d_traj_J} is an optimal trajectory for the optimal control problem~\eqref{eqt: result_optctrl2_1d}, whose optimal value equals $\valuefn(x,t)$ as defined in~\eqref{eqt: result_Lax2_1d}. Moreover, if $\initcond$ is convex, then the trajectory $s\mapsto\opttraj(s; x,t)$ is the unique optimal trajectory.
\end{proposition}
\begin{proof}
This is a corollary of Proposition~\ref{prop: optctrl2_hd}.
\end{proof}

\begin{proposition} \label{prop: HJ_1d}
Let $\initcond\colon \R\to \R$ be a continuous function and $\potentfn$ be the function defined in~\eqref{eqt: result_defV} with parameters $a,b>0$. 
Let $\valuefn$ be the function defined in~\eqref{eqt: result_Lax2_1d}. 
Then, the function $\valuefn$ is the unique viscosity solution to the HJ PDE~\eqref{eqt: result_HJ2_1d} in the solution set $C(\R\times [0,+\infty))$.
\end{proposition}
\begin{proof}
This is a corollary of Proposition~\ref{prop: HJ2_hd_2}.
\end{proof}

\subsection{Separable high-dimensional case}
\label{subsec: optctrlHJ_hd}
In this section, we consider the following high-dimensional optimal control problem with separable running cost:
\begin{equation}\label{eqt: result_optctrl2_hd}
\begin{split}
    \valuefn(\bx,t) = \min \Bigg\{\int_0^t \frac{\|\bx(s)\|^2}{2}  ds + \initcond(\bx(0)) \colon \bx(t) = \bx,\quad \\ \dot{\bx}(s)\in\prod_{i=1}^n[-b_i,a_i] \,\,\forall s\in(0,t)\Bigg\},
\end{split}
\end{equation}
where $\bx$ and $t$ are the terminal position and time horizon, respectively, $\{a_i\}$ and $\{b_i\}$ are positive scalars which provide the restrictions on the velocity $\dot{\bx}$, and the initial cost $\initcond\colon \Rn\to\R$ is a lower semi-continuous function. 
Note that we call the running cost separable since it is the sum of $n$ functions, the $j$-th of which only depends on the $j$-th component of $\bx(s)$. The constraint on the control is also separable since it can be written as $n$ constraints, the $j$-th of which only depends on the derivative of the $j$-th component of the trajectory $\bx(\cdot)$.
The corresponding HJ PDE reads:
\begin{equation} \label{eqt: result_HJ2_hd}
\begin{dcases} 
\frac{\partial \valuefn}{\partial t}(\bx,t)+ \sum_{i=1}^n \potentfn_i\left(\frac{\partial \valuefn(\bx,t)}{\partial x_i}\right) - \frac{1}{2}\|\bx\|^2 = 0 & \bx\in\mathbb{R}^n, t\in(0,+\infty),\\
\valuefn(\bx,0)=\initcond(\bx) & \bx\in\mathbb{R}^n,
\end{dcases}
\end{equation}
where each function $\potentfn_i\colon\R\to\R$ is the 1-homogeneous convex function defined in~\eqref{eqt: result_defV} with the positive constants $a = a_i$ and $b=b_i$ from~\eqref{eqt: result_optctrl2_hd} and the initial condition is given by the initial cost $\initcond$ in~\eqref{eqt: result_optctrl2_hd}.

For this problem, we denote the optimal value by $\valuefn(\bx,t)$ and the optimal trajectory by $s\mapsto \opttrajhd(s; \bx,t)$, where $\bx,t$ are parameters denoting the terminal position and the time horizon, respectively.
Define the function $\valuefn\colon \R^n\times\R\to\R$ by the following Lax-Oleinik-type representation formula:
\begin{equation}\label{eqt: result_Lax2_hd}
    \valuefn(\bx,t) := \inf_{\initposhd\in\prod_{i=1}^n [x_i-a_it, x_i+b_it]} \left\{\sum_{i=1}^n \valuefn(x_i, t; \initpos_i, a_i,b_i) + \initcond(\initposhd)\right\},
\end{equation}
for all $\bx\in\Rn, t\geq 0$ and
where each function $(x_i, t)\mapsto \valuefn(x_i, t; \initpos_i,a_i,b_i)$ on the right-hand side is the function defined by~\eqref{eqt: result_S2_1d} and~\eqref{eqt: result_S2_1d_negative}. 
Let $\initposhd^* = (\initpos_1^*, \dots, \initpos_n^*)\in\R^n$ be a minimizer in~\eqref{eqt: result_Lax2_hd}.
Note that the minimizer $\initposhd^*$ exists since the objective function in the minimization problem in~\eqref{eqt: result_Lax2_hd} is a lower semi-continuous function with compact domain $\prod_{i=1}^n[x_i-a_it,x_i+b_it]$ (see~\cite[Theorem~1.9]{Rockafellar1998Variational}). However, the minimizer may be not unique. When there are multiple minimizers, let $\initposhd^*$ be one such minimizer.
Define the trajectory $s\mapsto \opttrajhd(s; \bx,t)$ using $\initposhd^*$ as
\begin{equation}\label{eqt: optctrl2hd_traj_J}
    \opttrajhd(s; \bx,t) := \left(\opttraj(s; x_1,t,\initpos_1^*, a_1,b_1), \dots, \opttraj(s; x_n,t,\initpos_n^*, a_n,b_n)\right) \quad \forall s\in[0,t],
\end{equation}
where the function $\opttraj$ in the $i$-th component $\opttraj(s; x_i,t,\initpos_i^*, a_i,b_i)$ on the right-hand side is defined by \eqref{eqt: optctrl2_defx_1}, \eqref{eqt: optctrl2_defx_2}, \eqref{eqt: optctrl2_defx_3}, and \eqref{eqt: optctrl2_defx_neg} for different cases. Note that the $i$-th component of the minimizer $\initposhd^*$ satisfies $\initpos_i^*\in [x_i-a_it, x_i+b_it]$. 
In other words, $(x_i,t)$ is in the domain~\eqref{eqt:domS_1d_upos} with $\initpos=\initpos_i^*$, and hence, the $i$-th component on the right-hand side of~\eqref{eqt: optctrl2hd_traj_J} is well-defined for each $i\in\{1,\dots,n\}$.

Next, we present two propositions. In Proposition~\ref{prop: optctrl2_hd}, we show that the trajectory $s\mapsto \opttrajhd(s; \bx,t)$ is an optimal trajectory of the problem~\eqref{eqt: result_optctrl2_hd}, where the optimal value equals $\valuefn(\bx,t)$.
By Proposition~\ref{prop: optctrl2_hd}, if there are multiple minimizers in~\eqref{eqt: result_Lax2_hd}, for any minimizer $\initposhd^*$, the corresponding function $s\mapsto \opttrajhd(s;\bx,t)$ defined in~\eqref{eqt: optctrl2hd_traj_J} with $\initposhd^*$ is an optimal trajectory.
In Proposition~\ref{prop: HJ_hd}, we show that $\valuefn$ is the viscosity solution to the high-dimensional HJ PDE~\eqref{eqt: result_HJ2_hd}.

\begin{proposition} \label{prop: optctrl2_hd}
Let $\initcond\colon\R^n\to\R$ be a lower semi-continuous function. Let $\ba=(a_1,\dots,a_n)$ and $\bb=(b_1,\dots, b_n)$ be two vectors in $(0,+\infty)^n$. For any $\bx\in\R^n$ and $t>0$, the function $[0,t]\ni s\mapsto \opttrajhd(s; \bx,t)\in\R$ defined in~\eqref{eqt: optctrl2hd_traj_J} is an optimal trajectory of the optimal control problem~\eqref{eqt: result_optctrl2_hd}, where the optimal value equals $\valuefn(\bx,t)$, as defined in~\eqref{eqt: result_Lax2_hd}.
Moreover, if $\initcond$ is convex, then the trajectory $s\mapsto\opttraj(s; \bx,t)$ is the unique optimal trajectory.
\end{proposition}
\begin{proof}
In this proof, whenever there is no ambiguity, we use the notation $\opttrajhd(s)$ instead of $\opttrajhd(s;\bx,t)$.
Since $\initcond$ is lower semi-continuous and the domain of the optimization problem in~\eqref{eqt: result_Lax2_hd} is $\prod_{i=1}^n [x_i-a_it,x_i+b_it]$, which is compact, the minimizer $\initposhd^*$ exists (see~\cite[Theorem~1.9]{Rockafellar1998Variational}). Moreover, each component on the right-hand side of~\eqref{eqt: optctrl2hd_traj_J} is well-defined due to the constraints on $\initposhd^*$, i.e., that $\initpos_i^*\in [x_i-a_it, x_i+b_it]$ (see the discussion below~\eqref{eqt: optctrl2hd_traj_J}). Thus, $\opttrajhd$ is well-defined. Note that if $\initposhd^*$ is not unique, we will prove that any arbitrary minimizer $\initposhd^*$ of~\eqref{eqt: result_Lax2_hd} corresponds to an optimal trajectory $\opttrajhd$.

We prove by contradiction that $\opttrajhd$ is an optimal trajectory.
Assume $\opttrajhd$ is not an optimal trajectory. 
Then, there exists another trajectory $\tilde{\opttrajhd}$ satisfying 
\begin{equation}\label{eqt:prop23pf_feasible_traj}
\tilde{\opttrajhd}(t)=\bx,\quad\quad  \dot{\tilde{\opttrajhd}}(s)\in\prod_{i=1}^n[-b_i,a_i]\quad\forall\,s\in(0,t),
\end{equation}
and
\begin{equation} \label{eqt: prop35_tilde_traj}
    \initcond(\tilde{\opttrajhd}(0)) + \sum_{i=1}^n\int_0^t \frac{\tilde{\opttraj}_i(s)^2}{2} ds
    < \initcond(\opttrajhd(0)) + \sum_{i=1}^n\int_0^t \frac{\opttraj_i(s)^2}{2} ds,
\end{equation}
where $\opttraj_i$ and $\tilde{\opttraj}_i$ denote the $i$-th components of $\opttrajhd$ and $\tilde{\opttrajhd}$, respectively.
By Lemma~\ref{lem: optctrl2_traj_x0}, we have that
\begin{equation} \label{eqt: prop35_opt_i}
\int_0^t\frac{\opttraj(s; x_i,t, \tilde{\opttraj}_i(0), a_i,b_i)^2}{2} ds  \leq
\int_0^t\frac{\tilde{\opttraj}_i(s)^2}{2} ds , \quad \forall i\in\{1,\dots, n\},
\end{equation}
where $\opttraj(s; x_i,t, \tilde{\opttraj}_i(0), a_i,b_i)$ is the one-dimensional optimal trajectory of the problem~\eqref{eqt: optctrl_1d_indicatorcost} with parameters $x=x_i$, $a=a_i$, $b=b_i$, and $\initpos=\tilde{\opttraj}_i(0)$ and is defined by~\eqref{eqt: optctrl2_defx_1}, \eqref{eqt: optctrl2_defx_2}, \eqref{eqt: optctrl2_defx_3}, and~\eqref{eqt: optctrl2_defx_neg} for different cases. Note that each $\opttraj(s; x_i,t, \tilde{\opttraj}_i(0), a_i,b_i)$ is well-defined since we have $\tilde{\opttraj}_i(0)-b_it\leq x_i\leq \tilde{\opttraj}_i(0) + a_it$ by~\eqref{eqt:prop23pf_feasible_traj}.
Moreover, by Lemma~\ref{lem: appendix_optval2_equalS} (in Appendix \ref{sec: appendix_analyticsolns}), we have that
\begin{equation} \label{eqt: prop35_equalS}
\begin{split}
    &\int_0^t\frac{\opttraj(s; x_i,t, \tilde{\opttraj}_i(0), a_i,b_i)^2}{2} ds = \valuefn(x_i,t; \tilde{\opttraj}_i(0),a_i,b_i), \\
    &\int_0^t \frac{\opttraj_i(s)^2}{2} = \valuefn(x_i,t; \opttraj_i(0),a_i,b_i),
\end{split}
\end{equation}
for each $i\in\{1,\dots, n\}$. Then, by~\eqref{eqt: prop35_tilde_traj}, \eqref{eqt: prop35_opt_i}, and~\eqref{eqt: prop35_equalS}, we have
\begin{equation}\label{eqt:prop23pf_ineqt}
\begin{split}
    &\initcond(\tilde{\opttrajhd}(0)) + \sum_{i=1}^n \valuefn(x_i,t; \tilde{\opttraj}_i(0),a_i,b_i) \\
    =\,& \initcond(\tilde{\opttrajhd}(0)) + \sum_{i=1}^n\int_0^t\frac{\opttraj(s; x_i,t, \tilde{\opttraj}_i(0), a_i,b_i)^2}{2} ds\\
    \leq\,& \initcond(\tilde{\opttrajhd}(0)) + \sum_{i=1}^n\int_0^t \frac{\tilde{\opttraj}_i(s)^2}{2} ds\\
     <\,& \initcond(\opttrajhd(0)) + \sum_{i=1}^n\int_0^t \frac{\opttraj_i(s)^2}{2} ds\\
    =\,& \initcond(\opttrajhd(0)) + \sum_{i=1}^n \valuefn(x_i,t; \opttraj_i(0),a_i,b_i)\\
    =\,& \initcond(\initposhd^*) + \sum_{i=1}^n \valuefn(x_i,t; \initpos_i^*,a_i,b_i)\\
    =\,& \inf_{\initposhd\in\prod_{i=1}^n [x_i-a_it, x_i+b_it]} \left\{\sum_{i=1}^n \valuefn(x_i, t; \initpos_i, a_i,b_i) + \initcond(\initposhd)\right\} \\
    =\,& \valuefn(\bx,t),
\end{split}
\end{equation}
where the third equality holds since, by definition of $\opttrajhd$, we have that $\opttrajhd(0) = \initposhd^*$ and the fourth equality holds by definition of $\initposhd^*$. 
Note that by~\eqref{eqt:prop23pf_feasible_traj}, we have that $\tilde{\opttrajhd}(0)\in \prod_{i=1}^n[x_i-a_it, x_i+b_it]$, and hence, $\tilde{\opttrajhd}(0)$ satisfies the constraint in the minimization problem in~\eqref{eqt:prop23pf_ineqt}.
Thus, the strict inequality in~\eqref{eqt:prop23pf_ineqt} yields a contradiction, and we conclude that $\opttrajhd$ is an optimal trajectory, whose cost equals $\valuefn(\bx,t)$ by~\eqref{eqt:prop23pf_ineqt}.
Moreover, when $\initcond$ is convex, the optimal control problem~\eqref{eqt: result_optctrl2_hd} is a convex optimization problem with a strictly convex objective function, which implies the uniqueness of the optimal trajectory in this case.
\end{proof}

\begin{proposition}\label{prop: HJ_hd}
Let $\initcond\colon \R^n\to \R$ be a continuous function.
Let $\ba=(a_1,\dots,a_n)$ and $\bb=(b_1,\dots, b_n)$ be two vectors in $(0,+\infty)^n$ and $\potentfn_i\colon \R\to[0,+\infty)$ be defined by~\eqref{eqt: result_defV} with constants $a=a_i$ and $b=b_i$ for each $i\in\{1,\dots, n\}$. 
Then, the function $\valuefn$ defined in~\eqref{eqt: result_Lax2_hd} is the unique viscosity solution to the HJ PDE~\eqref{eqt: result_HJ2_hd} in the solution set $C(\Rn\times [0,+\infty))$.
\end{proposition}
\begin{proof}
This is a corollary of Proposition~\ref{prop: HJ2_hd_2}.
\end{proof}

\subsection{General high-dimensional case} \label{subsec:optctrlHJ_hd_general}
In this section, we provide the analytical solutions to the high-dimensional optimal control problem~\eqref{eqt: result_optctrl2_hd_general} and the corresponding HJ PDE~\eqref{eqt: HJhd_2_general}. Denote the optimal value by $\valuefn(\bx,t)$ and the optimal trajectory by $\opttrajgeneral(s; \bx,t)$.
We define the function $\valuefn\colon \Rn\times [0,+\infty)\to\R$ by the following Lax-Oleinik-type representation formula:
\begin{equation} \label{eqt: result_HJ2_Lax_hd_general}
    \valuefn(\bx,t) := 
    \inf_{\initposhd\in\prod_{i=1}^n [y_i-a_it, y_i+b_it]} \left\{\sum_{i=1}^n \valuefn(y_i, t; \initpos_i, a_i,b_i) + \initcond\left((\matP^T)^{-1}\initposhd+\vz\right)\right\},
\end{equation}
for all $\bx\in\Rn, t\geq 0$ and where 
the vector $\by\in\Rn$ is defined by:
\begin{equation}\label{eqt:general_hd_defy}
\by= (y_1,\dots,y_n):= \matP^{T}\bx-\matP^{T}\vz.
\end{equation}
In~\eqref{eqt: result_HJ2_Lax_hd_general},
each function $(y_i,t)\mapsto \valuefn(y_i, t; \initpos_i, a_i,b_i)$ on the right-hand side is the function defined in~\eqref{eqt: result_S2_1d} and~\eqref{eqt: result_S2_1d_negative}, and $\matP$, $\vz$, $\{a_i\}_{i=1}^n$, $\{b_i\}_{i=1}^n$ are the parameters in~\eqref{eqt: result_optctrl2_hd_general}. 

Let $\initposhd^* = (\initpos^*_1,\dots, \initpos^*_n)$ be a minimizer in the minimization problem~\eqref{eqt: result_HJ2_Lax_hd_general}. With a similar argument as in Section~\ref{subsec: optctrlHJ_hd}, we conclude that the minimizer $\initposhd^*$ exists but may be not unique. If it is not unique, let $\initposhd^*$ be one such minimizer.
Define the trajectory $[0,t]\ni s\mapsto \opttrajgeneral(s; \bx,t)\in\R^n$ by 
\begin{equation}\label{eqt: optctrl2hd_traj_general}
    \opttrajgeneral(s; \bx,t) :=
    (\matP^T)^{-1}\left(\opttraj(s; y_1, t, \initpos_1^*, a_1, b_1), \dots, \opttraj(s; y_n, t, \initpos_n^*, a_n, b_n)\right) + \vz,
\end{equation}
for all $s\in[0,t]$,
where $y_1,\dots,y_n$ are the components of the vector $\by$ defined in~\eqref{eqt:general_hd_defy} and
the $i$-th element $\opttraj(s; y_i,t,\initpos_i^*, a_i,b_i)$ on the right-hand side is the one-dimensional trajectory defined by~\eqref{eqt: optctrl2_defx_1},~\eqref{eqt: optctrl2_defx_2},~\eqref{eqt: optctrl2_defx_3}, and~\eqref{eqt: optctrl2_defx_neg} for different cases of $y_i$, $t$, and $\initpos_i^*$. Note that the $i$-th component of the minimizer $\initposhd^*$ satisfies $\initpos_i^*\in [y_i-a_it, y_i+b_it]$. 
In other words, $(y_i,t)$ is in the domain in~\eqref{eqt:domS_1d_upos} with $\initpos=\initpos_i^*$, and hence, the term $\opttraj(s; y_i, t, \initpos_i^*, a_i, b_i)$ on the right-hand side of~\eqref{eqt: optctrl2hd_traj_general} is well-defined for each $i\in\{1,\dots,n\}$.

Next, we present two propositions showing that the  functions $\valuefn$ and $\opttrajgeneral$ defined above do indeed solve the optimal control problem~\eqref{eqt: result_optctrl2_hd_general} and the corresponding HJ PDE~\eqref{eqt: HJhd_2_general}. More specifically, Proposition~\ref{prop: optctrl2_hd_general} shows that the trajectory $s\mapsto \opttrajgeneral(s; \bx,t)$ is an optimal trajectory of the problem~\eqref{eqt: result_optctrl2_hd_general}, whose optimal value equals $\valuefn(\bx,t)$. 
If there are multiple minimizers in~\eqref{eqt: result_HJ2_Lax_hd_general}, Proposition~\ref{prop: optctrl2_hd_general} shows that any minimizer $\bu^*$ defines an optimal trajectory $s\mapsto\opttrajgeneral(s;\bx,t)$ by~\eqref{eqt: optctrl2hd_traj_general}.
Proposition~\ref{prop: HJ2_hd_2} shows that the function $\valuefn$ defined in~\eqref{eqt: result_HJ2_Lax_hd_general} is the viscosity solution to the HJ PDE~\eqref{eqt: HJhd_2_general}.

\begin{proposition} \label{prop: optctrl2_hd_general}
Let $\initcond\colon \Rn\to\R$ be a lower semi-continuous function.
Let $\ba=(a_1,\dots,a_n)$ and $\bb=(b_1,\dots, b_n)$ be two vectors in $(0,+\infty)^n$, $\vz$ be a vector in $\Rn$, and $\matP$ be an invertible matrix with $n$ rows and $n$ columns.
Then, for any $\bx\in\R^n$ and $t>0$, the function $[0,t]\ni s\mapsto \opttrajgeneral(s; \bx,t)\in\R$ defined in~\eqref{eqt: optctrl2hd_traj_general} is an optimal trajectory for the optimal control problem~\eqref{eqt: result_optctrl2_hd_general}, where the matrix $M$ satisfies $M= \matP\matP^T$.
The optimal value of the problem~\eqref{eqt: result_optctrl2_hd_general} equals $\valuefn(\bx,t)$, as defined in~\eqref{eqt: result_HJ2_Lax_hd_general}.
Moreover, if $\initcond$ is convex, then the trajectory $s\mapsto\opttrajgeneral(s; \bx,t)$ is the unique optimal trajectory.
\end{proposition}
\begin{proof}
Fix $\bx\in\R^n$ and $t>0$. Define $\tilde{\initcond}\colon\R^n\to\R$ by:
\begin{equation}\label{eqt:prop25_def_newJ}
\tilde{\initcond}(\bz):= \initcond\left((\matP^T)^{-1}\bz+\vz\right) \quad \forall \, \bz\in\R^n.
\end{equation}
Let $\bx\colon [0,t]\to \R^n$ be a feasible trajectory in problem~\eqref{eqt: result_optctrl2_hd_general}, i.e., let $\bx$ satisfy the constraints in problem~\eqref{eqt: result_optctrl2_hd_general}, and let $\by(s) = \matP^T (\bx(s) - \vz)$ for all $s\in[0,t]$. Then, by straightforward computation, we have
\begin{equation*}
\begin{split}
    &\dot{\by}(s)=\matP^T \dot{\bx}(s)\in\prod_{i=1}^n [-b_i,a_i], \quad \forall s\in(0,t),\\
    &\by(t) = \matP^T(\bx(t)-\vz) = \matP^T(\bx-\vz) = \by,
\end{split}
\end{equation*}
where $\by$ is the vector defined in~\eqref{eqt:general_hd_defy}.
Therefore, $\by(\cdot)$ is a feasible trajectory for the following optimal control problem:
\begin{equation}\label{eqt: prop36_equiv_optctrl}
    \min \left\{\int_0^t \frac{1}{2}\|\by(s)\|^2  ds + \tilde{\initcond}(\by(0)) \colon \dot{\by}(s)\in\prod_{i=1}^n[-b_i,a_i] \,\,\forall s\in(0,t), \,\, \by(t) = \by\right\}.
\end{equation}
Moreover, by some computation involving the definitions of $M$ and $\tilde{\initcond}$,
the cost for the trajectory $\bx(\cdot)$ in problem~\eqref{eqt: result_optctrl2_hd_general} equals
\begin{equation*}
\begin{split}
    &\int_0^t \frac{1}{2}\|\bx(s) -\vz\|_M^2  ds + \initcond(\bx(0))\\
    =\,& \int_0^t \frac{1}{2}\|(\matP^T)^{-1}\by(s) \|_M^2  ds + \initcond\left((\matP^T)^{-1}\by(0)+\vz\right)\\
    =\,& \int_0^t \frac{1}{2}\|\by(s) \|^2  ds + \tilde{\initcond}(\by(0)),
\end{split}
\end{equation*}
which equals the cost of $\by(\cdot)$ in problem~\eqref{eqt: prop36_equiv_optctrl}.
Similarly, if $\by(\cdot)$ is a feasible trajectory in~\eqref{eqt: prop36_equiv_optctrl}, then $s\mapsto (\matP^T)^{-1}\by(s)+\vz$ is a feasible trajectory in~\eqref{eqt: result_optctrl2_hd_general} whose cost equals the cost of $\by(\cdot)$ in~\eqref{eqt: prop36_equiv_optctrl}. Therefore, the two optimal control problems~\eqref{eqt: result_optctrl2_hd_general} and~\eqref{eqt: prop36_equiv_optctrl} are equivalent to each other. By Proposition~\ref{prop: optctrl2_hd}, the trajectory 
$s\mapsto \by^*(s):=\left(\opttraj(s; y_1, t, \initpos_1^*, a_1, b_1), \dots, \opttraj(s; y_n, t, \initpos_n^*, a_n, b_n)\right)$
is an optimal trajectory for problem~\eqref{eqt: prop36_equiv_optctrl}, whose optimal value equals $\valuefn(\bx,t)$ in~\eqref{eqt: result_HJ2_Lax_hd_general}. Hence, the corresponding trajectory 
\begin{equation*}
\begin{split}
s\mapsto &(\matP^T)^{-1}\by^*(s)+\vz \\
=& (\matP^T)^{-1}\left(\opttraj(s; y_1, t, \initpos_1^*, a_1, b_1), \dots, \opttraj(s; y_n, t, \initpos_n^*, a_n, b_n)\right) + \vz\\
=& \opttrajgeneral(s; \bx,t)
\end{split}
\end{equation*}
is an optimal trajectory for problem~\eqref{eqt: result_optctrl2_hd_general}, whose optimal value also equals $\valuefn(\bx,t)$. 

Furthermore, if $\initcond$ is convex, then problem~\eqref{eqt: result_optctrl2_hd_general} is a convex optimization problem with a strictly convex objective function since the matrix $M$ is positive definite.
Thus, if $\initcond$ is convex, the minimizer is unique and the unique optimal trajectory is $s\mapsto\opttrajgeneral(s; \bx,t)$.
\end{proof}

\begin{proposition} \label{prop: HJ2_hd_2}
Let $\initcond\colon\R^n\to\R$ be a continuous function.
Let $\potentfn\colon \R^n\to [0,+\infty)$ be a piecewise affine 1-homogeneous convex function. Assume that there exist linearly independent vectors $\initposhd_1,\dots,\initposhd_n\in\Rn$ and positive scalars $a_i,b_i>0$ for each $i\in\{1,\dots, n\}$, such that the sublevel set of $\potentfn$ satisfies
\begin{equation*}
    \{\bx\in\Rn\colon \potentfn(\bx)\leq 1\} = \co \left(\bigcup_{j=1}^n \left\{\frac{1}{a_j}\initposhd_j, -\frac{1}{b_j}\initposhd_j\right\}\right),
\end{equation*}
where $\co E$ denotes the convex hull of a set $E$. Define the matrix $\matP$ to be the matrix whose columns are $\initposhd_1,\dots,\initposhd_n$, and define the matrix $M$ by $M:= \matP\matP^T$.
Then, the function $\valuefn\colon\R^n\times [0,+\infty)\to\R$  defined by~\eqref{eqt: result_HJ2_Lax_hd_general}
is the unique viscosity solution to the HJ PDE~\eqref{eqt: HJhd_2_general} in the solution set $C(\Rn\times [0,+\infty))$.
\end{proposition}
\begin{proof}
We prove this proposition using~\cite[Theorem~III.3.17]{Bardi1997Optimal}.
We write the optimal control problem~\eqref{eqt: result_optctrl2_hd_general} in the standard form in~\cite[Chapter~III]{Bardi1997Optimal}, which reads:
\begin{equation} \label{eqt: standardoptctrl}
    \inf\left\{\int_0^t
    \ell(\bx(s), \balp(s)) ds
    + \initcond(\bx(0)) \right\},
\end{equation}
subject to the constraint that 
$\balp(\cdot)\colon [0,t]\to A$ is a measurable function and $\bx(\cdot)\colon [0,t]\to\Rn$ is an absolutely continuous function satisfying the following ODE:
\begin{equation*}
    \begin{dcases}
    \dot{\bx}(s) = f(\bx(s),\balp(s)) & s\in (0,t),\\
    \bx(t) = \bx.
    \end{dcases}
\end{equation*}
In our case, the set $A$, the source term $f\colon \Rn\times A\to \Rn$, and the running cost $\ell\colon \Rn\times A\to \R$ are given by:
\begin{equation*}
\begin{split}
A = \prod_{i=1}^n[-b_i,a_i]\subset \Rn, \quad 
f(\bx,\balp) = (\matP^T)^{-1}\balp,\quad
\ell(\bx,\balp) = \frac{1}{2}\|\bx - \vz\|_M^2.
\end{split}
\end{equation*}
In~\cite[Chapter~III]{Bardi1997Optimal}, it is shown that the optimal control problem~\eqref{eqt: standardoptctrl} corresponds to the following HJ PDE:
\begin{equation}\label{eqt: standardHJPDE}
\begin{dcases} 
\frac{\partial \valuefn}{\partial t}(\bx,t) + H(\bx,\nabla_{\bx}\valuefn(\bx,t)) = 0 & \bx\in\mathbb{R}^n, t\in(0,+\infty),\\
\valuefn(\bx,0)=\initcond(\bx) & \bx\in\mathbb{R}^n,
\end{dcases}
\end{equation}
where $\initcond$ is the initial condition given by the initial cost in~\eqref{eqt: standardoptctrl} and $H$ is the Hamiltonian given by $A, f,\ell$ as follows:
\begin{equation*}
    H(\bx,\bp) = \sup_{\balp\in A} \{\langle f(\bx,\balp),\bp\rangle - \ell(\bx,\balp)\},
\end{equation*}
where $\langle\cdot,\cdot\rangle $ denotes the standard inner product in $\Rn$.
In our case, the Hamiltonian $H$ is
\begin{equation*}
\begin{split}
    H(\bx,\bp) &= \sup_{\balp\in A} \{\langle f(\bx,\balp),\bp\rangle - \ell(\bx,\balp)\}\\
    &= \sup_{\balp\in \prod_{i=1}^n[-b_i,a_i]} \left\{\langle (\matP^T)^{-1}\balp, \bp\rangle -  \frac{1}{2}\|\bx - \vz\|_M^2\right\}\\
    &= \sup_{\balp\in \prod_{i=1}^n[-b_i,a_i]} \langle \balp, \matP^{-1}\bp\rangle -  \frac{1}{2}\|\bx - \vz\|_M^2\\
    &= \sum_{i=1}^n \sup_{\alpha_i\in [-b_i,a_i]} \alpha_i (\matP^{-1}\bp)_i
    -  \frac{1}{2}\|\bx - \vz\|_M^2\\
    &=: h(\bp) - \frac{1}{2}\|\bx - \vz\|_M^2,
\end{split}
\end{equation*}
where, in the last line, we define the function $h\colon \Rn\to \R$ by:
\begin{equation*}
h(\bp):= \sum_{i=1}^n \sup_{\alpha_i\in [-b_i,a_i]} \alpha_i (\matP^{-1}\bp)_i\quad \forall \bp\in\Rn.
\end{equation*}
By definition, $h$ is a non-negative convex 1-homogeneous function, and hence, it is uniquely determined by its sublevel set $\{\bp\in\Rn\colon h(\bp)\leq 1\}$, which is computed as follows:
\begin{equation*}
\begin{split}
    &\{\bp\in\Rn\colon h(\bp)\leq 1\}\\
    =\,& \{\bp\in\Rn\colon \exists \bbet\in\unitsim_n \text{ s.t. }\alpha_i(\matP^{-1}\bp)_i\leq \beta_i \,\forall \alpha_i\in[-b_i,a_i] \,\forall i\in\{1,\dots,n\}\}\\
    =\,& \left\{\bp\in\Rn\colon \exists \bbet\in\unitsim_n \text{ s.t. }(\matP^{-1}\bp)_i\in \beta_i\left[-\frac{1}{b_i}, \frac{1}{a_i}\right] \,\forall i\in\{1,\dots,n\}\right\},
\end{split}
\end{equation*}
where $\unitsim_n$ denotes the standard simplex set defined by:
\begin{equation*}
    \unitsim_n:= \left\{(\beta_1,\dots,\beta_n)\in[0,1]^n\colon \sum_{i=1}^n \beta_i = 1\right\}.
\end{equation*}
Therefore, $h(\bp)\leq 1$ holds if and only if there exists $\bbet\in \unitsim_n$, such that
\begin{equation*}
    \bp = \matP(\matP^{-1}\bp) = \sum_{i=1}^n (\matP^{-1}\bp)_i \initposhd_i \in \sum_{i=1}^n\beta_i \co\left\{-\frac{\initposhd_i}{b_i},\frac{\initposhd_i}{a_i}\right\}.
\end{equation*}
Thus, we have 
\begin{equation*}
    \{\bp\in\Rn\colon h(\bp)\leq 1\} = \co\left(\bigcup_{i=1}^n \left\{\frac{\initposhd_i}{a_i}, -\frac{\initposhd_i}{b_i}\right\}\right) = \{\bx\in\Rn\colon \potentfn(\bx)\leq 1\},
\end{equation*}
which implies $h = \potentfn$, and hence, the corresponding HJ PDE~\eqref{eqt: standardHJPDE} is the HJ PDE in~\eqref{eqt: HJhd_2_general}. 

Now, we apply~\cite[Theorem~III.3.17]{Bardi1997Optimal} to prove the conclusion. If the assumptions are satisfied, then~\cite[Theorem~III.3.17]{Bardi1997Optimal} implies that the value function in the optimal control problem~\eqref{eqt: standardoptctrl} (which is~\eqref{eqt: result_optctrl2_hd_general} in our case) is the unique viscosity solution to the corresponding HJ PDE~\eqref{eqt: standardHJPDE} (which is~\eqref{eqt: HJhd_2_general} in our case). Our goal is to check the assumptions of~\cite[Theorem~III.3.17]{Bardi1997Optimal}, which include:
\begin{itemize}
    \item[($A_0$)] The set $A$ is a topological space, and the function $f\colon \Rn\times A\to \Rn$ is continuous.
    \item[($A_1$)] The function $f$ is bounded on $B_R(\Rn)\times A$ for all $R>0$. Here and after, $B_R(\Rn)$ denotes the closed ball in $\Rn$ centered at zero with radius $R$.
    \item[($A_2$)]
    There exists some positive constant $L_R$ depending on $R$, such that there holds
    $\|f(\by, \balp) -f(\bx, \balp)\|\leq L_R\|\bx - \by\|$ for all $\bx, \by \in B_R(\Rn)$, $\balp\in A$, and $R > 0$. 
    \item[(3.27)] There exists some positive constant $K$, such that
    $\|f(\bx,\balp)\|\leq K(\|\bx\|+1)$ holds for all $\bx\in\Rn$ and $\balp\in A$.
    \item[(3.40)] For all $R>0$, $\ell$ is continuous and bounded on $B_R(\Rn) \times A$ and there exists some positive constant $L_R$ depending on $R$, such that
    $|\ell(\bx, \balp) - \ell(\by, \balp)|\leq L_R\|\bx - \by\|$ holds for all $\bx, \by \in B_R(\Rn)$ and $\balp \in A$.
\end{itemize}
The assumption ($A_0$) is satisfied by definition. ($A_1$), ($A_2$), and (3.27) are satisfied because the set $A$ is compact and the function $f$ does not depend on $\bx$. Since $\ell$ does not depend on $\balp$ and $\ell$ is continuous by definition, the function $\ell$ is bounded in $B_R(\Rn) \times A$ for all $R>0$. By straightforward computation, for all $R>0$, we have that
\begin{equation*}
    |\ell(\bx,\balp) - \ell(\by,\balp)|= \frac{1}{2}\left|(\bx+\by -2\vz)^T M (\bx-\by)\right|\leq \|M\|(R +\|\vz\|)\|\bx-\by\|,
\end{equation*}
for all $\bx,\by\in B_R(\R^n)$ and $\balp\in A$. 
Hence, the assumptions of~\cite[Theorem~III.3.17]{Bardi1997Optimal} are all satisfied. As a result, the value function in~\eqref{eqt: result_optctrl2_hd_general} is the unique continuous viscosity solution to the HJ PDE~\eqref{eqt: HJhd_2_general}.
According to Proposition~\ref{prop: optctrl2_hd_general}, $\valuefn$ defined by~\eqref{eqt: result_HJ2_Lax_hd_general} is the value function of the problem~\eqref{eqt: result_optctrl2_hd_general}.
Therefore, $\valuefn$ is the unique continuous viscosity solution to the HJ PDE~\eqref{eqt: HJhd_2_general}.
\end{proof}

\section{Efficient algorithms}\label{sec: ADMM}

In this section, we present efficient algorithms for evaluating the optimal trajectory of the high-dimensional optimal control problem \eqref{eqt: result_optctrl2_hd} as well as the solution of the corresponding high-dimensional HJ PDE \eqref{eqt: result_HJ2_hd}.
We note that our algorithms may be easily adjusted to solve the more general problems given in \eqref{eqt: result_optctrl2_hd_general} and \eqref{eqt: HJhd_2_general}, but for simplicity of notation, we define our algorithms for \eqref{eqt: result_optctrl2_hd} and \eqref{eqt: result_HJ2_hd}, instead. 
To solve the more general problems in~\eqref{eqt: result_optctrl2_hd_general} and~\eqref{eqt: HJhd_2_general}, one can first compute the minimizer and the minimal value of the optimization problem in~\eqref{eqt: result_HJ2_Lax_hd_general} by applying our proposed algorithms to the new initial cost $\tilde{\initcond}$ defined in~\eqref{eqt:prop25_def_newJ} and the new terminal position $\by$ defined in~\eqref{eqt:general_hd_defy} and then compute the optimal values and optimal trajectories using~\eqref{eqt: result_HJ2_Lax_hd_general} and~\eqref{eqt: optctrl2hd_traj_general}.

Recall that the representation formulas for problems~\eqref{eqt: result_optctrl2_hd} and~\eqref{eqt: result_HJ2_hd} are provided in Section~\ref{sec: HJ_optctrl}.
Note that these problems are numerically solvable using the representation formulas if the optimization problem in~\eqref{eqt: result_Lax2_hd} is numerically solvable. 
Therefore, in this section, we provide different methods to solve~\eqref{eqt: result_Lax2_hd} for different classes of initial costs $\initcond$.
In Section~\ref{sec:quad_convex}, we consider quadratic initial costs and propose explicit formulas for solving~\eqref{eqt: result_Lax2_hd} exactly.
In Section~\ref{sec:numerical_convex}, we consider convex initial costs and apply ADMM to solve~\eqref{eqt: result_Lax2_hd}. Note that ADMM can be replaced by any other appropriate convex optimization algorithm, and we only apply ADMM in this paper for illustrative purposes. Furthermore, by applying optimization algorithms to the representation formula directly, we solve the optimal control problem exactly, without discretizations or approximations of the original problem.
In Section~\ref{sec: ADMM_nonconvex}, we extend our methods in the previous sections to address a class of nonconvex initial costs using a min-plus technique.
In each of these three sections, we provide numerical results using both a CPU implementation and an FPGA implementation to demonstrate the efficiency of our numerical solvers in various dimensions $n$ as well as their potential for real-time applications. 
All of our CPU results are run using an 11th Gen Intel Laptop Core i7-1165G7 with a 2.80GHz processor. All of our FPGA results are run using a Xilinx Alveo U280 board with a target design running at 300 MHz and double floating point precision. 
For a general overview of FPGAs, we refer the reader to \cite{KastnerFPGA}.

To avoid confusion, we use $\valuefn$ and $\opttrajhd$ to denote the analytical solutions to the HJ PDEs and optimal control problems, while we use $\Snum$ and $\gmnum$ to denote their numerical approximations obtained by our proposed methods.

\subsection{Quadratic initial costs}\label{sec:quad_convex}
\begin{figure}[htbp]
    \centering
    \begin{subfigure}{0.45\textwidth}
        \centering 
        \includegraphics[width=\textwidth]{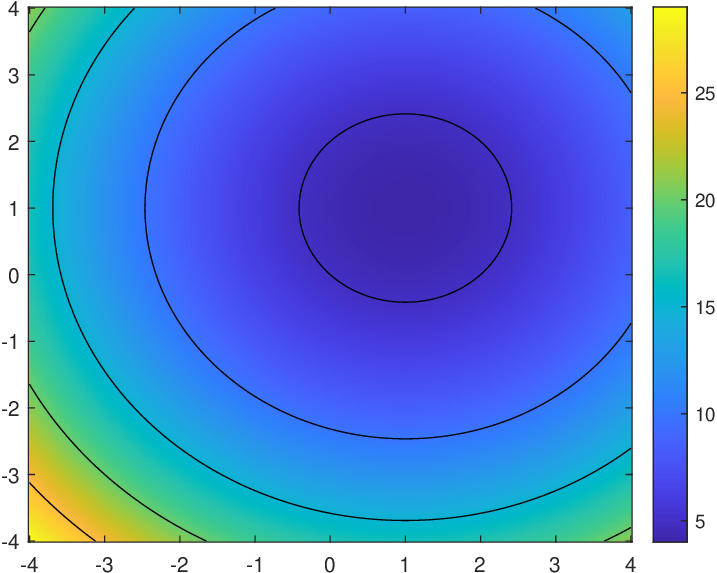}
        \caption{$t=0$}
    \end{subfigure}
    \hfill
    \begin{subfigure}{0.45\textwidth}
        \centering 
        \includegraphics[width=\textwidth]{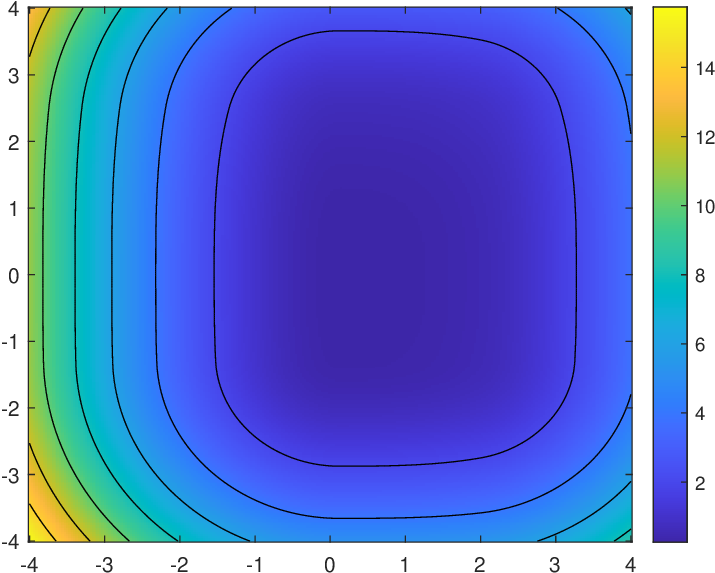}
        \caption{$t=0.25$}
    \end{subfigure}
    
    \begin{subfigure}{0.45\textwidth}
        \centering 
        \includegraphics[width=\textwidth]{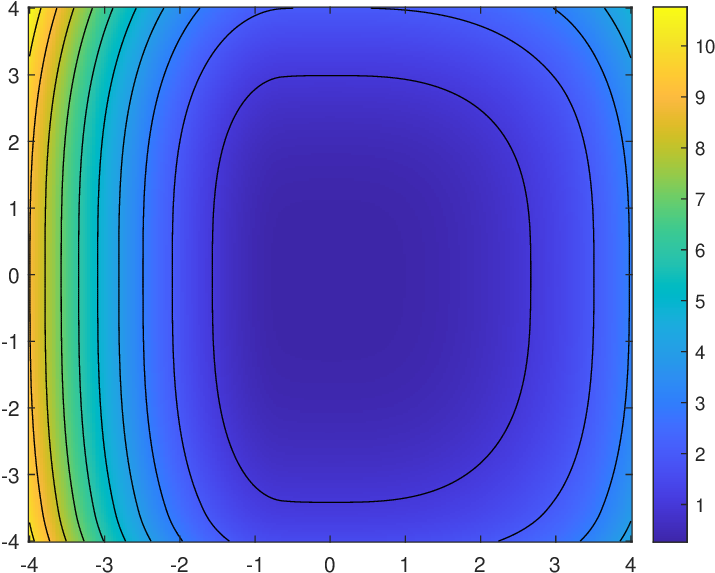}
        \caption{$t=0.5$}
    \end{subfigure}
    \hfill
    \begin{subfigure}{0.45\textwidth}
        \centering 
        \includegraphics[width=\textwidth]{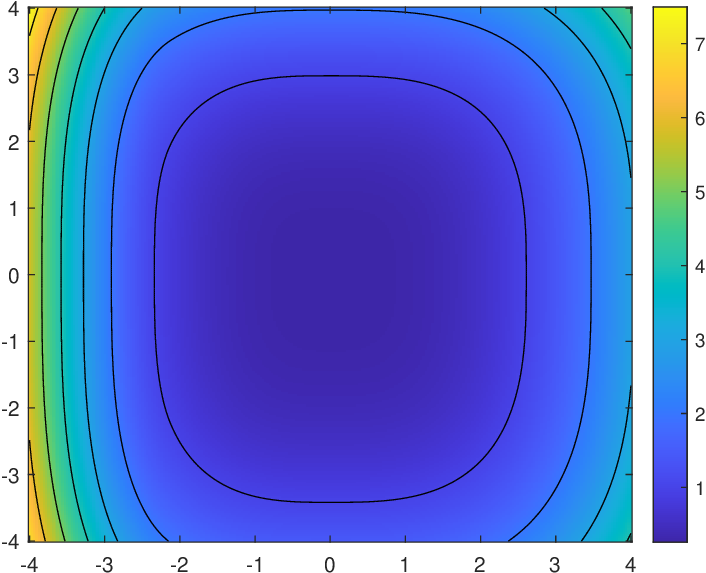}
        \caption{$t=1$}
    \end{subfigure}
    
    \caption{Evaluation of the solution $\Snum(\bx,t)$ of the $10$-dimensional HJ PDE (\ref{eqt: result_HJ2_hd}) with $\ba = (4, 6, 5, \dots, 5)$, $\bb = (3, 9, 6, \dots, 6)$, and initial condition $\initcond(\bx) = \frac{1}{2}\|\bx-\mathbf{1}\|^2$ for $\bx = (x_1, x_2, 0, \dots, 0)$, where $(x_1,x_2)\in [-4,4]^2$, and different times $t$. Plots for times $t = 0$, $0.25$, $0.5$, and $1$ are depicted in (a)-(d), respectively. Level lines are superimposed on the plots.}
    \label{fig: HJ2_quadratic_contour}
\end{figure}

In this section, we present an exact numerical solver for solving the optimal control problem~\eqref{eqt: result_optctrl2_hd} and the corresponding HJ PDE~\eqref{eqt: result_HJ2_hd} with quadratic initial costs $\initcond$. Assume that the function $\initcond$ is defined by:
\begin{equation}\label{eqt:example1_quadratic_ic}
    \initcond(\bx) = \frac{\lambda}{2}\|\bx - \by\|^2 + \alpha \quad\forall \bx\in\Rn,
\end{equation}
where $\by\in\Rn$, $\alpha\in\R$, and $\lambda>0$ are some parameters.
Recall that $\|\cdot\|$ denotes the $\ell^2$-norm in $\Rn$. 

For this quadratic initial cost, the optimization problem~\eqref{eqt: result_Lax2_hd} is equivalent to computing the proximal point of the function $\initposhd\mapsto \frac{1}{\lambda}\sum_{i=1}^n \valuefn(x_i, t; \initpos_i, a_i,b_i)$, which can be split into $n$ one-dimensional subproblems where the $i$-th subproblem reads:
\begin{equation} \label{eqt:quad_ic_subprob}
    u^*_i = \argmin_{u\in [x_i-a_it, x_i+b_it]} \left\{ \valuefn(x_i, t; u, a_i,b_i) + \frac{\lambda}{2}(u - y_i)^2\right\}.
\end{equation}
Hence, solving~\eqref{eqt: result_Lax2_hd} in this case is embarrassingly parallel; these $n$ subproblems can be solved independently in parallel using the analytical solution~\eqref{eqt: ADMM_dk_explicit_formula} in Appendix~\ref{sec:appendix_alg_proxV}, and the solution to the problems~\eqref{eqt: result_optctrl2_hd} and~\eqref{eqt: result_HJ2_hd} can be obtained directly from~\eqref{eqt: result_Lax2_hd} and~\eqref{eqt: optctrl2hd_traj_J} using the minimizer computed in~\eqref{eqt: ADMM_dk_explicit_formula}. We use our numerical solver for this case as a building block for our numerical methods in Sections~\ref{sec:numerical_convex} and~\ref{sec: ADMM_nonconvex}.

Since we implement our representation formula explicitly, the result of this numerical solver is exact up to machine precision. Furthermore, since the complexity for solving each subproblem~\eqref{eqt:quad_ic_subprob} is $\Theta(1)$ (see the discussion in Appendix~\ref{sec:appendix_alg_proxV}), the complexity for solving HJ PDE~\eqref{eqt: result_HJ2_hd} with quadratic initial cost is $\Theta(n)$ and the curse of dimensionality is avoided in this case.

Now, we apply our proposed method to solve the high-dimensional HJ PDE~\eqref{eqt: result_HJ2_hd} with quadratic initial cost $\initcond$ defined by:
\begin{equation} \label{eqt:numerical_example1_ic}
 \initcond(\bx) = \frac{1}{2}\|\bx -\mathbf{1}\|^2 \quad \forall \bx\in\Rn,
\end{equation}
i.e., we set $\by = \mathbf{1} = (1,1,\dots,1)\in\R^n$, $\alpha = 0$, and $\lambda=1$ in~\eqref{eqt:example1_quadratic_ic}.
We define the parameters $\ba=(a_1,\dots, a_n)\in \R^n$ and $\bb = (b_1,\dots, b_n)\in \R^{n}$ by:
\begin{equation} \label{eqt:numerical_def_ab}
    a_i = \begin{dcases}
    4 & \text{if } i=1,\\
    6 &\text{if } i=2,\\  
    5 &\text{if } i>2,
    \end{dcases}
    \quad\quad\text{ and }\quad\quad 
    b_i = \begin{dcases}
    3 & \text{if } i=1,\\
    9 &\text{if } i=2,\\  
    6 &\text{if } i>2.
    \end{dcases}
\end{equation}
Figure \ref{fig: HJ2_quadratic_contour} depicts two-dimensional contour plots of the numerical solution $\Snum(\bx,t)$ to this $10$-dimensional HJ PDE (i.e., $n=10$) with quadratic initial cost \eqref{eqt:numerical_example1_ic} at different positions $\bx=(x_1,x_2,0,\dots,0)$ and different times $t \in \{0,0.25, 0.5,1\}$.

\begin{table}[htbp]
    \centering
    \begin{tabular}{c|c|c|c}
        \hline
        \textbf{$\mathbf{n}$} & \textbf{CPU time (s)} & \textbf{FPGA time (s)} & \textbf{Speedup}  \\
        \hline
        4 & 6.4665e-08 & 1.334e-08 & 4.8475 \\
        8 & 1.6845e-07 & 2.667e-08 & 6.3161 \\
        12 & 4.6512e-07 & 4.000e-08 & 11.6280 \\
        16 & 7.4280e-07 & 5.334e-08 & 13.9258 \\
        \hline
    \end{tabular}
    \hfill 
    
    \caption{Comparison of the average time per call over $100,000$ runs for evaluating the solution of the HJ PDE \eqref{eqt: result_HJ2_hd} with quadratic initial condition~\eqref{eqt:numerical_example1_ic} for various dimensions $n$ using a CPU implementation on a single Intel Core i7-1165G7 versus an FPGA implementation on a Xilinx Alveo U280 board with a frequency of 300 MHz.}
    \label{tab:timing_quadratic}
\end{table}

\begin{table}[htbp]
    \centering
    \begin{tabular}{c|c|c|c|c|c}
        \hline
        \textbf{$\mathbf{n}$} & \textbf{Latency (ns)} & \textbf{BRAM} & \textbf{DSPs} & \textbf{FFs} & \textbf{LUTs} \\
        \hline
        4 & 400,224 (1.334e06) & 0 (0\%) & 847 (9\%) & 91,716 (3\%) & 55,345 (4\%) \\
        8 & 800,238 (2.667e06) & 0 (0\%) & 847 (9\%) & 92,105 (3\%) & 55,444 (4\%) \\
        12 & 1,200,246 (4.000e06) & 0 (0\%) & 847 (9\%) & 92,301 (3\%) & 55,467 (4\%) \\
        16 & 1,600,258 (5.334e06) & 0 (0\%) & 847 (9\%) & 92,626 (3\%) & 55,517 (4\%) \\
        \hline
    \end{tabular}
    \hfill 
    
    \caption{FPGA resources and latencies in cycles and nanoseconds (ns) for evaluating the solution of the HJ PDE \eqref{eqt: result_HJ2_hd} with quadratic initial condition~\eqref{eqt:numerical_example1_ic} at $100,000$ points $(\bx,t,\by)\in\Rn\times[0,\infty)\times\Rn$ for various dimensions $n$ using double precision floating points on a Xilinx Alveo U280 board with a frequency of 300 MHz.}
    \label{tab:quad_fpga_resources}
\end{table}

The running time using either a CPU or an FPGA implementation of our numerical solver in different dimensions is shown in Table~\ref{tab:timing_quadratic}. To compute the running time, we first compute the overall running time for computing the solution at $100,000$ random points $(\bx,t)\in [-4,4]^n\times [0,0.5]$ and then report the average running time for computing the solution $\Snum(\bx,t)$ at one point $(\bx,t)$ over these $100,000$ trials.
From Table~\ref{tab:timing_quadratic}, we see that, using a CPU implementation, it takes less than $8\times 10^{-7}$ seconds to compute the solution at one point in a $16$-dimensional problem, which demonstrates the efficiency of our proposed solver even in high dimensions. However, using our FPGA implementation, it takes less than $6\times 10^{-8}$ seconds to compute the solution at one point in a 16-dimensional problem, for approximately a $14\times$ speed up over the CPU implementation in dimension 16. 

We achieve this speedup by designing our FPGA implementation to have high throughput, where throughput refers to the amount of data that can be processed in a given amount of time. Specifically, we design our FPGA implementation to have an iteration interval (II) of 1, which means that we can begin processing a new input at every FPGA clock cycle (e.g., for our implementation, every 3.3333 nanoseconds). The inputs of our FPGA kernel are the points $(x_i, t, y_i)\in\R\times[0,\infty)\times\R$ as defined in \eqref{eqt:quad_ic_subprob}. In other words, our FPGA implementation streams the points $(\bx, t,\by)\in\Rn\times[0,\infty)\times\Rn$ elementwise. In contrast, the CPU implementation achieves its performance by relying on both elementwise (i.e., solving the one-dimensional subproblems \eqref{eqt:quad_ic_subprob} $n$ times) and pointwise (i.e., solving the $n$-dimensional problem \eqref{eqt: result_Lax2_hd} for multiple points $(\bx, t,\by)\in\Rn\times[0,\infty)\times\Rn$) parallelism, but must execute these parallelized tasks sequentially. Thus, as the dimension $n$ increases, the CPU is able to parallelize fewer points $(\bx, t,\by)\in\Rn\times[0,\infty)\times\Rn$ at a time, and its performance degrades by some multiplicative factor as $n$ increases. However, due to its elementwise streaming and II of 1, our FPGA implementation achieves average runtimes that only increase as $n$ times the length of one FPGA clock cycle, or, in other words, as the dimension $n$ increases, the performance of the FPGA implementation degrades only by some small additive amount. As a result, not only does our FPGA implementation achieve a speedup over the CPU implementation in lower dimensions (e.g., a speedup of about 5 in dimension 4), but this speedup becomes more pronounced as the dimension increases.

In Table~\ref{tab:quad_fpga_resources}, we present the amount of FPGA resources and latencies used to implement and run the FPGA implementation of our numerical solver for various dimemsions $n$ and 100,000 points $(\bx, t, \by)\in\Rn\times[0,\infty)\times\Rn$. We observe that since our design streams the points $(\bx, t, \by)\in\Rn\times[0,\infty)\times\Rn$ elementwise (i.e., our FPGA kernel takes the inputs $(x_i, t, y_i)\in\R\times[0,\infty)\times\R$) the latency of our FPGA implementation scales linearly in the dimension $n$ and the amount of FPGA resources used remains essentially constant in $n$. 

Note that the Alveo U280 board consists of three ``chiplets." Since routing resources between chiplets are limited, crossing chiplets can severely degrade performance \cite{Russo2020SLRCrossing,Prakash2021SLRCrossing}. As such, we design our FPGA implementation to use less than 30\% of any given type of FPGA resource (e.g., flip flops (FFs), lookup tables (LUTs), digital signal processing units (DSPs), block random access memory (BRAM), etc.) to ensure that no chiplet is crossed. Since our design uses less than 30\% of the FPGA resources available on the Xilinx Alveo U280 board, we could either use a smaller (i.e., cheaper) FPGA to implement our numerical solver with similar performance as we report here or we could parallelize by simply implementing multiple, independent copies of our FPGA kernel to maximize usage of the FPGA board. In the latter case, we could achieve a further speedup of $\times 9$ (i.e., 3 copies of our FPGA kernel per each of the 3 chiplets, ensuring that no kernel requires crossing chiplets) for a total speedup of about 44 to 125 over the CPU depending on the dimension $n$.

\subsection{Convex initial costs}\label{sec:numerical_convex}

\begin{algorithm}[htbp]
\SetAlgoLined
\SetKwInOut{Input}{Inputs}
\SetKwInOut{Output}{Outputs}
\Input{Parameters $\ba,\bb\in\Rn$, terminal position $\bx\in\Rn$, time horizon $t>0$, running time $s>0$ of the trajectory, and convex initial cost $\initcond\colon\Rn\to\R$ in the problem~\eqref{eqt: result_optctrl2_hd} and the PDE~\eqref{eqt: result_HJ2_hd}.
Parameter $\lambda>0$, initialization $\bd^0\in\Rn$, $\admmbb^0\in\Rn$, and error tolerance $\epsilon>0$ for the ADMM scheme.
}
\Output{The optimal trajectory $\gmnum(s;\bx,t)$ in the optimal control problem~\eqref{eqt: result_optctrl2_hd} and the solution value $\Snum(\bx,t)$ to the corresponding HJ PDE~\eqref{eqt: result_HJ2_hd}.}
 \For{$k = 1,2,\dots$}{
 Update $\bv^{k+1}\in\Rn$ by:
 \begin{equation}\label{eqt: ADMM_HJ2_vupdate2}
     \bv^{k+1} = \argmin_{\bv\in \Rn} \left\{ \initcond(\bv) + \frac{\lambda}{2} \left\| \bv - \bd^k + \admmbb^k\right\|^2\right\}.
 \end{equation}
 \\
 Update $\bd^{k+1}\in\Rn$, where the $i$-th element $d_i^{k+1}$ is updated by:
 \begin{equation}\label{eqt: ADMM_HJ2_dupdate}
     d_i^{k+1} = \argmin_{d_i\in \R} \left\{ \valuefn(x_i, t; d_i, a_i,b_i) + \frac{\lambda}{2} (v_i^{k+1} - d_i + \admmb_i^k)^2\right\}.
 \end{equation}
 \\
 Update $\admmbb^{k+1}\in\Rn$ by:
 \begin{equation*}
     \admmbb^{k+1} = \admmbb^k + \bv^{k+1} - \bd^{k+1}. 
 \end{equation*}
 \\
 \If{$\|\bv^{k+1} - \bv^{k}\|^2 \leq \epsilon$, $\|\bd^{k+1} - \bd^{k}\|^2 \leq \epsilon$, and $\|\bv^{k+1} - \bd^{k+1}\|^2 \leq \epsilon$}{
   set $N=k+1$ and $\initposhd^N = \bd^N$\;
   break\;
   }
 }
 Output the optimal trajectory by:
 \begin{equation*}
     \gmnum(s;\bx,t) = (\opttraj(s; x_1,t,\initpos_1^N,a_1,b_1), \dots, \opttraj(s; x_n,t,\initpos_n^N,a_n,b_n)), 
 \end{equation*}
 where the $i$-th component $\opttraj(s; x_i,t,\initpos_i^N, a_i,b_i)$ is defined by~\eqref{eqt: optctrl2_defx_1},~\eqref{eqt: optctrl2_defx_2},~\eqref{eqt: optctrl2_defx_3}, and~\eqref{eqt: optctrl2_defx_neg}.
 Also, output the solution to the HJ PDE by:
 \begin{equation*}
     \Snum(\bx,t)  = \sum_{i=1}^n \valuefn\left(x_i,t; \initpos_i^N, a_i,b_i\right) + \initcond\left(\initposhd^N\right),
 \end{equation*}
 where the $i$-th component $\valuefn\left(x_i,t; \initpos_i^N, a_i,b_i\right)$ in the summation is defined by~\eqref{eqt: result_S2_1d} and~\eqref{eqt: result_S2_1d_negative}.
 \caption{An ADMM algorithm for solving the optimal control problem~\eqref{eqt: result_optctrl2_hd} and the corresponding HJ PDE~\eqref{eqt: result_HJ2_hd} with convex initial cost. \label{alg:admm_ver2}}
\end{algorithm}

In this section, we solve~\eqref{eqt: result_optctrl2_hd} and~\eqref{eqt: result_HJ2_hd} with convex initial cost $\initcond$. To solve these problems, we need to solve the convex optimization problem in the representation formula~\eqref{eqt: result_Lax2_hd}, 
which can be solved using many possible convex optimization algorithms. Based on the discussion in Section \ref{sec:quad_convex}, proximal point-based methods would be a reasonable approach. 

For illustrative purposes, in this section, we apply ADMM (see~\cite{Glowinski2014Alternating,Boyd2011Distributed}) to~\eqref{eqt: result_Lax2_hd} with certain convex initial costs $\initcond$ whose proximal points are numerically computable.
The details of applying ADMM to this problem are given in Algorithm~\ref{alg:admm_ver2}. We emphasize that ADMM is not the only possible optimization algorithm that can be applied here. Rather, any appropriate optimization algorithm can be applied to~\eqref{eqt: result_Lax2_hd}, the choice of which depends on the properties of the function $\initcond$ and among which the use of ADMM in Algorithm~\ref{alg:admm_ver2} is simply one such possible choice.

In each iteration of ADMM in Algorithm~\ref{alg:admm_ver2}, we first update $\bv^{k+1}$ using~\eqref{eqt: ADMM_HJ2_vupdate2}. The vector $\bv^{k+1}$ is the proximal point of $\frac{\initcond}{\lambda}$ at $\bd^k-\admmbb^k$, which is assumed to be numerically computable. Then, we compute
$\bd^{k+1}$ componentwise using~\eqref{eqt: ADMM_HJ2_dupdate}, which can be solved in parallel. More specifically, we apply the solver in Appendix~\ref{sec:appendix_alg_proxV} to solve~\eqref{eqt: ADMM_HJ2_dupdate}, and hence, the complexity for computing $\bd^{k+1}$ is $\Theta(n)$. Note that the update step for $\bd^{k+1}$ in Algorithm~\ref{alg:admm_ver2} has the same form as solving~\eqref{eqt: result_Lax2_hd} with quadratic initial cost $\bx\mapsto \frac{\lambda}{2}\|\bx - \bv^{k+1} - \admmbb^{k}\|^2$. As a result, the solver proposed in Section~\ref{sec:quad_convex} serves as a building block in our ADMM algorithm (Algorithm \ref{alg:admm_ver2}), and the running time in Table~\ref{tab:timing_quadratic} underpins the running time for updating $\bd^{k+1}$ in each iteration of ADMM. Additionally, since we apply ADMM to the representation formula directly, we do not rely on discretizations or approximations of the optimal control problem. Instead, we solve the problem exactly.

In the following proposition, we prove that the optimal trajectory $\gmnum$ and the solution value $\Snum$ as computed by Algorithm \ref{alg:admm_ver2} do indeed converge to their analytical counterparts as the number of ADMM iterates $N$ approaches infinity.

\begin{proposition}\label{prop:convergence_ADMM_convexJ}
Let $\initcond\colon\Rn\to\R$ be a convex function and $\ba,\bb$ be two vectors in $(0,+\infty)^n$. Let $\bx$ be any vector in $\Rn$ and $t>0$ be any scalar. Let $\Sanaly$ and $\gmanaly$ be the functions defined in~\eqref{eqt: result_Lax2_hd} and~\eqref{eqt: optctrl2hd_traj_J}, respectively. Let $\lambda>0$ and the initialization $\bd^0,\admmbb^0\in\Rn$ be arbitrary parameters for Algorithm \ref{alg:admm_ver2}. Let $\Snum^N$ and $\gmnum^N$ be the output solution and trajectory, respectively, from Algorithm \ref{alg:admm_ver2} with iteration number $N$. Then, we have
\begin{equation}\label{eqt:prop31_conv}
    \lim_{N\to\infty} \Snum^N(\bx,t) = \Sanaly(\bx,t) \quad \text{ and } \quad \lim_{N\to\infty} \sup_{s\in[0,t]}\|\gmnum^N(s;\bx,t) - \gmanaly(s;\bx,t)\| = 0.
\end{equation}
\end{proposition}
\begin{proof}
The proof is provided in Appendix~\ref{subsec:appendix_pf_prop31}.
\end{proof}

The convergence of the output optimal trajectory $\gmnum^N$ and solution $\Snum^N$ from Algorithm~\ref{alg:admm_ver2} are proved in the proposition above. For a general convex function $\initcond$, the convergence rate of the output solution $\Snum^N(\bx,t)$ is $o(\frac{1}{\sqrt{N}})$ if the best iteration (in terms of having the smallest objective function value among the first $N$ iterations) is selected as the output. This convergence rate can be improved to $O(\frac{1}{N})$ if the output $\bu^N$ is chosen in an ergodic manner, i.e., by setting the output $\bu^N$ to be $\frac{1}{N}\sum_{k=1}^N\bd^k$.
Moreover, when the initial condition $\initcond$ satisfies stronger assumptions (for instance, if $\initcond$ is strongly convex and differentiable with Lipschitz gradient),
we obtain linear convergence for $\bu^N$, the output solution $\Snum^N$, and the output trajectory $\gmnum^N$.
For more details on the convergence rates, we refer readers to~\cite{Deng2016global,Davis2017Faster}.

\bigbreak 
Now, we show a numerical example solved using Algorithm~\ref{alg:admm_ver2} for the optimal control problem~\eqref{eqt: result_optctrl2_hd} and the HJ PDE~\eqref{eqt: result_HJ2_hd} with convex initial cost $\initcond$ defined by:
\begin{equation}\label{eqt: hd_initial_data_l12}
    \initcond(\bx) = \frac{1}{2} \|\bx - \mathbf{1}\|_1^2 \quad \forall \bx\in\R^n,
\end{equation}
where $\|\cdot\|_1$ denotes the $\ell^1$-norm in $\Rn$ and $\mathbf{1}$ is the $n$-dimensional vector whose components are all $1$'s. We set the parameters $\ba,\bb\in\Rn$ to be the vectors defined in~\eqref{eqt:numerical_def_ab}, i.e., $\ba = ( 4, 6, 5, \dots, 5)$ and $\bb = ( 3, 9, 6, \dots, 6)$.
With these parameters and initial cost, we apply the ADMM algorithm in Algorithm~\ref{alg:admm_ver2} to solve~\eqref{eqt: result_optctrl2_hd} and~\eqref{eqt: result_HJ2_hd}.
We set the parameters in Algorithm~\ref{alg:admm_ver2} to be $\lambda = 1$, $\bd^0 = \bx$, $\admmbb^0 = \mathbf{0}$, and $\epsilon = 10^{-8}$. 
In order to solve~\eqref{eqt: ADMM_HJ2_dupdate}, we apply the efficient method described in~\cite[Section~4.4]{Darbon2016Algorithms}, which has complexity $\Theta(n)$. Therefore, the complexity for each ADMM iteration in Algorithm~\ref{alg:admm_ver2} is also $\Theta(n)$. In other words, if the number of iterations is fixed, the curse of dimensionality is avoided in this example.

\begin{figure}[htbp]
    \centering
    \begin{subfigure}{0.45\textwidth}
        \centering 
        \includegraphics[width=\textwidth]{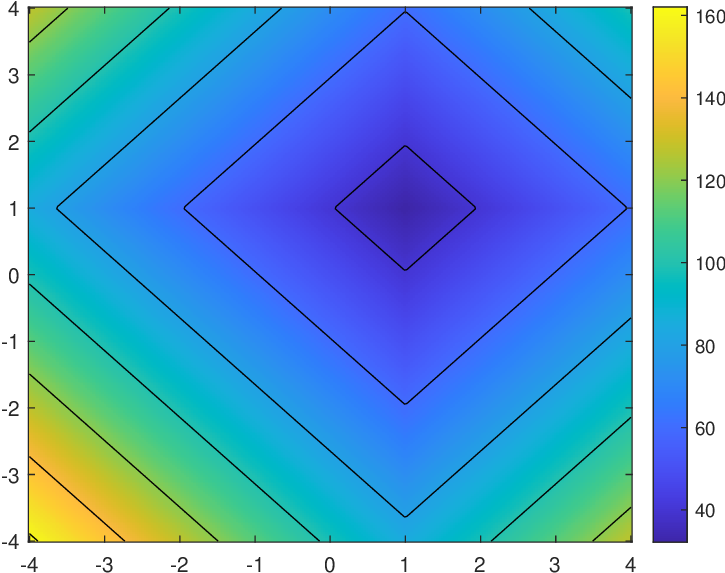}
        \caption{$t=0$}
    \end{subfigure}
    \hfill
    \begin{subfigure}{0.45\textwidth}
        \centering 
        \includegraphics[width=\textwidth]{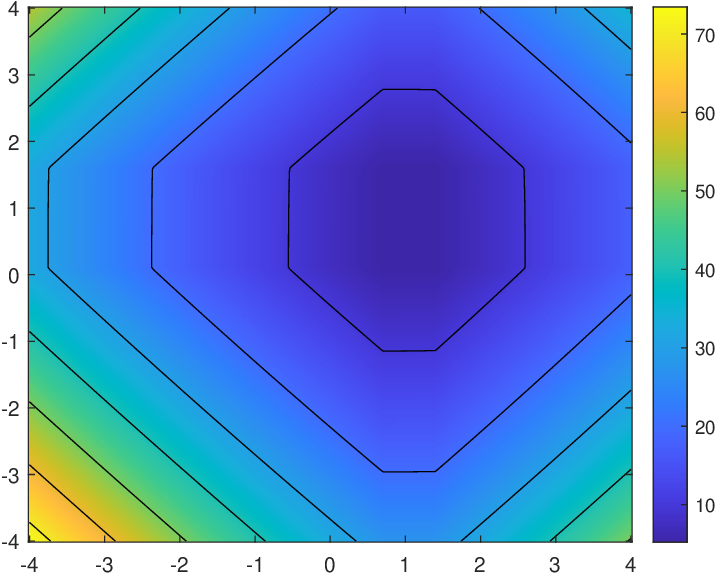}
        \caption{$t=0.1$}
    \end{subfigure}
    
    \begin{subfigure}{0.45\textwidth}
        \centering 
        \includegraphics[width=\textwidth]{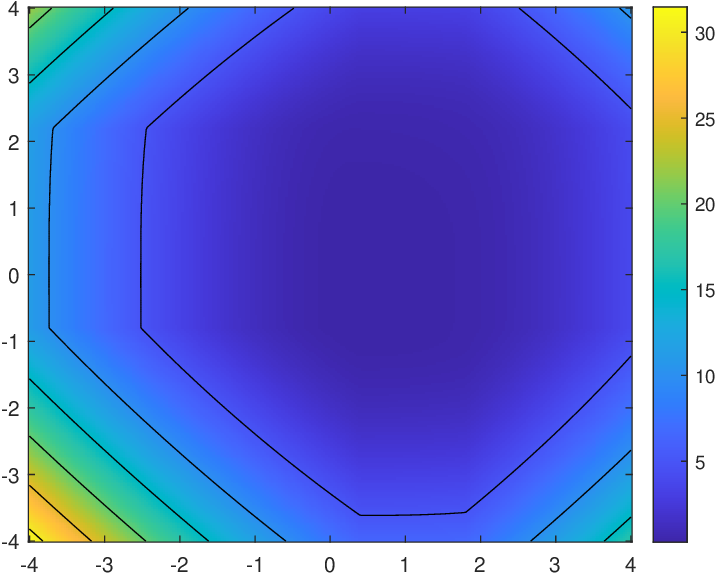}
        \caption{$t=0.2$}
    \end{subfigure}
    \hfill
    \begin{subfigure}{0.45\textwidth}
        \centering 
        \includegraphics[width=\textwidth]{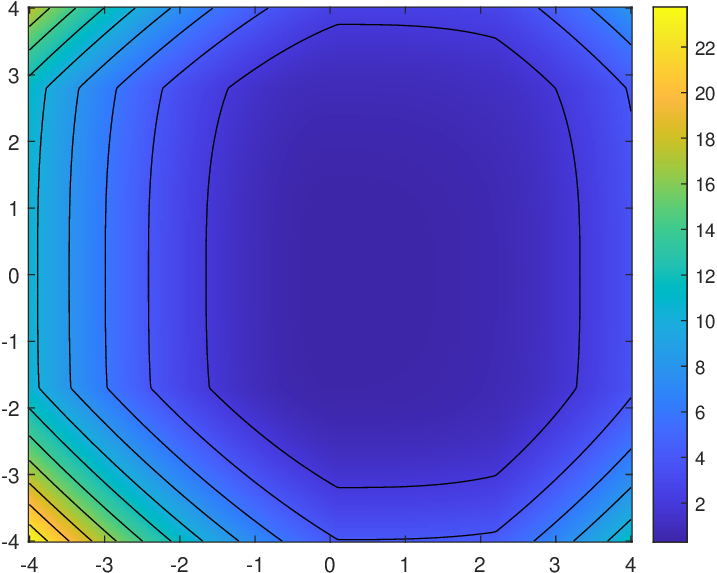}
        \caption{$t=0.3$}
    \end{subfigure}
    
    \begin{subfigure}{0.45\textwidth}
        \centering 
        \includegraphics[width=\textwidth]{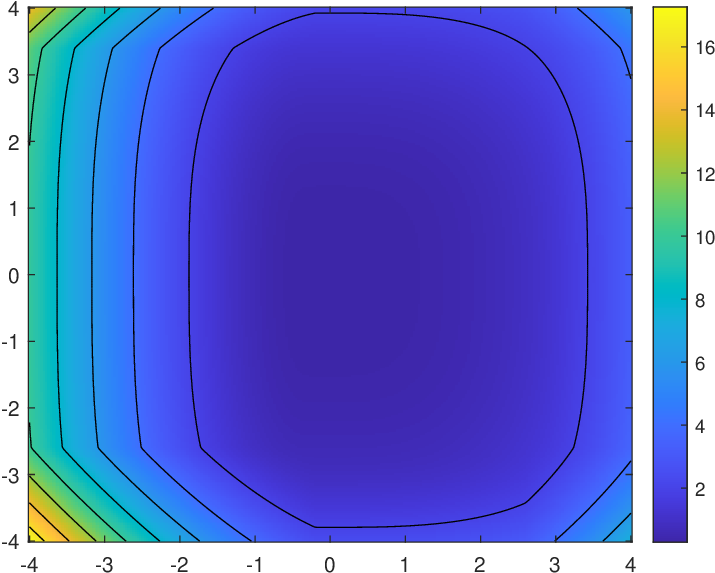}
        \caption{$t=0.4$}
    \end{subfigure}
    \hfill
    \begin{subfigure}{0.45\textwidth}
        \centering 
        \includegraphics[width=\textwidth]{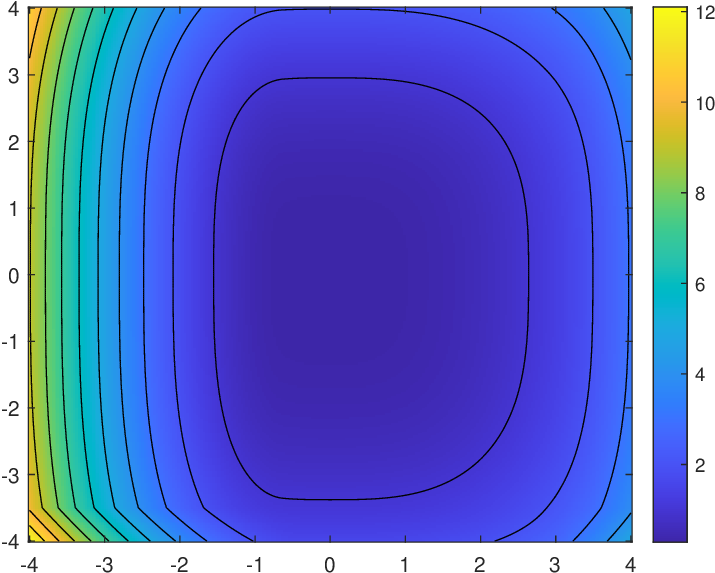}
        \caption{$t=0.5$}
    \end{subfigure}
    
    \caption{Evaluation of the solution $\Snum(\bx,t)$ of the 10-dimensional HJ PDE (\ref{eqt: result_HJ2_hd}) with $\ba = ( 4, 6, 5, \dots, 5)$, $\bb=( 3, 9, 6, \dots, 6)$, and initial condition $\initcond(\bx) = \frac{1}{2}\|\bx - \mathbf{1}\|_1^2$ for $\bx = (x_1, x_2, 0, \dots, 0)$, where $(x_1,x_2)\in [-4,4]^2$, and different times $t$. Plots for $t = 0$, $0.1$, $0.2$, $0.3$, $0.4$, and $0.5$ are depicted in (a)-(f), respectively. Level lines are superimposed on the plots.}
    \label{fig: HJ2_HD_l12}
\end{figure}

\begin{figure}[htbp]
    \centering
    \begin{subfigure}{0.32\textwidth}
        \centering \includegraphics[width=\textwidth]{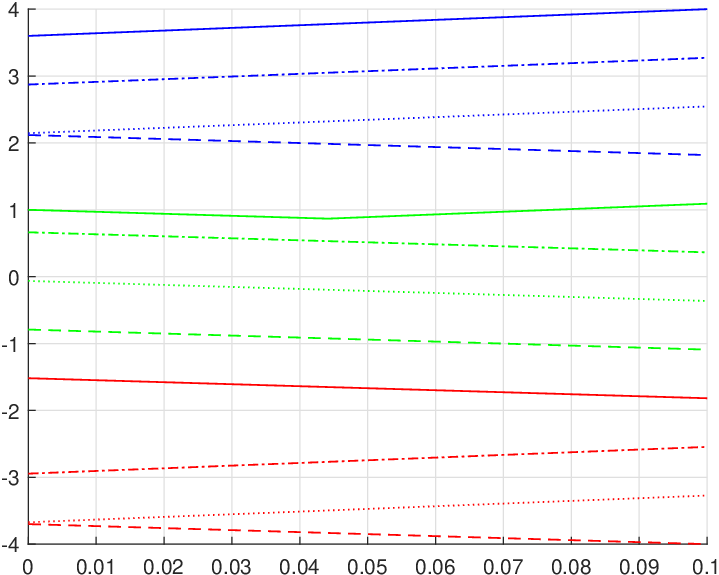}
        \caption{First component, $t=0.1$}
    \end{subfigure}
    \hfill
    \begin{subfigure}{0.32\textwidth}
        \centering \includegraphics[width=\textwidth]{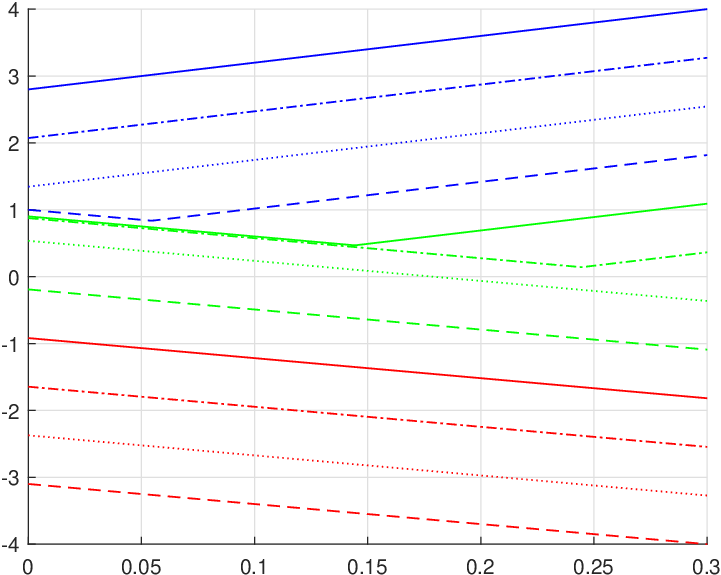}
        \caption{First component, $t=0.3$}
    \end{subfigure}
    \hfill
    \begin{subfigure}{0.32\textwidth}
        \centering \includegraphics[width=\textwidth]{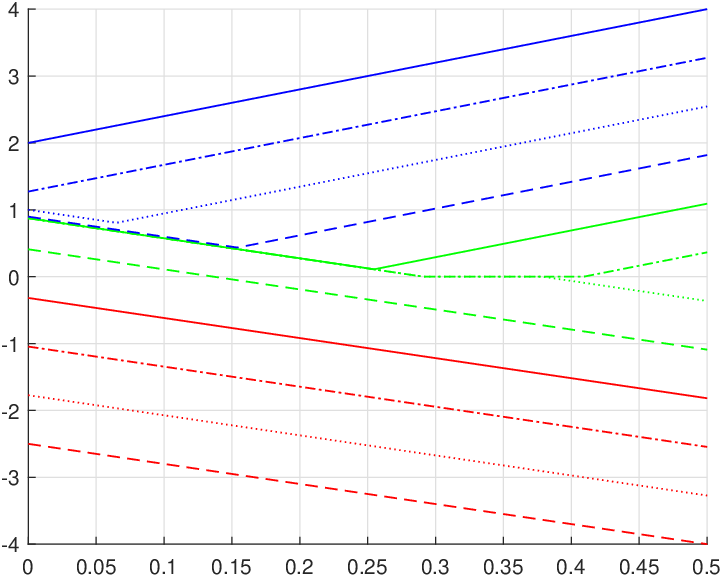}
        \caption{First component, $t=0.5$}
    \end{subfigure}
    
    \begin{subfigure}{0.32\textwidth}
        \centering \includegraphics[width=\textwidth]{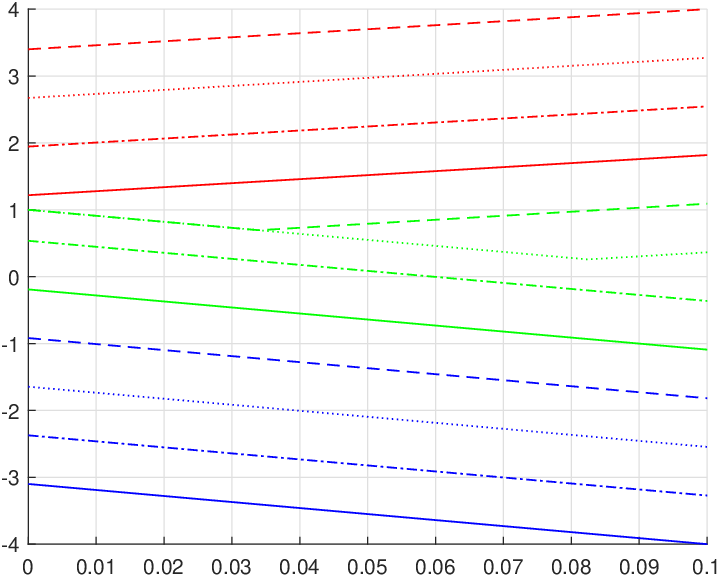}
        \caption{Second component, $t=0.1$}
    \end{subfigure}
    \hfill
    \begin{subfigure}{0.32\textwidth}
        \centering \includegraphics[width=\textwidth]{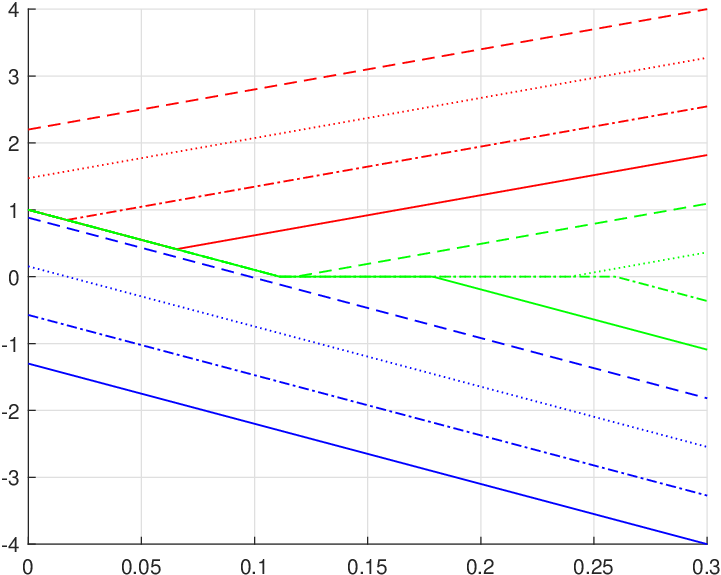}
        \caption{Second component, $t=0.3$}
    \end{subfigure}
    \hfill
    \begin{subfigure}{0.32\textwidth}
        \centering \includegraphics[width=\textwidth]{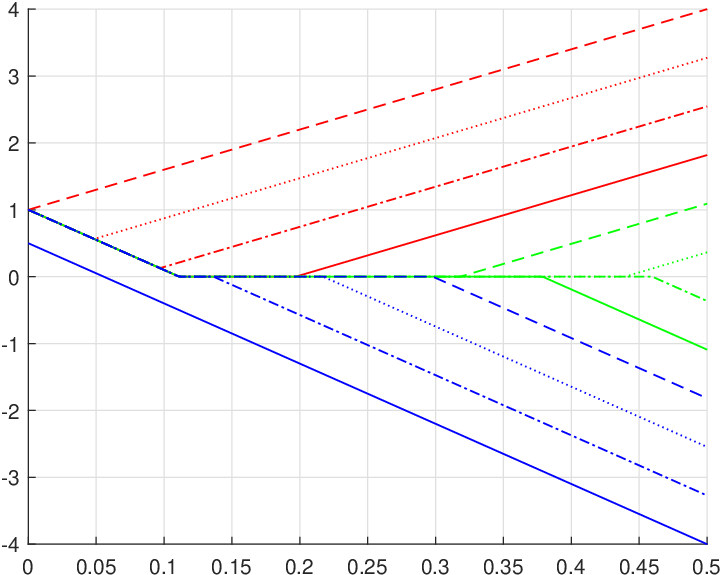}
        \caption{Second component, $t=0.5$}
    \end{subfigure}

    \caption{Evaluation of the optimal trajectory $\gmnum(s;(x,-x,0,\dots,0), t)$ of the $10$-dimensional optimal control problem (\ref{eqt: result_optctrl2_hd}) with $\ba = ( 4, 6, 5, \dots, 5)$, $\bb = ( 3, 9, 6, \dots, 6)$, and initial cost $\initcond(\bx) = \frac{1}{2}\|\bx - \mathbf{1}\|_1^2$ versus $s\in[0,t]$ for different terminal positions $(x,-x,0,\dots,0)$ ($x\in [-4,4]$) and different time horizons $t$. The different colors and line markers simply help differentiate between the different trajectories. Figures (a)-(c) depict the first component of the trajectory versus $s\in[0,t]$ with different time horizons $t$, while Figures (d)-(f) depict the second component of the trajectory versus $s\in[0,t]$ with different time horizons $t$. Plots for time horizons $t = 0.1$, $0.3$, and $0.5$ are depicted in Figures (a)/(d), (b)/(e), and (c)/(f), respectively. We note that because of our choice of $\ba$ and $\bb$, the piecewise slopes of our trajectories are not symmetric about $0$.}
    \label{fig: HJ2_HD_l12_trajectory}
\end{figure}

We solve the problem in $10$ dimensions (i.e., we set $n=10$) and plot the solution $\Snum$ and the optimal trajectories $\gmnum$ in Figure~\ref{fig: HJ2_HD_l12} and Figure~\ref{fig: HJ2_HD_l12_trajectory}, respectively. 
Figure \ref{fig: HJ2_HD_l12} depicts two-dimensional slices of the solution $\Snum(\bx,t)$, as computed using Algorithm~\ref{alg:admm_ver2}, of the HJ PDE \eqref{eqt: result_HJ2_hd}
at different positions $\bx=(x_1, x_2, 0, \dots, 0)$ and at different times $t$. As expected, in Figure~\ref{fig: HJ2_HD_l12}(a), we see that the initial condition $\initcond$ is not smooth, e.g., we see kinks in the contour plots near $(x_1,x_2) = (1,-1), (-2,1), (4,1),$ and $(1,4)$. In Figures~\ref{fig: HJ2_HD_l12}(b)-(f), we see that the solution continues to evolve with several kinks as well. These kinks help numerically verify that our algorithm does indeed compute the non-smooth viscosity solution to the corresponding HJ PDE. Overall, the solution appears to be continuous in $(x_1, x_2)$ at all times $t$, which is consistent with the results of Proposition \ref{prop: HJ_hd}.

Figure \ref{fig: HJ2_HD_l12_trajectory} depicts one-dimensional slices of the optimal trajectory $\gmnum(s;\bx,t)$ of the corresponding $10$-dimensional optimal control problem~\eqref{eqt: result_optctrl2_hd}, using different terminal positions $\bx=(x,-x,0,\dots,0)$ for $x\in[-4,4]$ and different time horizons $t$. We observe that the one-dimensional slices are piecewise linear and continuous in $s$, which is consistent with the properties of the formulas in~\eqref{eqt: optctrl2_defx_1},~\eqref{eqt: optctrl2_defx_2},~\eqref{eqt: optctrl2_defx_3}, and~\eqref{eqt: optctrl2_defx_neg}. We note that in each subplot, all line segments with positive slope are parallel with slope $a_i$ (i.e., the $i$-th component of the trajectory has velocity $a_i$), while all lines segments with negative slope are parallel with slope $-b_i$ (i.e., the $i$-th component of the trajectory has velocity $-b_i$). As such, in any given subplot, the piecewise slopes of the trajectories are not symmetric about $0$ due to our choice of $\ba$ and $\bb$. 

\begin{table}[htbp]
    \begin{subtable}[h]{\textwidth}
    \centering
    \begin{tabular}{c|c|c|c}
        \hline
        \textbf{$\mathbf{n}$} & \textbf{CPU time (s)} & \textbf{FPGA time (s)} & \textbf{Speedup}  \\
        \hline
        4 & 5.3711e-08 & 9.1888e-09 & 5.8453 \\
        8 & 1.1719e-07 & 3.6775e-08 & 3.1867 \\
        12 & 1.8880e-07 & 6.9133e-08 & 2.7310 \\
        16 & 2.7344e-07 & 1.1555e-07 & 2.3664 \\
        \hline
    \end{tabular}
    \caption{Comparison of the average time per 1 iteration of ADMM over 100,000 runs for various dimensions $n$. }
    \label{tab:timing_l1sq_highthroughput}
    \end{subtable}
    \vspace{5pt}
    
    \begin{subtable}[h]{\textwidth}
    \centering
    \begin{adjustbox}{width=\textwidth}
    \begin{tabular}{c|c|c|c|c|c|c|c}
        \hline
        \textbf{$\mathbf{n}$} & \textbf{$N$} &  \textbf{Latency (ns)} & \textbf{Interval (ns)} & \textbf{BRAM} & \textbf{DSPs} & \textbf{FFs} & \textbf{LUTs} \\
        \hline
        4 & 8 & 2,205,620 (7.351e06) & 2,200,087 (7.334e06) & 0 (0\%)  & 5,840 (64\%) & 2,606,704  (99\%) & 1,160,176 (88\%) \\
        8 & 4 & 4,413,012 (1.471e07) & 4,400,085 (1.467e07) & 0 (0\%)  & 4,955 (54\%) & 2,630,997 (100\%) & 870,327 (66\%) \\
        12 & 3& 6,222,468 (2.074e07) & 6,200,127 (2.067e07) & 0 (0\%)  & 4,949 (54\%) & 2,515,392 (96\%) & 948,691 (72\%) \\
        16 & 2 & 6,933,965 (2.311e07) & 6,900,222 (2.300e07)& 0 (0\%)  & 3,341 (37\%) & 2,149,905 (82\%) & 740,551 (56\%) \\
        \hline
    \end{tabular}
    \end{adjustbox}
    \caption{Latency and interval in cycles and nanoseconds (ns) and FPGA resources used to implement $N$ ADMM iterations for evaluating the solution at 100,000 points $(\bx,t)\in\Rn\times[0,\infty)$.}
    \label{tab:l1sq_fpga_resources_highthroughput}
    \end{subtable}
    
    \caption{Results for a high throughput FPGA implementation of ADMM for evaluating the solution of the HJ PDE (\ref{eqt: result_HJ2_hd}) with initial condition $\initcond(\bx) = \frac{1}{2}\|\bx - \mathbf{1}\|_1^2$. The FPGA implementation uses a Xilinx Alveo U280 board with a frequency of 300 MHz. In (a), the CPU implementation uses a single Intel Core i7-1165G7, and the number of ADMM iterates used per run is the same as the number of iterates $N$ used in (b) for each dimension $n$.}
    \label{tab:l1sq_fpga_highthroughput}
\end{table}

\begin{table}[htbp]
    \begin{subtable}[h]{\textwidth}
    \centering
    \begin{tabular}{c|c|c|c}
        \hline
        \textbf{$\mathbf{n}$} & \textbf{CPU time (s)} & \textbf{FPGA time (s)} & \textbf{Speedup}  \\
        \hline
        4 & 5.4254e-08 & 8.4816e-08 & 0.6397 \\
        8 & 1.1161e-07 & 1.1857e-07 & 0.9413 \\
        12 & 1.8359e-07 & 1.8067e-07 & 1.0162 \\
        16 & 2.5879e-07 & 2.3167e-07 & 1.1171 \\
        \hline
    \end{tabular}
    \caption{Comparison of the average time per 1 iteration of ADMM over 100,000 runs for various dimensions $n$. }
    \label{tab:timing_l1sq_lowlatency}
    \end{subtable}
    \vspace{5pt}

    \begin{subtable}[h]{\textwidth}
    \centering
    \begin{adjustbox}{width=\textwidth}
    \begin{tabular}{c|c|c|c|c|c|c|c}
        \hline
        \textbf{$\mathbf{n}$} & \textbf{$N$} &  \textbf{Latency (ns)} & \textbf{Interval (ns)} & \textbf{BRAM} & \textbf{DSPs} & \textbf{FFs} & \textbf{LUTs} \\
        \hline
        4 & 9 & 433 (1.443e03)
        & 229 (7.633e02)
        & 0 (0\%)  & 2,608 (28\%) & 612,037 (23\%) & 390,421 (29\%) \\
        8 & 7 & 405 (1.350e03)
        & 249 (8.300e02) 
        & 0 (0\%) & 2,211 (24\%) & 616,758 (23\%) & 377,262 (28\%) \\
        12 & 5 & 377 (1.257e03)
        & 271 (9.033e02)
        & 0 (0\%) & 1,949 (21\%) & 629,957 (24\%) & 387,615 (29\%) \\
        16 & 4 & 371 (1.237e03)
        & 278 (9.267e02)
        & 0 (0\%) & 1,898 (21\%) & 620,606 (23\%) & 382,723 (29\%) \\
        \hline
    \end{tabular}
    \end{adjustbox}
    \caption{Latency and interval in cycles and nanoseconds (ns) and FPGA resources used to implement $N$ ADMM iterations for evaluating the solution at 1 point $(\bx,t)\in\Rn\times[0,\infty)$. }
    \label{tab:l1sq_fpga_resources_lowlatency}
    \end{subtable}

    \caption{Results for a low latency FPGA implementation of ADMM for evaluating the solution of the HJ PDE (\ref{eqt: result_HJ2_hd}) with initial condition $\initcond(\bx) = \frac{1}{2}\|\bx - \mathbf{1}\|_1^2$. The FPGA implementation uses a Xilinx Alveo U280 board with a frequency of 300 MHz. In (a), the CPU implementation uses a single Intel Core i7-1165G7, and the number of ADMM iterates used per run is the same as the number of iterates $N$ used in (b) for each dimension $n$.}
    \label{tab:l1sq_fpga_lowlatency}
\end{table}

Next, we describe two different FPGA implementations of ADMM for this example. Specifically, we present a high throughput implementation and a low latency implementation, the results for which are shown in Tables \ref{tab:l1sq_fpga_highthroughput} and \ref{tab:l1sq_fpga_lowlatency}, respectively. Latency refers to the amount of time that it takes for an input to finish being processed. The ``optimality" of a given FPGA implementation is usually determined by its latency (where lower latency is more optimal), its throughput (where higher throughput is more optimal), or the amount of resources used (where fewer resources is more optimal). However, optimizing for one of these criteria usually competes with the optimization of another criteria, and thus, one can at best only expect a Pareto optimal implementation. For example, a high throughput implementation typically has a relatively high latency and vice versa. In fact, we observe that this trend holds for our FPGA implementations. For instance, our high throughput implementation has latencies (in cycles) per ADMM iteration of approximately 694.4, 3,242.8, 7,467.7, and 16,906.0 for dimensions $n=4, 8, 12,$ and 16, respectively. Meanwhile, our low latency implementation has latencies (in cycles) per ADMM iteration of approximately 48.1, 57.9, 75.4, and 92.8 for dimensions $n=4, 8, 12$, and 16, respectively. Note that the latency and interval listed in Tables \ref{tab:l1sq_fpga_resources_highthroughput} and \ref{tab:l1sq_fpga_resources_lowlatency} correspond to the latency and interval for processing $N$ ADMM iterations for all 100,000 points (in Table \ref{tab:l1sq_fpga_resources_highthroughput}) or 1 point (in Table \ref{tab:l1sq_fpga_resources_lowlatency}), respectively. Thus, we can compute the latency per ADMM iteration using the values in Tables \ref{tab:l1sq_fpga_resources_highthroughput} and \ref{tab:l1sq_fpga_resources_lowlatency} as follows:
$$\frac{\text{latency} - \text{interval} +  \text{interval}/\text{\# points}}{N},$$
where \# points = 100,000 in Table \ref{tab:l1sq_fpga_resources_highthroughput}, \# points = 1 in Table \ref{tab:l1sq_fpga_resources_lowlatency}, and the numerator in the above formula corresponds to the latency (for $N$ ADMM iterations) per point.

The benefit of a high throughput implementation is that it achieves low average runtime. Hence, high throughput implementations are best suited for offline computations, where a computational experiment needs to be run many times for many different inputs. We can expect a high throughput implementation to achieve the best speedup compared to a CPU implementation, whose performance is typically also measured by average runtime.   For example, in Table \ref{tab:timing_l1sq_highthroughput}, we see that our high throughput implementation achieves a speedup of about 2-6 times the average runtime of the CPU depending on the dimension $n$, whereas in Table \ref{tab:timing_l1sq_lowlatency}, we see that the average runtime of our low latency implementation is approximately the same as that of the CPU for each $n$. Here, we compute the average runtime using the procedure described in Section \ref{sec:quad_convex}, except we average over the number of ADMM iterations in addition to the number of runs.

In contrast, low latencies are best suited for online computations, where results are required to be available within some fixed short period of time after an input is provided. Such online computations are critical for real-time optimal control applications. While it is impossible to measure the latency of a CPU implementation, FPGA implementations have guaranteed latencies.  For example, as computed above, our low latency implementation computes one iteration of ADMM in 92.8 clock cycles (or approximately $3.0933 \times 10^{-7}$ seconds) for the 16-dimensional problem, but this quantity would not be able to be measured on a CPU. Hence, our low latency FPGA implementation achieves similar performance (in terms of throughput; e.g., see Table \ref{tab:timing_l1sq_lowlatency}) as the CPU but with a guaranteed (low) latency. 

In both FPGA implementations, we stream both the points $(\bx, t)\in\Rn\times[0,\infty)$ and the ADMM iterates $\bv^k, \bd^k, \admmbb^k\in\Rn$ between consecutive ADMM iterations, where the number of ADMM iterations per implementation is determined by the the amount of resources available. For the high throughput implementation, the amount of resources used per ADMM iteration is high (and increases with the dimension $n$), and we aim to use as many of the resources as possible in order to maximize the throughput. Thus, our high throughput FPGA kernel must cross chiplets. In order to ensure that the performance does not degrade due to the crossing of chiplets, we stream $(\bx, t, \bv^k, \bd^k, \admmbb^k)$ in a single stream of doubles. As a result, (in contrast with the high throughput FPGA implementation in Section~\ref{sec:quad_convex}) our high throughput implementation for ADMM cannot achieve an II of 1. Instead, the II of our high throughput FPGA kernel (per point) is lower bounded by $4n + 1$, the dimension of the concatenated vector $(\bx, t, \bv^k, \bd^k, \admmbb^k)$, which means that its II per ADMM iteration is lower bounded by $(4n + 1)/N$, where $N$ is the total number of ADMM iterations implemented. For example, for our high throughput implementation, the II per ADMM iteration is about 2.8, 11.0, 20.7, and 34.5 cycles for dimensions $n = 4, 8, 12,$ and 16, respectively. Note that the II per ADMM iteration can be computed using the quantities in Table \ref{tab:l1sq_fpga_resources_highthroughput} as $\frac{\text{interval}}{\text{\# points} \times \text{N}}.$ We also note that the average runtime of our high throughput FPGA implementation would be improved if we could implement more ADMM iterations on the FPGA. However, as we observe in Table \ref{tab:l1sq_fpga_resources_highthroughput}, our current implementation is heavily limited by the number of flip flops used. Theoretically, we should be able to reduce the number of flip flops using other FPGA resources, such as BRAM or URAM, instead, but we leave this for future research.

In contrast, for the low latency FPGA implementation, we stream each of the quantities $\bx, t, \bv^k, \bd^k,$ and $\admmbb^k$ separately. Using separate streams for these quantities ensures that these quantities are available more immediately for computations, which is critical for achieving a low latency. However, using multiple streams also means that crossing chiplets is likely to cause a significant decrease in performance. Hence, we aim for designs that use less than 30\% of any given FPGA resource to ensure that no chiplet is crossed. In Table \ref{tab:l1sq_fpga_resources_lowlatency}, we see that our low latency FPGA implementation meets this constraint.

\subsection{Certain nonconvex initial costs}\label{sec: ADMM_nonconvex}

In this section, we use min-plus techniques to extend our Lax-Oleinik-type representation formulas~\eqref{eqt: result_Lax2_1d},~\eqref{eqt: result_Lax2_hd}, and~\eqref{eqt: result_HJ2_Lax_hd_general} to handle a certain class of nonconvex initial costs.
Moreover, we propose an algorithm based on the resulting extended representation formulas, which uses the numerical methods in Sections~\ref{sec:quad_convex} and~\ref{sec:numerical_convex} (or any possible algorithm for solving~\eqref{eqt: result_Lax2_1d},~\eqref{eqt: result_Lax2_hd}, and~\eqref{eqt: result_HJ2_Lax_hd_general}) as building blocks. 
From the numerical results and resulting running times, our proposed algorithm is shown to be able to solve the optimal control problems and corresponding HJ PDEs with these nonconvex initial costs efficiently.

We have already shown in Sections~\ref{sec:quad_convex} and~\ref{sec:numerical_convex} that the solutions in~\eqref{eqt: result_Lax2_1d}, \eqref{eqt: result_Lax2_hd}, and~\eqref{eqt: result_HJ2_Lax_hd_general} are computable using convex optimization methods, such as Algorithm~\ref{alg:admm_ver2} if the initial cost $\initcond\colon\Rn\to\R$ is a convex function. 
However, these representation formulas are also computable for a broader class of initial costs $\initcond$. Consider the following nonconvex initial condition:
\begin{equation}\label{eqt: J_minplus}
    \initcond(\bx) = \min_{j\in\{1,\dots, m\}} \initcond_j(\bx) \quad\forall \bx\in\Rn,
\end{equation}
where $\initcond_j\colon\Rn\to\R$ is a convex function for each $j\in\{1,\dots,m\}$.
In this case, the min-plus technique (see~\cite{Kolokoltsov1997Idempotent,McEneaney2006maxplus}) is applied, 
and the optimization problem in~\eqref{eqt: result_HJ2_Lax_hd_general} can be written as
\begin{equation}\label{eqt: pde2_SJ_general_minplus}
\begin{split}
    &\valuefn(\bx,t)\\
    =& 
    \inf_{\initposhd\in\prod_{i=1}^n [y_i-a_it, y_i+b_it]} \left\{\sum_{i=1}^n \valuefn(y_i, t; \initpos_i, a_i,b_i) + \min_{j\in\{1,\dots, m\}}\initcond_j\left((\matP^T)^{-1}\initposhd+\vz\right)\right\}\\
    =& \min_{j\in\{1,\dots, m\}}\left\{
    \inf_{\initposhd\in\prod_{i=1}^n [y_i-a_it, y_i+b_it]} \left\{\sum_{i=1}^n \valuefn(y_i, t; \initpos_i, a_i,b_i) + \initcond_j\left((\matP^T)^{-1}\initposhd+\vz\right)\right\}\right\}\\
    =&:\min_{j\in\{1,\dots, m\}} \valuefn_j(\bx,t),
\end{split}
\end{equation}
where $\by=(y_1,\dots,y_n)$ is the vector defined in~\eqref{eqt:general_hd_defy}, and in the last line, the function $\valuefn_j\colon \Rn\times [0,+\infty) \to\R$ for each $j\in\{1,\dots, m\}$ is defined by:
\begin{equation}\label{eqt: minplus_PDE_jsubproblem}
    \valuefn_j(\bx,t):= \inf_{\initposhd\in\prod_{i=1}^n [y_i-a_it, y_i+b_it]} \left\{\sum_{i=1}^n \valuefn(y_i, t; \initpos_i, a_i,b_i) + \initcond_j\left((\matP^T)^{-1}\initposhd+\vz\right)\right\},
\end{equation}
which is the solution to the corresponding HJ PDE with convex initial cost $\initcond_j$.
Therefore, to compute $\valuefn(\bx,t)$, this problem is divided into $m$ subproblems. In the $j$-th subproblem, convex optimization methods (e.g., the solver in Section~\ref{sec:quad_convex} or  Algorithm~\ref{alg:admm_ver2} in Section \ref{sec:numerical_convex}) are applied to solve~\eqref{eqt: minplus_PDE_jsubproblem} and to compute the optimal value $\valuefn_j(\bx,t)$.
Then, by~\eqref{eqt: pde2_SJ_general_minplus}, the solution $\valuefn(\bx,t)$ is the minimum among these $m$ optimal values.

In addition, a minimizer $\initposhd^*$ of~\eqref{eqt: result_HJ2_Lax_hd_general} can be computed as
\begin{equation} \label{eqt: opttraj_ustar_minplus}
    \initposhd^*\in \bigcup_{r\in \mathcal{J}} \argmin_{\initposhd\in\prod_{i=1}^n [y_i-a_it, y_i+b_it]} \left\{\sum_{i=1}^n \valuefn(y_i, t; \initpos_i, a_i,b_i) + \initcond_r\left((\matP^T)^{-1}\initposhd+\vz\right)\right\},
\end{equation}
where the set $\mathcal{J}$ is defined by:
\begin{equation*}
    \mathcal{J}:= \argmin_{j\in\{1,\dots,m\}}
    \valuefn_j(\bx,t),
\end{equation*}
with $\valuefn_j$ defined by~\eqref{eqt: minplus_PDE_jsubproblem}. In other words, $\initposhd^*$ can be computed using~\eqref{eqt: opttraj_ustar_minplus} by applying convex optimization methods to solve each of the $m$ subproblems, where the $j$-th subproblem is defined in~\eqref{eqt: minplus_PDE_jsubproblem}. Finally, an optimal trajectory $\opttrajgeneral$ is computed using~\eqref{eqt: optctrl2hd_traj_general} as long as the minimizer $\initposhd^*$ is obtained. Note that $\opttrajgeneral$ may not be unique since we are solving a nonconvex optimization problem and hence $\bu^*$ is no longer necessarily unique.

Similarly, the min-plus technique can also be applied to the representation formulas~\eqref{eqt: result_Lax2_1d} and~\eqref{eqt: result_Lax2_hd} to solve the corresponding optimal control problems and HJ PDEs with nonconvex initial condition $\initcond$ of the form~\eqref{eqt: J_minplus}. 
To be specific,~\eqref{eqt: result_Lax2_1d} and~\eqref{eqt: result_Lax2_hd} are the one-dimensional and high-dimensional cases, respectively, of~\eqref{eqt: result_HJ2_Lax_hd_general} when $\matP$ is the identity matrix and $\vz$ is the zero vector. Therefore, in these cases, the $j$-th subproblem becomes
\begin{equation}\label{eqt: minplus_PDE_jsubproblem_separable}
    \valuefn_j(\bx,t):= \inf_{\initposhd\in\prod_{i=1}^n [x_i-a_it, x_i+b_it]} \left\{\sum_{i=1}^n \valuefn(x_i, t; \initpos_i, a_i,b_i) + \initcond_j\left(\initposhd\right)\right\},
\end{equation}
and an optimal trajectory $\opttrajhd$ and the optimal value $\valuefn$ are computed using the minimizers and the minimal values of these subproblems, similarly to before. The details of the proposed algorithm for solving the optimal control problem~\eqref{eqt: result_optctrl2_hd} and the corresponding HJ PDE~\eqref{eqt: result_HJ2_hd} are given in Algorithm~\ref{alg:admm_minplus}. The more general problems in~\eqref{eqt: result_optctrl2_hd_general} and~\eqref{eqt: HJhd_2_general} can be solved by replacing the $j$-th subproblem in Algorithm~\ref{alg:admm_minplus} with~\eqref{eqt: minplus_PDE_jsubproblem} for each $j\in\{1,\dots,m\}$. Note that each subproblem is solved independently from each other, and hence, each subproblem can be solved in parallel and/or using different numerical methods, when convenient. The complexity of Algorithm~\ref{alg:admm_minplus} is the sum of the complexity of solving all $m$ subproblems. 

\begin{algorithm}[htbp]
\SetAlgoLined
\SetKwInOut{Input}{Inputs}
\SetKwInOut{Output}{Outputs}
\Input{Parameters $\ba,\bb\in\Rn$, terminal position $\bx\in\Rn$, time horizon $t>0$, running time $s>0$ of the trajectory, and the convex functions $\initcond_1,\dots,\initcond_m$ in~\eqref{eqt: J_minplus}.
}
\Output{An optimal trajectory $\gmnum(s;\bx,t)$ in the optimal control problem~\eqref{eqt: result_optctrl2_hd} and the solution value $\Snum(\bx,t)$ to the corresponding HJ PDE~\eqref{eqt: result_HJ2_hd} with nonconvex initial cost $\initcond$ of the form~\eqref{eqt: J_minplus}.}
 \For{$j = 1,2,\dots,m$}{
 Numerically solve the $j$-th subproblem~\eqref{eqt: minplus_PDE_jsubproblem_separable} using any appropriate method (e.g., the solver in Section~\ref{sec:quad_convex} or Algorithm~\ref{alg:admm_ver2} in Section~\ref{sec:numerical_convex}), and get the optimal trajectory $\gmnum_j(s;\bx,t)$ and the solution $\Snum_j(\bx,t)$ to the problems~\eqref{eqt: result_optctrl2_hd} and~\eqref{eqt: result_HJ2_hd} with initial cost $\initcond_j$\;
 }
 Compute the index $r$ by:
 \begin{equation}\label{eqt:minplus_defr}
    r \in \argmin_{j\in\{1, \dots, m\}} \Snum_j (\bx,t).
\end{equation}
\\
 Output an optimal trajectory $\gmnum(s;\bx,t)$ and the solution value $\Snum(\bx,t)$ using
 \begin{equation}\label{eqt:alg_minplus_output}
 \begin{split}
    \gmnum(s;\bx,t) & = \gmnum_r(s;\bx,t), \quad\quad    \Snum(\bx,t)  = \Snum_r(\bx,t) = \min_{j\in\{1, \dots, m\}} \Snum_j (\bx,t).
\end{split}
\end{equation}
 \caption{An optimization algorithm for solving the optimal control problem~\eqref{eqt: result_optctrl2_hd} and the corresponding HJ PDE~\eqref{eqt: result_HJ2_hd} with nonconvex initial cost $\initcond$ of the form~\eqref{eqt: J_minplus}. \label{alg:admm_minplus}}
\end{algorithm}

In the following proposition, we show some error analysis for Algorithm~\ref{alg:admm_minplus}, given the error of each subproblem. To be specific, if each subproblem converges, then the value function $\Snum$ converges pointwise and any cluster point of the numerical optimal trajectory $\gmnum$ is an optimal trajectory in~\eqref{eqt: result_optctrl2_hd}. Note that the convergence of the numerical optimal trajectories is not guaranteed, due to the non-uniqueness of the optimal trajectory $\gmanaly$.
Moreover, we prove error bounds for the output solution and output trajectory of Algorithm~\ref{alg:admm_minplus}. 
With this error analysis, the convergence rate of any convergent subsequence is determined by the convergence rate of each of the $m$ subproblems.
For the special case when each subproblem is solved using Algorithm~\ref{alg:admm_ver2}, the convergence rate is given by the slowest convergence rate of each subproblem, which is discussed in Section~\ref{sec:numerical_convex}.

\begin{proposition}\label{prop:convergence_ADMM_nonconvexJ}
Let $\initcond\colon\Rn\to\R$ be a function satisfying~\eqref{eqt: J_minplus} for some convex functions $\initcond_1,\dots,\initcond_m\colon\Rn\to\R$ and $\ba,\bb$ be two vectors in $(0,+\infty)^n$. Let $\bx$ be any vector in $\Rn$ and $t>0$ be any scalar. Let $\Sanaly$ be the function defined in~\eqref{eqt: result_Lax2_hd} with initial condition $\initcond$. Denote the output solution and trajectory of Algorithm~\ref{alg:admm_minplus} by $\Snum$ and $\gmnum$, respectively.
For the $j$-th subproblem, denote the analytical solution and the numerical solution by $\Sanaly_j$ and $\Snum_j$, respectively, and denote the analytical optimal trajectory and the numerical trajectory by $\gmanaly_j$ and $\gmnum_j$, respectively.
Assume there exists $\epsilon>0$, such that there holds
\begin{equation}
    |\Sanaly_j(\bx,t) - \Snum_j(\bx,t)| \leq \epsilon \quad \forall j\in\{1,\dots,m\}.
\end{equation}
Then, we have
\begin{equation}\label{eqt:prop_err_minplus_S}
    |\Sanaly(\bx,t) - \Snum(\bx,t)| \leq \epsilon.
\end{equation}
Further assume $\Sanaly_j(\bx,t) > \Sanaly(\bx,t) + 2\epsilon$ for each index $j$ satisfying $\Sanaly_j(\bx,t) \neq \Sanaly(\bx,t)$. Then, there exists an optimal trajectory $\gmanaly$ of the optimal control problem~\eqref{eqt: result_optctrl2_hd} with initial cost $\initcond$, such that there holds
\begin{equation}\label{eqt:prop_err_minplus_traj}
     \sup_{s\in[0,t]}\|\gmnum(s;\bx,t) - \gmanaly(s;\bx,t)\| \leq \sup_{s\in[0,t]}\|\gmnum_r(s;\bx,t) - \gmanaly_r(s;\bx,t)\|,
\end{equation}
where $r$ is the index in~\eqref{eqt:minplus_defr}.
\end{proposition}
\begin{proof}
The proof is provided in Appendix~\ref{subsec:appendix_pf_prop32}.
\end{proof}

Now, we present a high-dimensional numerical example using nonconvex initial cost $\initcond$ defined by:
\begin{equation}\label{eqt: hd_initial_data_minofquad}
    \initcond(\bx) = \min_{j\in\{1,2,3\}} \initcond_j(\bx) = \min_{j\in\{1,2,3\}} \left\{\frac{1}{2}\|\bx - \by_j\|^2 + \alpha_j\right\},
\end{equation}
where $\by_1 = (-2,0,\dots,0)$, $\by_2 = (2,-2,-1,0, \dots, 0)$, and $\by_3 = (0, 2, 0, \dots, 0)$ are vectors in $\R^n$ and $\alpha_1 = -0.5$, $\alpha_2 = 0$, and $\alpha_3 = -1$ are scalars in $\R$. Recall that $\|\cdot\|$ denotes the $\ell^2$-norm in the Euclidean space $\Rn$. We also use the parameters $\ba,\bb$ defined in~\eqref{eqt:numerical_def_ab}, i.e., we set $\ba = ( 4, 6, 5, \dots, 5)$ and $\bb = ( 3, 9, 6, \dots, 6)$.
Note that each subproblem has quadratic initial cost and thus can be solved using the solver in Section~\ref{sec:quad_convex}. In other words, the solver proposed in Section~\ref{sec:quad_convex} serves as the building block for Algorithm~\ref{alg:admm_minplus} in this example. Recall that the solver in Section~\ref{sec:quad_convex} has complexity $\Theta(n)$. Therefore, Algorithm~\ref{alg:admm_minplus} has complexity $\Theta(mn)$ in this example, and hence, it overcomes the curse of dimensionality.

\begin{figure}[htbp]
    \centering
    \begin{subfigure}{0.45\textwidth}
        \centering 
        \includegraphics[width=\textwidth]{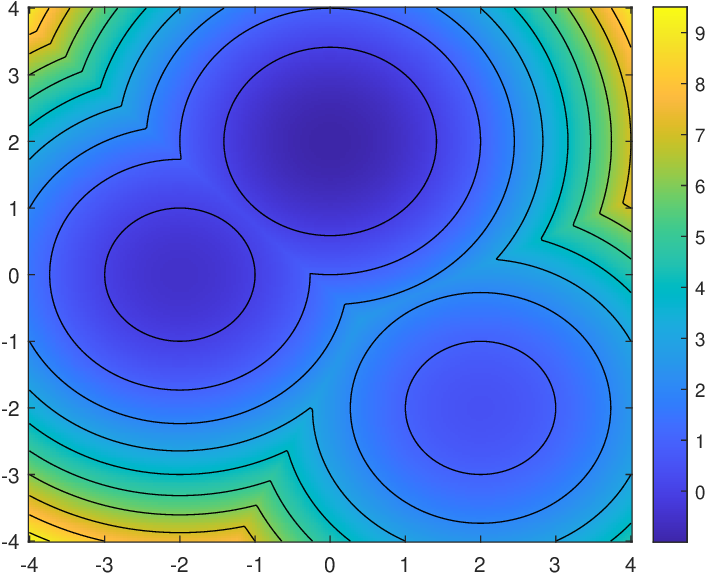}
        \caption{$t=0$}
    \end{subfigure}
    \hfill
    \begin{subfigure}{0.45\textwidth}
        \centering 
        \includegraphics[width=\textwidth]{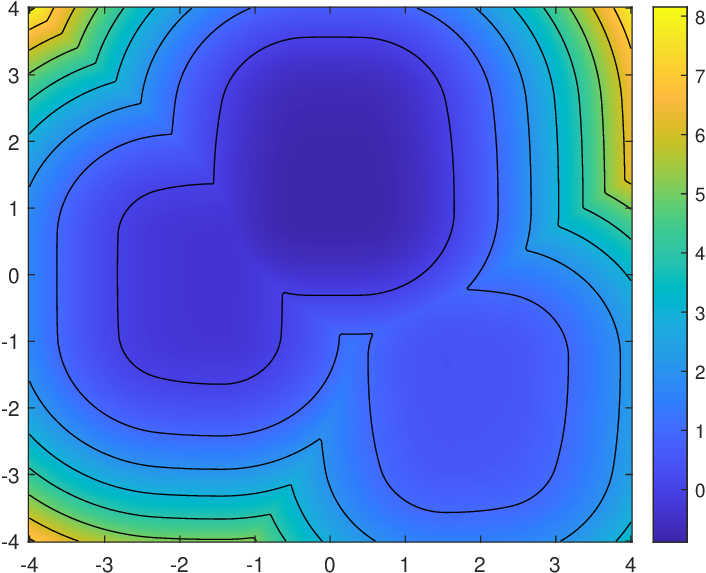}
        \caption{$t=0.1$}
    \end{subfigure}
    
    \begin{subfigure}{0.45\textwidth}
        \centering 
        \includegraphics[width=\textwidth]{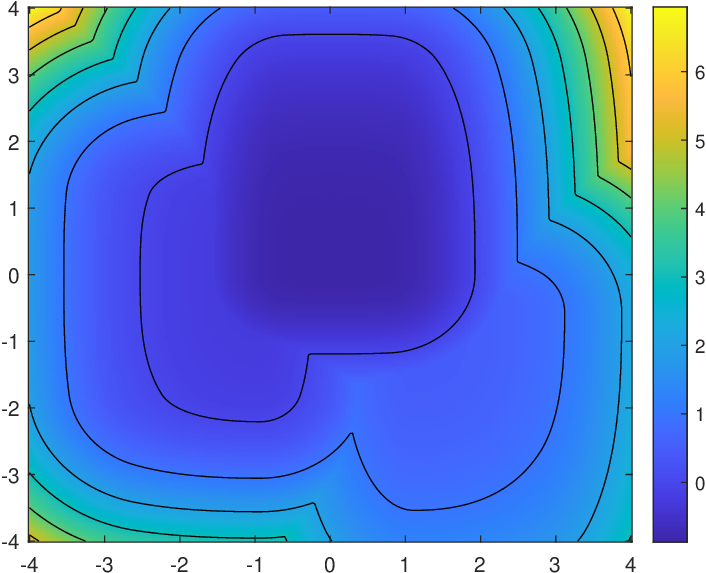}
        \caption{$t=0.2$}
    \end{subfigure}
    \hfill
    \begin{subfigure}{0.45\textwidth}
        \centering 
        \includegraphics[width=\textwidth]{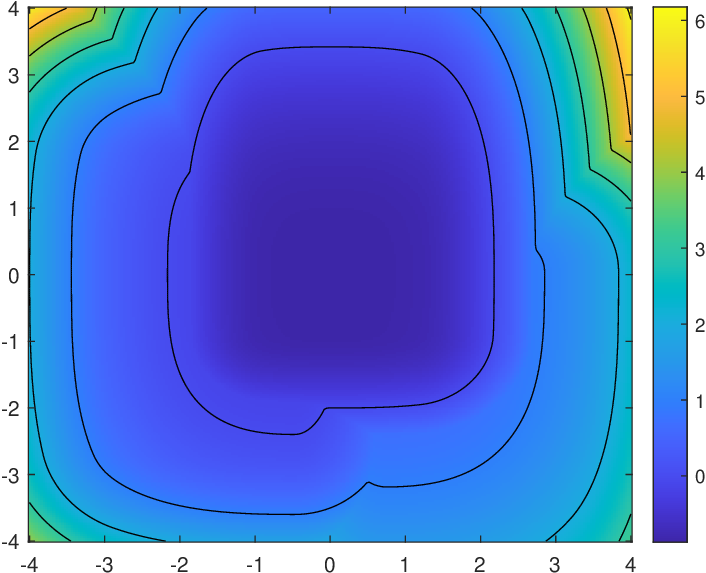}
        \caption{$t=0.3$}
    \end{subfigure}
    
    \begin{subfigure}{0.45\textwidth}
        \centering 
        \includegraphics[width=\textwidth]{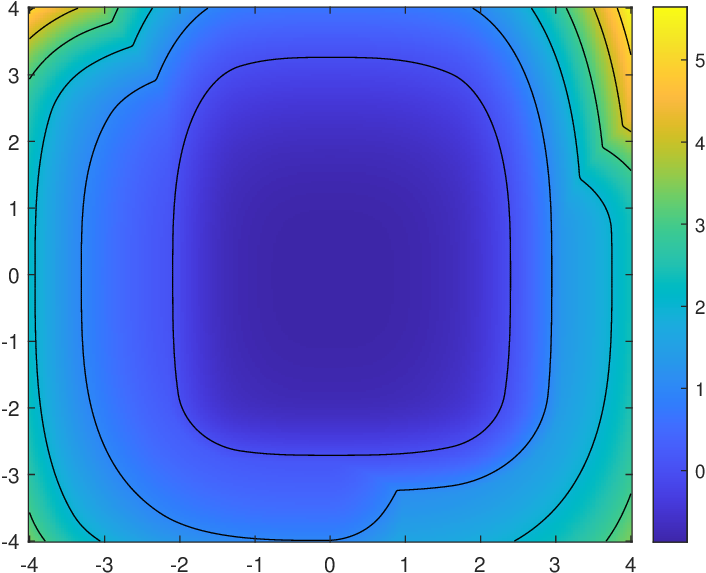}
        \caption{$t=0.4$}
    \end{subfigure}
    \hfill
    \begin{subfigure}{0.45\textwidth}
        \centering 
        \includegraphics[width=\textwidth]{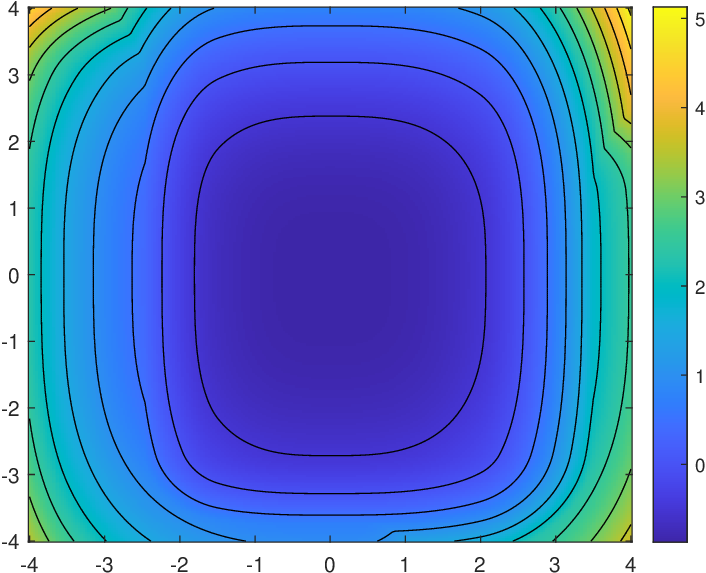}
        \caption{$t=0.5$}
    \end{subfigure}
    
    \caption{Evaluation of the solution $\Snum(\bx,t)$ of the 10-dimensional HJ PDE (\ref{eqt: result_HJ2_hd}) with $\ba = ( 4, 6, 5, \dots, 5)$, $\bb=( 3, 9, 6, \dots, 6)$, and initial condition $\initcond(\bx) = \min_{j\in\{1,2,3\}} \initcond_j(\bx) = \min_{j\in\{1,2,3\}} \{\frac{1}{2}\|\bx - \by_j\|^2 + \alpha_j\}$, where $\by_1 = (-2,0,\dots,0)$, $\by_2 = (2,-2,-1,0, \dots, 0)$, $\by_3 = (0, 2, 0, \dots, 0)$, $\alpha_1 = -0.5$, $\alpha_2 = 0$, and $\alpha_3 = -1$, for $\bx = (x_1, x_2, 0, \dots, 0)$, where $(x_1,x_2)\in [-4,4]^2$, and different times $t$. Plots for $t = 0$, $0.1$, $0.2$, $0.3$, $0.4$, and $0.5$ are depicted in (a)-(f), respectively. Level lines are superimposed on the plots.}
    \label{fig: HJ2_HD_quad}
\end{figure}

\begin{figure}[htbp]
    \centering
    \begin{subfigure}{0.32\textwidth}
        \centering \includegraphics[width=\textwidth]{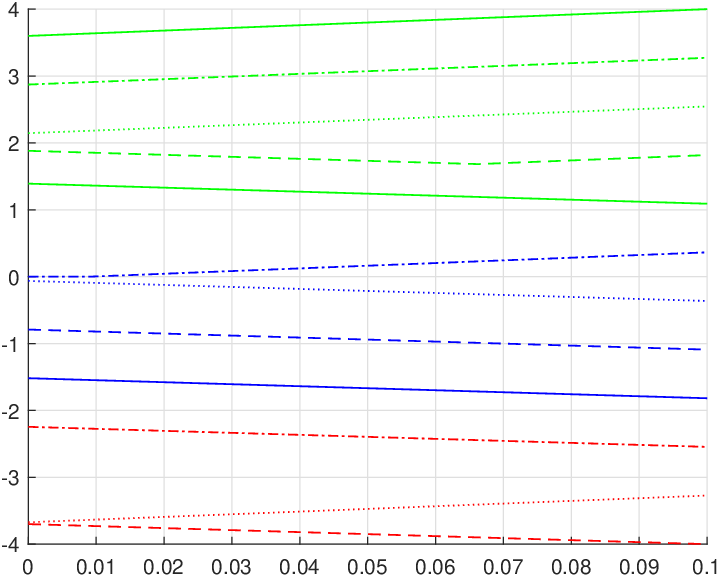}
        \caption{First component, $t=0.1$}
    \end{subfigure}
    \hfill
    \begin{subfigure}{0.32\textwidth}
        \centering \includegraphics[width=\textwidth]{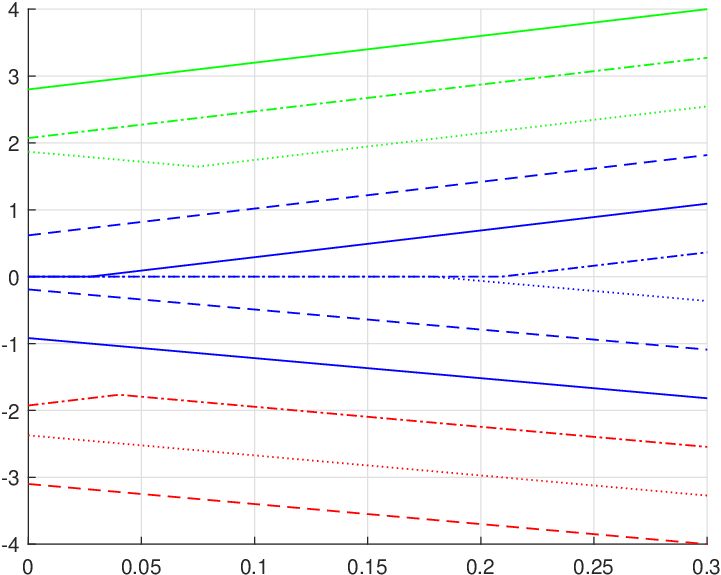}
        \caption{First component, $t=0.3$}
    \end{subfigure}
    \hfill
    \begin{subfigure}{0.32\textwidth}
        \centering \includegraphics[width=\textwidth]{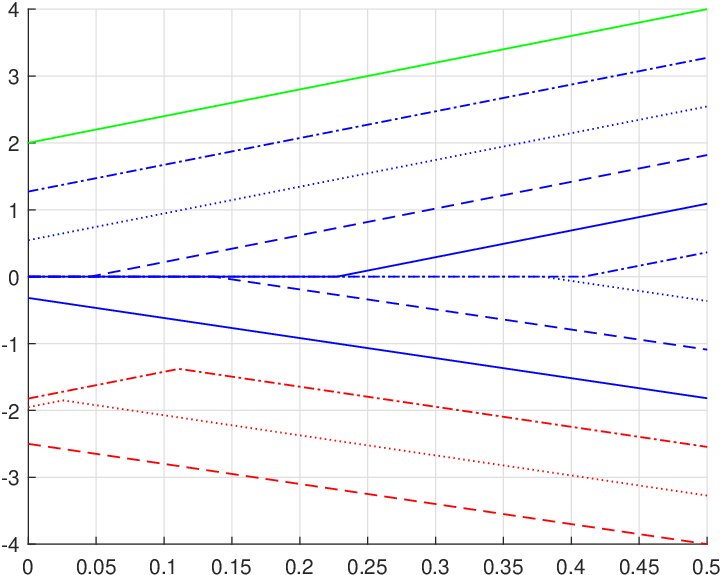}
        \caption{First component, $t=0.5$}
    \end{subfigure}
    
    \begin{subfigure}{0.32\textwidth}
        \centering \includegraphics[width=\textwidth]{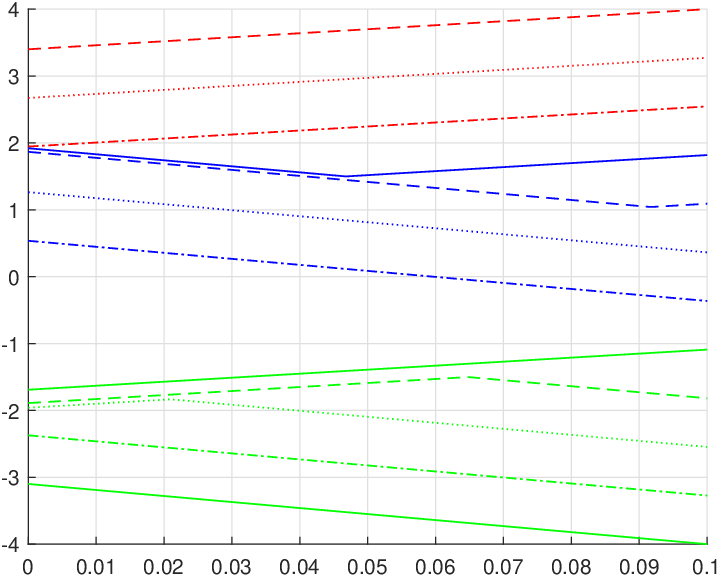}
        \caption{Second component, $t=0.1$}
    \end{subfigure}
    \hfill
    \begin{subfigure}{0.32\textwidth}
        \centering \includegraphics[width=\textwidth]{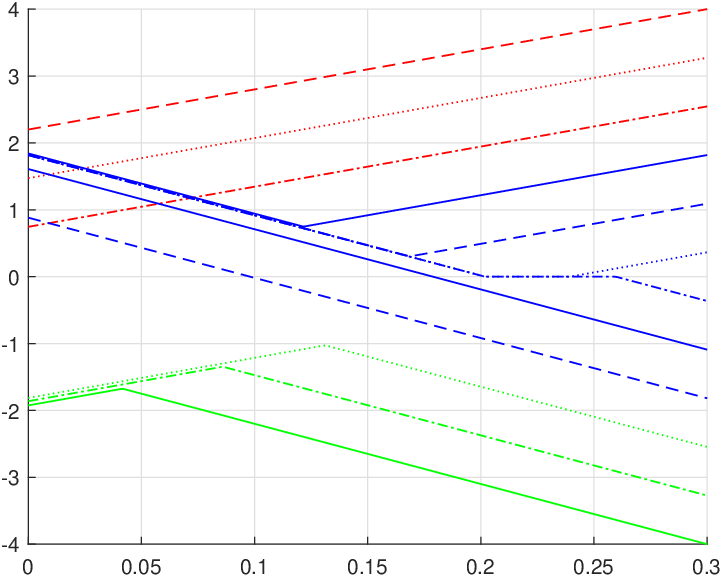}
        \caption{Second component, $t=0.3$}
    \end{subfigure}
    \hfill
    \begin{subfigure}{0.32\textwidth}
        \centering \includegraphics[width=\textwidth]{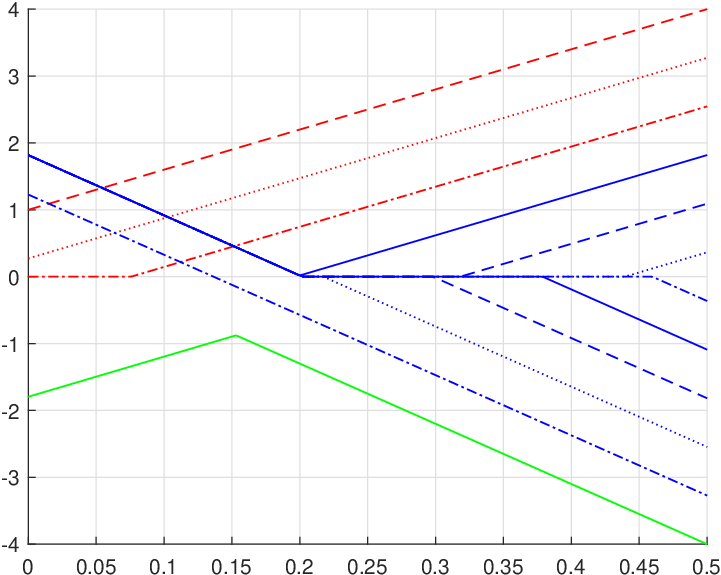}
        \caption{Second component, $t=0.5$}
    \end{subfigure}
    
    \caption{Evaluation of an optimal trajectory $\gmnum(s;(x,-x,0,\dots,0), t)$ of the $10$-dimensional optimal control problem (\ref{eqt: result_optctrl2_hd}) with $\ba = ( 4, 6, 5, \dots, 5)$, $\bb = ( 3, 9, 6, \dots, 6)$, and initial cost $\initcond(\bx) = \min_{j\in\{1,2,3\}} \initcond_j(\bx) = \min_{j\in\{1,2,3\}} \{\frac{1}{2}\|\bx - \by_j\|^2 + \alpha_j\}$, where $\by_1 = (-2,0,\dots,0)$, $\by_2 = (2,-2,-1,0, \dots, 0)$, $\by_3 = (0, 2, 0, \dots, 0)$, $\alpha_1 = -0.5$, $\alpha_2 = 0$, and $\alpha_3 = -1$, versus $s\in[0,t]$ for different terminal positions $(x,-x,0,\dots,0)$ ($x\in[-4,4]$) and different time horizons $t$. The color of the lines denotes which initial cost was used, i.e., $r \in \argmin_{j\in\{1,2,3\}} \Snum_j(\bx,t)$ for $r = 1,2,3$ corresponds to {\textcolor{red}{red}}, {\textcolor{green}{green}}, and {\textcolor{blue}{blue}}, respectively. The different line markers simply help to differentiate between the different trajectories. Figures (a)-(c) depict the first component of the trajectory versus $s\in[0,t]$ with different time horizons $t$, while Figures (d)-(f) depict the second component of the trajectory versus $s\in[0,t]$ with different time horizons $t$. Plots for time horizons $t = 0.1$, $0.3$, and $0.5$ are depicted in Figures (a)/(d), (b)/(e), and (c)/(f), respectively. We note that because of our choice of $\ba = ( 4, 6, 5, \dots, 5)$ and $\bb = ( 3, 9, 6, \dots, 6)$, the piecewise slopes of our trajectories are not symmetric about $0$.}
    \label{fig: HJ2_HD_quad_trajectory}
\end{figure}

We solve the $10$-dimensional problem (i.e., we set $n=10$) and plot the solution $\Snum$ and the optimal trajectories $\gmnum$ in Figures~\ref{fig: HJ2_HD_quad} and~\ref{fig: HJ2_HD_quad_trajectory}, respectively.
Figure \ref{fig: HJ2_HD_quad} depicts two-dimensional slices of the numerical solution $\Snum(\bx,t)$, as computed using Algorithm~\ref{alg:admm_minplus}, to HJ PDE~\eqref{eqt: result_HJ2_hd} for different positions $\bx=(x_1, x_2, 0, \dots, 0)$ and different times $t$. In Figure \ref{fig: HJ2_HD_quad}(a), we can clearly see that the initial condition $\initcond$ is not smooth at the interfaces of the quadratics $\initcond_i$; e.g., there are obvious kinks near $(x_1,x_2) = (0,0)$, $(-2,2)$, $(0,-2)$, $(2.5,0)$. In Figures~\ref{fig: HJ2_HD_quad}(b)-(f), we see that, over time, the solution also evolves with several kinks. These kinks provide numerical validation that our algorithm does indeed provide the non-smooth viscosity solution to the corresponding HJ PDE. Overall, the solution appears to be continuous (but not necessarily differentiable) in $(x_1,x_2)$ at all times $t$, which is consistent with the results of Proposition \ref{prop: HJ_hd}.

Figure \ref{fig: HJ2_HD_quad_trajectory} depicts one-dimensional slices of an optimal trajectory 
$\gmnum(s;\bx, t)$ of the corresponding $10$-dimensional optimal control problem for different terminal positions $\bx=(x,-x,0,\dots,0)$ and different time horizons~$t$. We observe that the one-dimensional slices are piecewise linear and continuous in $s$, which is consistent with the properties of the formulas in~\eqref{eqt: optctrl2_defx_1},~\eqref{eqt: optctrl2_defx_2},~\eqref{eqt: optctrl2_defx_3}, and~\eqref{eqt: optctrl2_defx_neg}. We also note that in each subplot, all line segments with positive slope are parallel with slope $a_i$ (i.e. the $i$-th component of the trajectory has velocity $a_i$), while all lines segments with negative slope are parallel with slope $-b_i$ (i.e. the $i$-th component of the trajectory has velocity $-b_i$). As such, in any given subplot, the piecewise slopes of the trajectories are not symmetric about $0$ due to our choice of $\ba$ and $\bb$. We also observe different patterns of trajectories depending on the terminal positions and time horizons.

\begin{table}[htbp]
    \centering
    \begin{tabular}{c|c|c|c}
        \hline
        \textbf{$\mathbf{n}$} & \textbf{CPU time (s)} & \textbf{FPGA time (s)} & \textbf{Speedup}  \\
        \hline
        4 & 1.7887e-07 & 1.334e-08 & 13.4085 \\
        8 & 4.5562e-07 & 2.667e-08 & 17.0836 \\
        12 & 1.3138e-06 & 4.001e-08 & 32.8370 \\
        16 & 2.1028e-06 & 5.334e-08 & 39.4226 \\
        \hline
    \end{tabular}
    \hfill 
    
    \caption{Comparison of the average time per call over $100,000$ runs for evaluating the solution of the HJ PDE~\eqref{eqt: result_HJ2_hd}  with initial condition $\initcond(\bx) = \min_{j\in\{1,2,3\}} \initcond_j(\bx) = \min_{j\in\{1,2,3\}} \{\frac{1}{2}\|\bx - \by_j\|^2 + \alpha_j\}$, where $\by_1 = (-2,0,\dots,0)$, $\by_2 = (2,-2,-1,0, \dots, 0)$, $\by_3 = (0, 2, 0, \dots, 0)$, $\alpha_1 = -0.5$, $\alpha_2 = 0$, and $\alpha_3 = -1$, for various dimensions $n$ using a CPU implementation on a single Intel Core i7-1165G7 versus an FPGA implementation on a Xilinx Alveo U280 board with a frequency of 300 MHz.}
    \label{tab:timing_minplus_quad}
\end{table}

\begin{table}[htbp]
    \centering
    \begin{tabular}{c|c|c|c|c|c}
        \hline
        \textbf{$\mathbf{n}$} & \textbf{Latency (ns)} & \textbf{BRAM} & \textbf{DSPs} & \textbf{FFs} & \textbf{LUTs} \\
        \hline
        4 & 400,244 (1.334e06) & 0 (0\%) & 2,541 (28\%) & 278,904 (10\%) & 168,376 (12\%) \\
        8 & 800,264 (2.667e06) & 0 (0\%) & 2,541 (28\%) & 279,557 (10\%) & 169,091 (12\%) \\
        12 & 1,200,272 (4.001e06) & 0 (0\%) & 2,541 (28\%) & 280,337 (10\%) & 169,759 (13\%) \\
        16 & 1,600,288 (5.334e06) & 0 (0\%) & 2,541 (28\%) & 280,548 (10\%) & 170,441 (13\%) \\
        \hline
    \end{tabular}
    \hfill 
    
    \caption{FPGA resources and latencies in cycles and nanoseconds (ns) for evaluating the solution of the HJ PDE \eqref{eqt: result_HJ2_hd} with initial condition $\initcond(\bx) = \min_{j\in\{1,2,3\}} \initcond_j(\bx) = \min_{j\in\{1,2,3\}} \{\frac{1}{2}\|\bx - \by_j\|^2 + \alpha_j\}$, where $\by_1 = (-2,0,\dots,0)$, $\by_2 = (2,-2,-1,0, \dots, 0)$, $\by_3 = (0, 2, 0, \dots, 0)$, $\alpha_1 = -0.5$, $\alpha_2 = 0$, and $\alpha_3 = -1$, at $100,000$ points $(\bx,t,\by)\in\Rn\times[0,\infty)\times\Rn$ for various dimensions $n$ using double precision floating points on a Xilinx Alveo U280 board with a frequency of 300 MHz.}
    \label{tab:minplus_quad_fpga_resources}
\end{table}

In Table~\ref{tab:timing_minplus_quad}, we show the running time of a CPU and an FPGA implementation of Algorithm~\ref{alg:admm_minplus} for this example for different dimensions $n$. We measure the running time using the same method described in Section~\ref{sec:quad_convex}, and we use the solver and corresponding implementations from Section \ref{sec:quad_convex} in line 2 of Algorithm~\ref{alg:admm_minplus} to solve each of the $j$-th subproblems, $j = 1, 2, 3$. 
For this example, we design a high throughput FPGA implementation with an II of 1 and that streams the points $(\bx, t)\in\Rn\times[0,\infty)$ elementwise (i.e., our FPGA kernel takes input $(x_i, t)\in\R\times[0,\infty)$). Specifically, our FPGA implementation essentially replicates the high-throughput implementation of the numerical solver from Section \ref{sec:quad_convex} three times (once per $j$-th subproblem), where the three copies of the building block run in parallel. The outputs of the three building blocks are then combined using~\eqref{eqt:minplus_defr}. 

As such, in Table~\ref{tab:timing_minplus_quad}, we observe that the average runtimes for our FPGA implementation for this example are nearly identical to those for the FPGA implementation of the building block in Table \ref{tab:timing_quadratic}. In contrast, the CPU implementation now has to parallelize computations for the $m=3$ sub-problems in addition to its elementwise (i.e., for each $(x_i, t)\in\R\times[0,\infty)$, $i = 1, \dots, n)$ and pointwise (i.e., for various points $(\bx, t)\in\Rn\times[0,\infty)$) parallelizations. As a result, the average CPU runtimes in Table \ref{tab:timing_minplus_quad} are almost three times slower than those for the CPU implementation of the building block in Table \ref{tab:timing_quadratic}.

Now comparing the CPU and FPGA timing results for this example, in Table~\ref{tab:timing_minplus_quad}, we see that the CPU implementation takes less than $3\times 10^{-6}$ seconds on average to compute the solution at one point for the 16-dimensional problem, which shows the efficiency of our proposed algorithm even in high dimensions. Meanwhile, our FPGA implementation takes less than $6\times 10^{-8}$ seconds on average to compute the solution at one point for the 16-dimensional problem, for a speedup of about $40$ compared to the CPU. These results highlight the promising performance boosts FPGAs are able to achieve over CPUs. 

Table \ref{tab:minplus_quad_fpga_resources} shows the FPGA resources and latencies used by our FPGA implementation for this example for different dimensions $n$. We observe that, due to our use of elementwise streaming, the latency of our FPGA implementation scales linearly in the dimension $n$, while the amount of FPGA resources used remains essentially constant in $n$. We also observe that the latencies in Table  \ref{tab:minplus_quad_fpga_resources} are nearly identical to those for the FPGA implementation of the building block in Table \ref{tab:quad_fpga_resources}. Meanwhile, the FPGA resources used in Table \ref{tab:minplus_quad_fpga_resources} are approximately 3 times greater than those used by the FPGA implementation of the building block in Table \ref{tab:quad_fpga_resources}. These results are consistent with the fact that our FPGA implementation of the algorithm for this example consists of three copies of the building block from Section \ref{sec:quad_convex} running in parallel.

In Table \ref{tab:minplus_quad_fpga_resources}, we also see that our design uses less than 30\% of the FPGA resources available on the Xilinx Alveo U280 board. This means that we could either use a smaller (i.e., cheaper) FPGA to implement our numerical solver with similar performance as we report here or we could parallelize by simply implementing multiple, independent copies of our FPGA kernel to maximize usage of the FPGA board. In the latter case, we could achieve a further speedup of $\times 3$ (i.e., 1 copy of our FPGA kernel per each of the 3 chiplets, ensuring that no kernel requires crossing chiplets) for a total speedup of about 40 to 118 over the CPU depending on the dimension $n$.

\section{Conclusion}\label{sec:conclusion}
In this paper, we present analytical solutions to certain classes of control-constrained optimal control problems and the corresponding HJ PDEs where the associated running cost and Hamiltonian have state-dependence. Moreover, we provide efficient numerical methods for these problems and describe both CPU and FPGA implementations for these methods. While our CPU implementations already demonstrate the efficiency of our solvers in high dimensions, our FPGA implementations demonstrate the additional performance boosts and benefits that FPGAs can achieve over CPUs.
Our numerical results provide several examples for which our numerical algorithms overcome the curse of dimensionality and demonstrate that our algorithms have potential for real-time high-dimensional optimal control applications.
An interesting future research direction would be to combine our algorithms with other methods, such as numerical algorithms involving the LQR solver and/or more complicated state or control constraints (see, for instance,~\cite{dower2019game}), to address broader classes of optimal control problems.

\section*{Acknowledgements}
This research is supported by 
 DOE-MMICS SEA-CROGS DE-SC0023191
and AFOSR MURI FA9550-20-1-0358. P.C. is supported by the SMART Scholarship, which is funded by USD/R\&E (The Under Secretary of Defense-Research and Engineering), National Defense Education Program (NDEP) / BA-1, Basic Research. We thank Peter Dower for his useful feedback.

\section*{Statements and Declarations}
The authors declare that they have no known competing financial interests or personal relationships that could have influenced or appeared to have influenced the work reported in this paper. Furthermore, the authors declare that they have no known conflicts of interest.

\bibliographystyle{spmpsci}      
\bibliography{biblist_new}

\appendix
\section{Some technical lemmas for the analytical solutions}\label{sec: appendix_analyticsolns}

\begin{lemma}\label{lem: appendix_optval2_equalS}
Let $a,b,t$ be positive scalars and $x,\initpos$ be real numbers satisfying $\initpos-bt\leq x\leq \initpos+at$. 
Let $\valuefn$ be the function defined in~\eqref{eqt: result_S2_1d} and~\eqref{eqt: result_S2_1d_negative} and $[0,t]\ni s\mapsto \opttraj(s;x,t,\initpos,a,b)\in\R$ be the trajectory defined in~\eqref{eqt: optctrl2_defx_1},~\eqref{eqt: optctrl2_defx_2},~\eqref{eqt: optctrl2_defx_3}, and~\eqref{eqt: optctrl2_defx_neg} for different cases. Then, there holds
\begin{equation}\label{eqt: lemB1_optval_equals_S}
    \int_0^t \frac{1}{2} \left(\opttraj(s;x,t,\initpos,a,b)\right)^2 ds = \valuefn(x,t; \initpos,a,b).
\end{equation}
\end{lemma}

\begin{proof}
If $(x,t,\initpos)\in \Omega_1$ holds, we have
\begin{equation*}
\begin{split}
    \int_0^t \frac{1}{2} \left(\opttraj(s;x,t,\initpos,a,b)\right)^2 ds
    &= \int_0^{\frac{-x+\initpos + at}{a+b}} \frac{1}{2}(\initpos - bs)^2 ds + \int_{\frac{-x+\initpos + at}{a+b}}^t \frac{1}{2}(as - at+x)^2 ds\\
    &= -\left.\frac{1}{6b}(\initpos-bs)^3\right|_0^{\frac{-x+\initpos + at}{a+b}} + \left.\frac{1}{6a}(as-at+x)^3\right|_{\frac{-x+\initpos + at}{a+b}}^t\\
    &= -\left(\frac{1}{6b}+\frac{1}{6a}\right)\left(\frac{a\initpos+bx-abt}{a+b}\right)^3 + \frac{\initpos^3}{6b} + \frac{x^3}{6a}\\
    &= \valuefn(x,t;\initpos,a,b).
\end{split}
\end{equation*}
If $(x,t,\initpos)\in\Omega_2$ holds, we have
\begin{equation*}
\begin{split}
    \int_0^t \frac{1}{2} \left(\opttraj(s;x,t,\initpos,a,b)\right)^2 ds
    &= \int_0^{\frac{\initpos}{b}} \frac{1}{2}(\initpos - bs)^2 ds + \int_{t-\frac{x}{a}}^t \frac{1}{2}(as - at+x)^2 ds\\
    &= -\left.\frac{1}{6b}(\initpos-bs)^3\right|_0^{\frac{\initpos}{b}} + \left.\frac{1}{6a}(as-at+x)^3\right|_{t-\frac{x}{a}}^t\\
    &= \frac{\initpos^3}{6b} + \frac{x^3}{6a}\\
    &= \valuefn(x,t;\initpos,a,b).
\end{split}
\end{equation*}
If $(x,t,\initpos)\in \Omega_3$ holds, we have
\begin{equation*}
\begin{split}
    \int_0^t \frac{1}{2} \left(\opttraj(s;x,t,\initpos,a,b)\right)^2 ds
    &= \int_0^{\frac{\initpos}{b}} \frac{1}{2}(\initpos - bs)^2 ds + \int_{t-\frac{-x}{b}}^t \frac{1}{2}(-bs +bt+x)^2 ds\\
    &= -\left.\frac{1}{6b}(\initpos-bs)^3\right|_0^{\frac{\initpos}{b}} - \left.\frac{1}{6b}(-bs+bt+x)^3\right|_{t-\frac{-x}{b}}^t\\
    &= \frac{\initpos^3}{6b} - \frac{x^3}{6b}\\
    &= \valuefn(x,t;\initpos,a,b).
\end{split}
\end{equation*}
If $\initpos <0$, we have 
\begin{equation*}
\begin{split}
    \int_0^t \frac{1}{2} \left(\opttraj(s;x,t,\initpos,a,b)\right)^2 ds
    &= \int_0^t \frac{1}{2} \left(\opttraj(s;-x,t,-\initpos,b,a)\right)^2 ds\\
    &= \valuefn(-x,t;-\initpos,b,a)
    = \valuefn(x,t;\initpos,a,b).
\end{split}
\end{equation*}
Therefore,~\eqref{eqt: lemB1_optval_equals_S} holds for any $(x,t,\initpos)\in\R\times(0,+\infty)\times \R$ satisfying $\initpos-bt\leq x\leq \initpos+at$.
\end{proof}

\begin{lemma}\label{lem: appendix_conv_S_x0}
Let $a,b$ be positive scalars. 
Let $\valuefn$ be the function defined in~\eqref{eqt: result_S2_1d} and~\eqref{eqt: result_S2_1d_negative}. Then, for any $x\in\R$, $t>0$, the function $\R\ni \initpos\mapsto \valuefn(x,t;\initpos,a,b)\in \R\cup\{+\infty\}$ is strictly convex and twice continuously differentiable in its domain. 
\end{lemma}
\begin{proof}
In this proof, we regard the function $\valuefn(x,t;\initpos,a,b)$ as a function of $\initpos$ from its domain $[x-at,x+bt]$ to $\R$, and we use its derivative to mean the derivative of $\valuefn$ with respect to $\initpos$, by default.
To prove the statement, we need to prove that $\valuefn$ is twice continuously differentiable and that the second-order derivative is positive almost everywhere in the domain. 
We consider the following cases.

First, assume $x\geq at$ holds. After some computation, the function $\initpos\mapsto \valuefn(x,t;\initpos,a,b)$ can be written as follows:
\begin{equation*}
    \valuefn(x,t;\initpos,a,b) = \frac{\initpos^3}{6b} + \frac{x^3}{6a} - \left(\frac{1}{6a} + \frac{1}{6b}\right)\left(\frac{a\initpos + bx - abt}{a+b}\right)^3 \quad \forall \initpos\in [x-at, x+bt],
\end{equation*}
which is twice continuously differentiable. The second-order derivative is given by
\begin{equation}\label{eqt:lemA2_pf_2der_om1}
\begin{split}
    \frac{\partial^2 \valuefn(x,t;\initpos,a,b)}{\partial \initpos^2} &= \frac{\initpos}{b} - \frac{a}{b(a+b)^2} (a\initpos + bx-abt) = \frac{(b^2+2ab)\initpos - ab(x -at)}{b(a+b)^2}\\
    &\geq \frac{(b^2+ab)\initpos }{b(a+b)^2}\geq 0,
\end{split}
\end{equation}
where the first and second inequalities hold since we have $\initpos\geq x-at\geq 0$. Moreover, the second inequality becomes equality if and only if $\initpos$ is zero. In other words, the second-order derivative in~\eqref{eqt:lemA2_pf_2der_om1} is positive almost everywhere, and hence, the conclusion holds in this case.

Next, assume that $x$ is a point in $[0,at)$. In this case, the function $\initpos\mapsto \valuefn(x,t;\initpos,a,b)$ can be written as follows:
\begin{equation*}
    \valuefn(x,t;\initpos,a,b) = \begin{dcases}
    -\frac{\initpos^3}{6a} + \frac{x^3}{6a} & x-at\leq \initpos <0,\\
    \frac{\initpos^3}{6b} + \frac{x^3}{6a} & 0\leq \initpos <bt-\frac{bx}{a},\\
    \frac{\initpos^3}{6b} + \frac{x^3}{6a} - \left(\frac{1}{6a} + \frac{1}{6b}\right)\left(\frac{a\initpos + bx - abt}{a+b}\right)^3 & bt-\frac{bx}{a}\leq \initpos\leq x+bt.
    \end{dcases}
\end{equation*}
It is straightforward to check that this function is twice continuously differentiable and that the second-order derivative reads:
\begin{equation}\label{eqt:lemA2_pf_2der_om2}
    \frac{\partial^2 \valuefn(x,t;\initpos,a,b)}{\partial \initpos^2} =
    \begin{dcases}
    -\frac{\initpos}{a} & x-at< \initpos< 0,\\
    \frac{\initpos}{b} & 0\leq  \initpos < bt-\frac{bx}{a},\\
    \frac{(b^2+2ab)\initpos - ab(x -at)}{b(a+b)^2}
    & bt-\frac{bx}{a}\leq  \initpos< x+bt,
    \end{dcases}
\end{equation}
where the first line is positive since $\initpos<0$ holds in the first line, the second line is positive almost everywhere since $\initpos >0$ holds almost everywhere in the second line, and the third line is positive since the inequalities in~\eqref{eqt:lemA2_pf_2der_om1} also hold according to the condition on $u$ (there holds $u \geq bt-\frac{bx}{a}>0>x-at$). Therefore, the conclusion follows in this case.

Finally, we consider the case when $x<0$. By definition, we have that $\valuefn(x,t;\initpos,a,b) = \valuefn(-x,t;-\initpos,b,a)$, where the right-hand side is twice continuously differentiable and whose second-order derivative with respect to $-\initpos$ is positive almost everywhere by the same argument above. Therefore, the function $\valuefn(x,t;\initpos,a,b)$ is also strictly convex and twice continuously differentiable with respect to $\initpos$, and the conclusion holds.
\end{proof}

\section{Some computations for the numerical implementation}\label{sec: ADMMcomp_HJ2}

\subsection{A numerical method for computing the proximal point of $u\mapsto \frac{1}{\lambda}\valuefn(x,t;u,a,b)$} \label{sec:appendix_alg_proxV}
Here, we discuss how to compute the proximal point of the function $\R\ni u\mapsto \frac{1}{\lambda}\valuefn(x,t;u,a,b)\in \R\cup\{+\infty\}$, i.e., how to solve the following convex optimization problem:
\begin{equation}\label{eqt: ADMM_HJ2_di_update_simplified}
    u^* = \argmin_{u\in\mathbb{R}} \left\{\valuefn(x,t; u,a,b) + \frac{\lambda}{2}(u - \proxpty)^2\right\} = \argmin_{u\in [x-at,x+bt]} \left\{\valuefn(x,t; u,a,b) + \frac{\lambda}{2}(u - \proxpty)^2\right\},
\end{equation}
for any $\lambda,t,a,b>0$, and $x,\proxpty\in\R$.
We consider the following two cases for the variable $u$.

If $u \geq 0$, after some computation, the objective function in~\eqref{eqt: ADMM_HJ2_di_update_simplified} can be written as
\begin{equation*}
    \begin{split}
    F(u;x,t,a,b) 
     :=& \valuefn(x,t; u,a,b) + \frac{\lambda}{2}(u - \proxpty)^2 \\
    =& \begin{dcases}
        \frac{u^3}{6b} + \frac{x^3}{6a} - \left(\frac{1}{6a} + \frac{1}{6b}\right)\left(\frac{au + bx - abt}{a+b}\right)^3 + \frac{\lambda}{2}(u - \proxpty)^2 & u\in \Omega_1(x,t,a,b),\\
        \frac{u^3}{6b} + \frac{x^3}{6a} + \frac{\lambda}{2}(u - \proxpty)^2 & u\in \Omega_2(x,t,a,b),\\ 
        \frac{u^3}{6b} - \frac{x^3}{6b} + \frac{\lambda}{2}(u - \proxpty)^2 & u\in \Omega_3(x,t,a,b),\\
        +\infty & \text{otherwise},
    \end{dcases}    
    \end{split}
\end{equation*}
where the three regions $\Omega_1(x,t,a,b), \Omega_2(x,t,a,b), \Omega_3(x,t,a,b)\subset [0,+\infty)$ are defined by:
\begin{equation*}
\begin{split}
    \Omega_1(x,t,a,b) &:= \left\{u\in (bt,+\infty)\colon x-at\leq u\leq x+bt\right\} \bigcup \left\{u\in[0,bt]\colon u\geq x-at,\,\, u\geq bt-\frac{bx}{a}\right\},\\
    \Omega_2(x,t,a,b) &:= \begin{dcases}
    \left[0,bt-\frac{bx}{a}\right) & x\geq 0, \\
    \emptyset & x<0,
    \end{dcases}\\
    \Omega_3(x,t,a,b) &:= \begin{dcases}
    [0,x+bt] & x< 0,\\
    \emptyset & x\geq 0.
    \end{dcases}
\end{split}
\end{equation*}
In this case, the derivative of $F$ with respect to $u$ is given by:
\begin{equation}\label{eqt:appendixB_der1}
\begin{split}
    &\frac{\partial}{\partial u} F(u;x,t,a,b) \\
    =& \begin{dcases}
        \frac{(2a+b)u^2 - 2a(x-at)u - b(x-at)^2}{2(a+b)^2} + \lambda (u - \proxpty) & u\in \Omega_1(x,t,a,b),\\
        \frac{u^2}{2b} + \lambda(u - \proxpty) & u\in \Omega_2(x,t,a,b) \cup \Omega_3(x,t,a,b),
    \end{dcases}
\end{split}
\end{equation}
and the second derivative of $F$ with respect to $u$ can be easily computed using \eqref{eqt:lemA2_pf_2der_om1} and \eqref{eqt:lemA2_pf_2der_om2}, for different cases. 
To get possible candidates for the minimizer $u^*$ of $F$ in this case, we compute the roots of the functions in the two lines of~\eqref{eqt:appendixB_der1} and select the roots where the second derivative of $F$ is non-negative. After some calculations, the candidates are given by $u_1$ and $u_2$, which are defined as follows:
\begin{equation}\label{eqt:appendixB1_def_u1_u2}
\begin{split}
    u_1 & 
    := -\frac{\lambda(a+b)^2 - a(x-at)}{2a+b} + \sqrt{\left(\frac{\lambda(a+b)^2-a(x-at)}{2a+b}\right)^2 + \frac{b(x-at)^2}{2a+b} + \frac{2\lambda (a+b)^2\proxpty}{2a+b}},\\
    u_2 & := -\lambda b + \sqrt{(\lambda b)^2 + 2 \lambda b \proxpty}.
\end{split}
\end{equation}
Note that $u_1$ and $u_2$ may be not well-defined if the term under the square root is negative, in which case the corresponding function does not provide a possible candidate for $u^*$. Therefore, we assign $u_i$ ($i=1,2$) to be an arbitrary point in $[x-at,x+bt]$ if it is not well-defined.

If $u < 0$, after some computation, the objective function in~\eqref{eqt: ADMM_HJ2_di_update_simplified} can be written as
\begin{align*}
    F(u;x,t,a,b) & := \valuefn(x,t; u,a,b) + \frac{\lambda}{2}(u - \proxpty)^2 \\
    & = \valuefn(-x,t; -u,b,a) + \frac{\lambda}{2}(u - \proxpty)^2 \\
    & = \begin{dcases}
        -\frac{u^3}{6a} - \frac{x^3}{6b} + \left(\frac{1}{6a} + \frac{1}{6b}\right)\left(\frac{bu + ax + abt}{a+b}\right)^3 + \frac{\lambda}{2}(u - \proxpty)^2 & -u\in \Omega_1(-x,t,b,a),\\
        -\frac{u^3}{6a} - \frac{x^3}{6b} + \frac{\lambda}{2}(u - \proxpty)^2 & -u\in \Omega_2(-x,t,b,a),\\ 
        -\frac{u^3}{6a} + \frac{x^3}{6a} + \frac{\lambda}{2}(u - \proxpty)^2 & -u\in \Omega_3(-x,t,b,a),\\
        +\infty & \text{otherwise}.
    \end{dcases}
\end{align*}
Thus, for $u < 0$, the derivative of $F$ with respect to $u$ is given by:
\begin{equation}\label{eqt:appendixB_der2}
\begin{split}
    &\frac{\partial}{\partial u} F(u;x,t,a,b) \\
    =& \begin{dcases}
        \frac{-(a+2b)u^2 + 2b(x+bt)u + a(x+bt)^2}{2(a+b)^2} + \lambda (u - \proxpty) & -u\in \Omega_1(-x,t,b,a),\\
        -\frac{u^2}{2a} + \lambda(u - \proxpty) & -u\in \Omega_2(-x,t,b,a) \cup \Omega_3(-x,t,b,a).
    \end{dcases}
\end{split}
\end{equation}
Similarly as in the first case, we take the roots of the two functions in~\eqref{eqt:appendixB_der2}, such that the second order derivative of $F$ is non-negative. These roots provide possible candidates for $u^*$. We denote these candidates by $u_1'$ and $u_2'$, which are defined by:
\begin{equation}\label{eqt:appendixB1_def_u1p_u2p}
\begin{split}
    u_1' & 
    := \frac{\lambda(a+b)^2 + b(x+bt)}{a+2b} - \sqrt{\left(\frac{\lambda(a+b)^2 + b(x+bt)}{a+2b}\right)^2 + \frac{a(x+bt)^2}{a+2b} - \frac{2\lambda (a+b)^2\proxpty}{a+2b}}, \\
    u_2' & := \lambda a - \sqrt{(\lambda a)^2 - 2\lambda a \proxpty}.
\end{split}
\end{equation}
Similarly, if $u_1'$ or $u_2'$ is not well-defined, we set it to be any point in $[x-at,x+bt]$.

Note that the objective function $F$ is strictly convex and twice continuously differentiable with respect to $u$ by Lemma~\ref{lem: appendix_conv_S_x0}. Then, by the first and second derivative tests, the minimizer $u^*$ in~\eqref{eqt: ADMM_HJ2_di_update_simplified} is selected among the possible candidates $u_1,u_2,u'_1, u_2'$ defined in~\eqref{eqt:appendixB1_def_u1_u2} and~\eqref{eqt:appendixB1_def_u1p_u2p}, as well as the boundary points $x-at$ and $x+bt$. In other words, the minimizer $u^*$ satisfies
\begin{equation}\label{eqt: ADMM_dk_explicit_formula}
u^* =  
\argmin_{u\in\{u_1,u_2,u'_1,u'_2,x - at, x+bt\}} F(u;x,t,a,b).
\end{equation}
Numerically, we solve the optimization problem~\eqref{eqt: ADMM_HJ2_di_update_simplified} by computing the six candidates $u_1,u_2,u_1',u_2',x-at,x+bt$ and comparing the objective function values at those points. Then, the minimizer $u^*$ is selected using~\eqref{eqt: ADMM_dk_explicit_formula}. Therefore, the complexity of solving~\eqref{eqt: ADMM_HJ2_di_update_simplified} is $\Theta(1)$.

\subsection{An equivalent expression for $\valuefn(x,t;u,a,b)$ and $\opttraj(s;x,t,u,a,b)$}
Let $\valuefn$ be the function defined by~\eqref{eqt: result_S2_1d} and~\eqref{eqt: result_S2_1d_negative}, and let $\opttraj$ be the function defined 
by~\eqref{eqt: optctrl2_defx_1},~\eqref{eqt: optctrl2_defx_2}, \eqref{eqt: optctrl2_defx_3}, and~\eqref{eqt: optctrl2_defx_neg} for different cases.
Now, we present an equivalent expression for $\valuefn$ and $\opttraj$, which is used in our numerical implementation.

By straightforward calculation, the function $\valuefn$ can equivalently be expressed as
\begin{equation*}
\begin{split}
    &\valuefn(x,t;u,a,b) \\
    =& 
    \begin{dcases}
    \max\{\valuefn_3(x,t;u,a,b), \min\{\valuefn_1(x,t;u,a,b), \valuefn_2(x,t;u,a,b)\}\} & \text{if } u\in [0, x+bt],\\
    \max\{\valuefn_3(-x,t;-u,b,a), \min\{\valuefn_1(-x,t;-u,b,a), \valuefn_2(-x,t;-u,b,a)\}\} 
    & \text{if } u\in [x-at, 0),\\
    +\infty &\text{otherwise},
    \end{dcases}
\end{split}
\end{equation*}
where $\valuefn_1,\valuefn_2$, and $\valuefn_3$ are the functions in the first, second, and third lines of~\eqref{eqt: result_S2_1d}, respectively.

Similarly, assuming $u\in [x-at,x+bt]$ holds, the function $\opttraj$ can be expressed as
\begin{equation*}
    \opttraj(s;x,t,u,a,b)
    = \begin{dcases}
    \max\{u-bs, a(s-t)+x,0\} & \text{if } x\geq 0, 0\leq u\leq x+bt,\\
    \max\{u-bs,0\} + \min\{-b(s-t)+x,0\}
    & \text{if } x< 0, 0\leq u\leq x+bt,\\
    \min\{u+as, -b(s-t)+x,0\} & \text{if } x< 0, x-at\leq u<0,\\
    \min\{u+as,0\} + \max\{a(s-t)+x,0\}
    & \text{if } x\geq 0, x-at\leq u<0.
    \end{dcases}
\end{equation*}

Compared to the definitions~\eqref{eqt: result_S2_1d},~\eqref{eqt: result_S2_1d_negative},~\eqref{eqt: optctrl2_defx_1},~\eqref{eqt: optctrl2_defx_2},~\eqref{eqt: optctrl2_defx_3}, and~\eqref{eqt: optctrl2_defx_neg}, these equivalent formulas involve less conditional branching, and hence, they are more favorable for the performance of our numerical implementation.
\section{Proofs of convergence results in Section~\ref{sec: ADMM}}
In this section, we provide the proof of Proposition~\ref{prop:convergence_ADMM_convexJ} in Section~\ref{subsec:appendix_pf_prop31} and the proof of Proposition~\ref{prop:convergence_ADMM_nonconvexJ} in Section~\ref{subsec:appendix_pf_prop32}.
\subsection{Proof of Proposition~\ref{prop:convergence_ADMM_convexJ}} \label{subsec:appendix_pf_prop31}
Let $\bv^N$, $\bd^N$, and $\initposhd^N$ be the corresponding vectors in the algorithm at the $N$-th iteration. Let $\initposhd^*$ be the minimizer of the minimization problem in~\eqref{eqt: result_Lax2_hd}, which is unique since $\initcond$ is convex and 
each $\initpos_i\mapsto \valuefn(x_i,t;\initpos_i,a,b)$ is strictly convex by Lemma~\ref{lem: appendix_conv_S_x0}. 
According to~\cite[Theorem~2.2]{Deng2016global} whose assumptions are proved using~\cite[Remark~2.2]{Deng2016global}, both $\bv^N$ and $\bd^N$ converge to the point $\initposhd^*$ as $N$ approaches infinity, and hence, $\initposhd^N$ in Algorithm~\ref{alg:admm_ver2} also converges to $\initposhd^*$. 
Since $\initcond$ is a real-valued convex function, it is continuous in $\Rn$ and we have
\begin{equation}\label{eqt:prop31pf_convJ}
    \lim_{N\to\infty} \initcond(\initposhd^N) = \initcond(\initposhd^*).
\end{equation}
Note that the domain of the function $\initposhd \mapsto \sum_{i=1}^n\valuefn(x_i,t;\initpos_i, a_i,b_i)$ equals
\begin{equation}\label{eqt:prop31pf_domS}
    \prod_{i=1}^n [x_i-a_it,x_i+b_it],
\end{equation}
and the point $\initposhd^N = \bd^N$ is in the set in~\eqref{eqt:prop31pf_domS} by definition of $\bd^N$. Thus, the point $\initposhd^N$ is in the domain of the function $\initposhd \mapsto \sum_{i=1}^n\valuefn(x_i,t;\initpos_i, a_i,b_i)$.
By Lemma~\ref{lem: appendix_conv_S_x0}, the function $\initposhd \mapsto \sum_{i=1}^n\valuefn(x_i,t;\initpos_i, a_i,b_i)$ is continuous in its domain. It is also straightforward to check that the function $\initposhd \mapsto \sum_{i=1}^n\valuefn(x_i,t;\initpos_i, a_i,b_i)$ is Lipschitz in its domain, and we denote its Lipschitz constant by $L_i$. Thus, we have that
\begin{equation}\label{eqt:prop31pf_convS}
    \left|\sum_{i=1}^n\valuefn(x_i,t;\initpos^N_i, a_i,b_i)
    -\sum_{i=1}^n\valuefn(x_i,t;\initpos^*_i, a_i,b_i)\right| \leq \left(\sum_{i=1}^n L_i\right) \|\bu^N - \bu^*\|.
\end{equation}
Then, the convergence of $\Snum^N(\bx,t)$ to $\valuefn(\bx,t)$ follows from~\eqref{eqt:prop31pf_convJ} and~\eqref{eqt:prop31pf_convS}.

Now, it remains to prove the second formula in~\eqref{eqt:prop31_conv}. We have proved that $\initposhd^N$ and $\initposhd^*$ are both in the set in~\eqref{eqt:prop31pf_domS}. Let $\{\initposhd^{N_j}\}_j$ be a subsequence of $\{\initposhd^{N}\}_N$ (i.e., we assume $N_1<N_2<\cdots$ and $\lim_{j\to\infty}N_j=+\infty$), such that for each $i\in\{1,\dots, n\}$, the $i$-th component $\{\initpos_i^{N_j}\}_j$ of the subsequence satisfies one of the following assumptions:
\begin{itemize}
    \item[(i)] There exists an index $r_i$ in $\{1,2,3\}$, such that there hold
    \begin{equation*}
    \initpos_i^{N_j}\geq 0 \quad \text{ and }\quad  (x_i,t,\initpos_i^{N_j})\in \bar{\Omega}_{r_i}(a_i,b_i)\quad    \forall\, j\in\N.
    \end{equation*}
    \item[(ii)] There exists an index $r_i$ in $\{1,2,3\}$, such that there hold 
    \begin{equation*}
    \initpos_i^{N_j}< 0 \quad \text{ and }\quad  (-x_i,t,-\initpos_i^{N_j})\in \bar{\Omega}_{r_i}(b_i,a_i)\quad    \forall\, j\in\N.
    \end{equation*}
\end{itemize}
Here, to emphasize the dependence on $a_i$ and $b_i$, we use $\Omega_{r_i}(a_i,b_i)$ and $\bar{\Omega}_{r_i}(a_i,b_i)$ to respectively denote the set defined in~\eqref{eqt:def_dom_omegai} with constants $a=a_i$ and $b=b_i$ and its closure. Note that the situations considered in cases (i) and (ii) give a partition of the set in~\eqref{eqt:prop31pf_domS} (where some sets in the partition may be empty and the sets may overlap on the boundary, but neither of these possibilities affect the result). Hence, if the statement is proved for any such subsequence $\{\initposhd^{N_j}\}_j$, then the statement also holds for the whole sequence $\{\initposhd^{N}\}_N$. Thus, it suffices to prove the statement for the subsequence $\{\initposhd^{N_j}\}_j$.

Let $i\in\{1,\dots,n\}$ be any index. Assume case (i) holds for $\{\initpos_i^{N_j}\}_j$ with the index $r_i$. 
In other words, we assume $(x_i,t,\initpos_i^{N_j})\in \bar{\Omega}_{r_i}(a_i,b_i)$ holds for any $j\in\N$. Since the set $\bar{\Omega}_{r_i}(a_i,b_i)$ is closed and the subsequence $\{\initpos_i^{N_j}\}_j$ converges to $\initpos_i^*$, we conclude that $(x_i,t,\initpos_i^*)\in \bar{\Omega}_{r_i}(a_i,b_i)$ also holds.
Then, by definition of $\opttraj$ in $\Omega_{r_i}(a_i,b_i)$, 
it is straightforward to check that
\begin{equation}\label{eqt:prop31pf_opttraj_ineqt}
    \sup_{s\in[0,t]} \left|\opttraj(s;x_i,t, \initpos_i^{N_j},a_i,b_i) - \opttraj(s;x_i,t, \initpos_i^*,a_i,b_i)\right|\leq \left|\initpos_i^{N_j}-\initpos_i^*\right|.
\end{equation}
The proof for case (ii) is similar, so we omit it.
Note that~\eqref{eqt:prop31pf_opttraj_ineqt} holds for any arbitrary index $i\in\{1,\dots,n\}$. Hence, we have that
\begin{equation*}
\begin{split}
    \sup_{s\in[0,t]}\left\|\gmnum^{N_j}(s;\bx,t) - \opttrajhd(s;\bx,t)\right\|^2
    &= \sup_{s\in[0,t]}\sum_{i=1}^n\left|\opttraj(s;x_i,t, \initpos_i^{N_j},a_i,b_i) - \opttraj(s;x_i,t, \initpos_i^*,a_i,b_i)\right|^2\\
    &\leq \sum_{i=1}^n \sup_{s\in[0,t]}\left|\opttraj(s;x_i,t, \initpos_i^{N_j},a_i,b_i) - \opttraj(s;x_i,t, \initpos_i^*,a_i,b_i)\right|^2\\
    &\leq \sum_{i=1}^n\left|\initpos_i^{N_j}-\initpos_i^*\right|^2\\
    &=\left\|\initposhd^{N_j}-\initposhd^*\right\|^2,
\end{split}
\end{equation*}
where the second inequality holds by~\eqref{eqt:prop31pf_opttraj_ineqt}. Thus, the second formula in~\eqref{eqt:prop31_conv} holds for the subsequence by the convergence of $\initposhd^{N_j}$ to $\initposhd^*$. Moreover, the argument holds for any such subsequence, and hence, the statement holds for the whole sequence.
\qed

\subsection{Proof of Proposition~\ref{prop:convergence_ADMM_nonconvexJ}}
\label{subsec:appendix_pf_prop32}

Let $r$ be the index defined in~\eqref{eqt:minplus_defr}, and let the index set $\mathcal{J}\subseteq \{1,\dots,m\}$ be defined by:
\begin{equation*}
    \mathcal{J} := \argmin_{i\in\{1,\dots,m\}} \Sanaly_i(\bx,t).
\end{equation*}
Then, we have that
\begin{equation*}
    \Sanaly(\bx,t) - \Snum(\bx,t) = \Sanaly(\bx,t) - \Snum_r(\bx,t) \leq  \Sanaly_r(\bx,t) - \Snum_r(\bx,t) \leq \epsilon,
\end{equation*}
where the first equality holds by definition of $r$ and the first inequality holds since $\Sanaly$ satisfies~\eqref{eqt: pde2_SJ_general_minplus}. Similarly, for any $j\in\mathcal{J}$, we have that
\begin{equation*}
    \Sanaly(\bx,t) - \Snum(\bx,t) = \Sanaly_j(\bx,t) - \Snum(\bx,t) \geq  \Sanaly_j(\bx,t) - \Snum_j(\bx,t) \geq -\epsilon,
\end{equation*}
where the first equality holds by definition of $\mathcal{J}$ and the first inequality holds since $\Snum$ satisfies~\eqref{eqt:alg_minplus_output}.
Therefore,~\eqref{eqt:prop_err_minplus_S} holds.

Now, assume $\Sanaly_j(\bx,t) > \Sanaly(\bx,t) + 2\epsilon$ holds for each index $j$ satisfying $\Sanaly_j(\bx,t) \neq \Sanaly(\bx,t)$. We prove $r\in \mathcal{J}$ by contradiction. Assume $r$ is not in $\mathcal{J}$. Then, we have $\Sanaly_r(\bx,t) \neq \Sanaly(\bx,t)$. However, from straightforward calculation, we also have
\begin{equation*}
\begin{split}
    \Sanaly_r(\bx,t) - \Sanaly(\bx,t) &\leq (\Sanaly_r(\bx,t) - \Snum_r(\bx,t)) + (\Snum_r(\bx,t) - \Snum(\bx,t)) + (\Snum(\bx,t)- \Sanaly(\bx,t))\\
    &\leq \epsilon + 0 + \epsilon = 2\epsilon,
\end{split}
\end{equation*}
which leads to a contradiction with our assumption.
Therefore, we have $r\in\mathcal{J}$, and hence~\eqref{eqt:prop_err_minplus_traj} holds by definition of $r$ and $\mathcal{J}$.
\qed

\end{document}